\documentclass[reqno,11pt]{amsart}
\usepackage[margin=1in]{geometry}
\usepackage[english]{babel}
\usepackage{amsmath,amsthm,amssymb,mathtools}
\usepackage{dsfont}
\usepackage{tikz-cd}
\usepackage{xcolor}
\usepackage{hyperref}
\hypersetup{colorlinks=true, urlcolor=blue, linkcolor=blue, citecolor=black}

\makeatletter
\let\markboth\@gobbletwo
\let\markright\@gobble
\makeatother

\numberwithin{equation}{section}
\theoremstyle{plain}
\newtheorem{thmctr}{}[section] %
\newtheorem{lemma}[thmctr]{Lemma}
\newtheorem{theorem}[thmctr]{Theorem}
\newtheorem{corollary}[thmctr]{Corollary}
\newtheorem{proposition}[thmctr]{Proposition}
\theoremstyle{definition}
\newtheorem{definition}[thmctr]{Definition}
\theoremstyle{remark}
\newtheorem{remark}[thmctr]{Remark}

\allowdisplaybreaks

\newcommand{\cA}{\mathcal{A}}
\newcommand{\B}{\mathcal{B}}
\newcommand{\C}{\mathcal{C}}
\newcommand{\D}{\mathcal{D}}
\newcommand{\E}{\mathcal{E}}
\newcommand{\cL}{\mathcal{L}}
\newcommand{\M}{\mathcal{M}}
\newcommand{\N}{\mathcal{N}}
\newcommand{\Z}{\mathcal{Z}}
\newcommand{\oB}{\overline{\mathcal{\B}}}
\newcommand{\Vect}{\mathsf{Vec}}

\newcommand{\bbA}{\mathbb{A}}
\newcommand{\bbL}{\mathbb{L}}
\newcommand{\bbX}{\mathbb{X}}
\newcommand{\fg}{\mathfrak{g}}
\newcommand{\kk}{\Bbbk}
\newcommand{\unit}{\mathds{1}}

\providecommand{\hcomp}{\mathbin{\circ_{\mathrm{h}}}}

\DeclareMathOperator{\rank}{rank}
\newcommand{\ev}{\mathrm{ev}}
\newcommand{\coev}{\mathrm{coev}}
\newcommand{\id}{\mathrm{id}}
\newcommand{\br}{\mathrm{br}}
\newcommand{\Rep}{\mathrm{Rep}}
\newcommand{\Hom}{\mathrm{Hom}}
\newcommand{\iHom}{\underline{\Hom}}
\newcommand{\ad}{\mathrm{ad}}
\newcommand{\loc}{\mathrm{loc}}
\newcommand{\rev}{\mathrm{rev}}
\newcommand{\op}{\mathrm{op}}
\newcommand{\inn}{\mathrm{inn}}
\newcommand{\Fun}{\mathrm{Fun}}
\newcommand{\Aut}{\mathrm{Aut}}
\newcommand{\Alg}{{\sf Alg}}
\newcommand{\Rex}{{\sf Rex}}
\newcommand{\ra}{{\sf ra}}
\newcommand{\Inv}{{\sf Inv}}
\newcommand{\FPdim}{\mathsf{d}^{\scriptscriptstyle\mathrm{FP}}}
\newcommand{\tr}{\triangleright}
\newcommand{\tl}{\triangleleft}

\newcommand{\BrPic}{\mathrm{BrPic}}
\newcommand{\uBrPic}{\underline{\BrPic}}
\newcommand{\Pic}{\mathrm{Pic}}
\newcommand{\uPic}{\underline{\Pic}}
\newcommand{\uAut}{\underline{\Aut}}
\newcommand{\EqBr}{\mathrm{EqBr}}
\newcommand{\uEqBr}{\underline{\smash{\EqBr}}}

\newcommand{\Forg}{\mathrm{U}}
\newcommand{\Induc}{\mathrm{I}}
\newcommand{\Schau}{\Theta}
\newcommand{\Kappa}{\mathrm{K}}

\theoremstyle{plain}
\newtheorem{introthm}{Theorem}

\title[Lagrangian algebras and Brauer--Picard groupoids]
{Lagrangian Algebras and Brauer--Picard Groupoids of Finite Tensor Categories}

\author{Harshit Yadav}
\address{Max-Planck Institute for Mathematics, Vivatsgasse 7, 53111 Bonn, Germany}
\email{harshit11y@gmail.com}

\author{Kenichi Shimizu}
\address{Department of Mathematical Sciences, Shibaura Institute of Technology,\linebreak
307 Fukasaku, Minuma-ku, Saitama-shi, Saitama 337-8570, Japan}
\email{kshimizu@shibaura-it.ac.jp}

\begin{document}
\begin{abstract}
Let $\C$ be a finite tensor category over an algebraically closed field. We
show that indecomposable exact $\C$-module categories correspond bijectively
to Lagrangian algebras in the Drinfeld center $\Z(\C)$, through the adjoint
algebra of a module category. Under this bijection, unimodular module
categories correspond to Frobenius algebras, and for pivotal $\C$, pivotal
unimodular ones correspond to symmetric Frobenius algebras. We then prove
that taking the Drinfeld center is a fully faithful pseudofunctor from the
Brauer--Picard 2-groupoid of finite tensor categories to the 2-groupoid of
braided finite tensor categories. In particular, we recover the
isomorphism $\BrPic(\C)\cong\Aut_{\otimes}^{\br}(\Z(\C))$ of Davydov and
Nikshych, and its Picard analogue for braided $\C$. We also classify
the finite tensor categories Morita equivalent to $\C$ by orbits of
Lagrangian algebras. Finally, for $\C=\Rep(u_q(\fg))$ with $u_q(\fg)$ a small
quantum group at a suitable root of unity, we show that there are infinitely many
$\Aut_{\otimes}(\C)$-orbits of indecomposable exact $\C$-module categories,
all with the same dual category.
\end{abstract}

\maketitle

\setcounter{tocdepth}{1}
\tableofcontents

\section{Introduction}
\label{sec:introduction}

Finite tensor categories are non-semisimple generalizations of fusion
categories, and they appear throughout representation theory,
low-dimensional topology and conformal field theory. As for many algebraic
structures, a finite tensor category $\C$ is studied through its
`representations', the $\C$-module categories, and among these the exact
module categories of Etingof and Ostrik \cite{etingof2004finite} play the
central role. An indecomposable exact $\C$-module category $\M$ has a dual
finite tensor category $\C^*_{\M}$, which is Morita equivalent to $\C$
\cite[\S7.12]{etingof2016tensor}, and the Drinfeld centers of $\C$ and
$\C^*_{\M}$ are braided equivalent \cite{schauenburg2001monoidal}. So the
center $\Z(\C)$ is a Morita invariant of $\C$, and one may ask how much of
the module theory of $\C$ it remembers.

For fusion categories over a field of characteristic zero, the center
determines this module theory completely, by three results. Davydov,
M\"uger, Nikshych and Ostrik showed that indecomposable $\C$-module
categories correspond bijectively to Lagrangian algebras in $\Z(\C)$
\cite[Proposition~4.8]{davydov2013witt}. Etingof, Nikshych and Ostrik proved
that two fusion categories are Morita equivalent if and only if their
centers are braided equivalent \cite[Theorem~3.1]{etingof2011weakly}.
M\"uger had observed earlier that Morita equivalent spherical fusion categories have braided equivalent centers \cite[Remark~3.18]{muger2003subfactorsI}. More
precisely, taking the center is a fully faithful 2-functor from the
Brauer--Picard 2-groupoid to the 2-groupoid of braided fusion categories
\cite[Theorem~1.1]{etingof2010fusion}. In particular, the Brauer--Picard
group $\BrPic(\C)$ of invertible $\C$-bimodule categories is isomorphic to
the group $\Aut_{\otimes}^{\br}(\Z(\C))$ of braided autoequivalences of the
center. This leads to the question that drives this paper.

\smallskip
\noindent\textbf{Question.} \emph{Do these three results hold for finite
tensor categories over an algebraically closed field $\kk$ that need not be
semisimple?}
\smallskip

Here the indecomposable semisimple module categories of the fusion case are replaced by
the indecomposable exact ones. Davydov and Nikshych studied the
Brauer--Picard group in this generality and stated the isomorphism
$\BrPic(\C)\cong\Aut_{\otimes}^{\br}(\Z(\C))$ for finite tensor categories
\cite[Theorem~4.1]{davydov2013picard}. The semisimple arguments do not carry
over directly, however: they use that certain algebras built from the data
are exact, which is automatic for fusion categories but not true in general
(Remark~\ref{rem:DN-gap}). For the Morita result
\cite[Theorem~8.12.3]{etingof2016tensor}, Ostrik and Utiralova found a way
around this \cite[Corollary~5.14]{ostrik2026non}, and their idea is the
starting point of our approach. A key ingredient, for them and for us, is
the theorem of Coulembier, Stroi\'nski and Zorman that an algebra in a
finite tensor category is exact if and only if it is a finite product of
simple algebras \cite[Theorem~7.1]{coulembier2025simple}. It is what makes many arguments of the fusion case available for finite tensor categories.

The main purpose of this paper is to answer the question affirmatively, and
to refine the answer by pivotal and unimodular structures and by the 2-cells
of the Brauer--Picard 2-groupoid. We also apply the results to a small quantum group $u_q(\fg)$: the category $\Rep(u_q(\fg))$ has infinitely many indecomposable exact module categories with the same dual category, no two of which are related by a tensor autoequivalence of $\Rep(u_q(\fg))$.

\subsection{Lagrangian algebras and module categories}
A Lagrangian algebra in $\Z(\C)$ is a simple commutative algebra $\bbL$ with
$\FPdim(\bbL)=\FPdim(\C)$, where $\FPdim$ means the Frobenius--Perron dimension (Definition~\ref{def:lagrangian-algebra}). A
basic example is the adjoint algebra $\Induc_{\C}(\unit)$, where $\Induc_{\C}$
is right adjoint to the forgetful functor $\Forg_\C:\Z(\C)\to\C$. More generally, an
indecomposable exact left $\C$-module category $\M$ has an adjoint algebra
$\bbA_{\M}\in\Z(\C)$, obtained by transporting the adjoint algebra of the
dual category $\C^*_{\M}$ along the Schauenburg equivalence
$\Z(\C)\simeq\Z(\C^*_{\M})$ (\S\ref{subsec:adjoint-algebra}).
Our first main result says that every Lagrangian algebra arises in this
way, from a unique module category.

\begin{introthm}\label{thm:intro-A}
  Let $\C$ be a finite tensor category.
  \begin{enumerate}
  \item The assignment $\M\mapsto\bbA_{\M}$ induces a bijection between
    equivalence classes of indecomposable exact $\C$-module categories
    and isomorphism classes of Lagrangian algebras in $\Z(\C)$.
  \item The algebra $\bbA_{\M}$ is Frobenius if and only if $\M$ is
    unimodular, that is, $\C^*_{\M}$ is unimodular.
  \item If $\C$ is pivotal and $\Z(\C)$ carries the induced pivotal
    structure, then $\bbA_{\M}$ is symmetric Frobenius if and only if $\M$
    is pivotal and unimodular.
  \end{enumerate}
\end{introthm}

Part~(a) of Theorem~\ref{thm:intro-A} is the main new contribution of this paper. It is Theorem~\ref{thm:classification-Lagrangian}. It generalizes
\cite[Proposition~4.8]{davydov2013witt} from fusion to finite tensor
categories, with a different proof, and it also matches morphisms: every
algebra isomorphism $\bbA_{\M}\to\bbA_{\N}$ is induced by a module
equivalence $\M\to\N$, unique up to isomorphism
(Lemma~\ref{lem:isomorphisms-between-Lagrangian}). Part~(b),
(Proposition~\ref{prop:Lagrangian-uni}) together with~(a), shows that $\C$
is Morita equivalent to a unimodular finite tensor category if and only if
$\Z(\C)$ contains a Frobenius Lagrangian algebra. Proposition~\ref{prop:Lagrangian-uni} is deduced from \cite[Theorem~A]{yadav2023unimodular}. Part~(c) is Proposition~\ref{prop:Lagrangian-piv}, whose `if' direction is \cite[Corollary~6.15(ii)]{yadav2026Frobenius}.

\subsection{The Brauer--Picard 2-groupoid}
Our second main result says that the center also remembers invertible
bimodule categories, together with their autoequivalences. The
Brauer--Picard 2-groupoid $\uBrPic$ has finite tensor categories as 0-cells,
invertible bimodule categories as 1-cells, and isomorphism classes of
bimodule equivalences as 2-cells \cite[\S4]{etingof2010fusion},
\cite[\S2G]{davydov2013picard}. The 2-groupoid $\uEqBr$ has braided finite
tensor categories as 0-cells, braided tensor equivalences as 1-cells, and
monoidal natural isomorphisms as 2-cells.

\begin{introthm}\label{thm:intro-B}
  The assignment $\C\mapsto\Z(\C)$ extends to a fully faithful pseudofunctor
  $\Z:\uBrPic\to\uEqBr$. That is, for finite tensor categories $\C$ and
  $\D$, the center functor from the groupoid of invertible $\C$-$\D$-bimodule
  categories to the groupoid of braided tensor equivalences
  $\Z(\C)\to\Z(\D)$ is an equivalence. In particular,
  $\BrPic(\C)\cong\Aut_{\otimes}^{\br}(\Z(\C))$.
\end{introthm}

This is Theorem~\ref{thm:BrPic} together with~\eqref{eq:BrPic-center-iso}.
It generalizes \cite[Theorem~1.1]{etingof2010fusion} from fusion to finite
tensor categories. Kong and Zheng proved a related full faithfulness in
characteristic zero, for indecomposable multi-fusion categories with
equivalence classes of nonzero semisimple bimodule categories as morphisms
\cite[Theorem~3.3.7]{kong2018center}. For finite tensor categories, Davydov
and Nikshych stated the group isomorphism
\cite[Theorem~4.1]{davydov2013picard} and its upgrade to a monoidal
equivalence of groupoids \cite[Remark~4.2]{davydov2013picard}.
Theorem~\ref{thm:intro-B} recovers both,
and extends them to a fully faithful
pseudofunctor on the whole Brauer--Picard 2-groupoid.
The statement on 2-cells says that the isomorphism classes of bimodule autoequivalences of an invertible
bimodule category $\M$ correspond to the monoidal automorphisms of the
braided equivalence attached to it.

Theorem~\ref{thm:intro-B} turns the classification of finite tensor
categories Morita equivalent to $\C$ into a question about Lagrangian
algebras. Let $\bbX_{\C}$ be the set of equivalence classes of indecomposable
exact left $\C$-module categories. The group $\BrPic(\C)$ acts on $\bbX_{\C}$
by $[\mathcal P]\cdot[\M]=[\mathcal P\boxtimes_{\C}\M]$, and
$\Aut_{\otimes}^{\br}(\Z(\C))$ acts on the set $\mathrm{Lag}(\Z(\C))$ of
isomorphism classes of Lagrangian algebras. The following is
Theorem~\ref{thm:MME} and Corollary~\ref{cor:BrPic-stabilizer-Out}.

\begin{introthm}\label{thm:intro-C}
  Let $\C$ be a finite tensor category, and let $\M$ and $\N$ be
  indecomposable exact left $\C$-module categories.
  \begin{enumerate}
      \item $\C^*_{\M}\simeq\C^*_{\N}$ as tensor categories if and only if $[\M]$ and $[\N]$ lie in the same $\BrPic(\C)$-orbit.
      \item The assignment
      $\M\mapsto\C^*_{\M}$ induces a bijection from the set of
      $\BrPic(\C)$-orbits in $\bbX_{\C}$ to the set of tensor equivalence
      classes of finite tensor categories Morita equivalent to $\C$.
      \item  Moreover, the assignment $\M\mapsto\bbA_{\M}$ induces a bijection from the set of $\BrPic(\C)$-orbits to the
      set of $\Aut_{\otimes}^{\br}(\Z(\C))$-orbits in $\mathrm{Lag}(\Z(\C))$.
      \item   The stabilizer of $[\M]$ is isomorphic to $\mathrm{Out}_{\otimes}(\C^*_{\M})$, the group of tensor autoequivalences of $\C^*_{\M}$ modulo conjugation by invertible objects.
  \end{enumerate}
\end{introthm}

In the fusion case, the description of Morita equivalent categories by
orbits of Lagrangian algebras can be deduced from
\cite[Theorem~3.1]{etingof2011weakly} and the 2-equivalence between fusion
categories and quantum Manin pairs \cite[Proposition~4.4]{davydov2013witt}.
Theorem~\ref{thm:intro-C} generalizes it to finite tensor categories.

\subsection{Picard groupoids and small quantum groups}
Now let $\B$ be a braided finite tensor category with braiding $c$. Through the braiding, a
left $\B$-module category $\M$ becomes a $\B$-bimodule category $B(\M)$, and
the Picard groupoid $\uPic(\B)$ consists of those $\M$ for which $B(\M)$ is
invertible (\S\ref{subsec:picard-groupoid}). The center contains two
copies of $\B$, embedded by $G^+:\B\to\Z(\B)$, $X\mapsto(X,c_{X,-})$, and
$G^-:\oB\to\Z(\B)$, $X\mapsto(X,c_{-,X}^{-1})$. Our next result
describes $\uPic(\B)$ by braided autoequivalences of $\Z(\B)$ that are
trivialized on the second copy.

\begin{introthm}\label{thm:intro-D}
  Let $\B$ be a braided finite tensor category. The monoidal groupoid
  $\uPic(\B)$ is equivalent to the monoidal groupoid
  $\uAut^{\br}(\Z(\B);\oB)$ of braided autoequivalences of $\Z(\B)$ equipped
  with a trivialization on $G^-(\oB)$. In particular,
  $\Pic(\B)\cong\Aut_{\otimes}^{\br}(\Z(\B);\oB)$.
\end{introthm}

This is Theorem~\ref{thm:Pic} and
Corollary~\ref{cor:Pic-center-autoequivalences}. Davydov and Nikshych also had the group isomorphism \cite[Theorem~4.3]{davydov2013picard} and they deduced it
by restricting the isomorphism of \cite[Theorem~4.1]{davydov2013picard}.
Theorem~\ref{thm:intro-D} recovers it, with the convention for $B(\M)$ of \S\ref{subsec:picard-groupoid}, and lifts it to monoidal groupoids.
For $\M\in\uPic(\B)$, the braided autoequivalence of $\Z(\B)$
attached to $\M$ preserves $G^+(\B)$, and its restriction there is the
boundary autoequivalence $\partial_{\M}$ of $\B$ defined by Davydov and Nikshych
(Lemma~\ref{lem:boundary-map}). The resulting boundary map
$\partial:\Pic(\B)\to\Aut_{\otimes}^{\br}(\B)$ is bijective when $\B$ is
factorizable (Proposition~\ref{prop:boundary-bijective}; compare
\cite[Corollary~4.6]{davydov2013picard}).

The group $\Aut_{\otimes}(\B)$ acts on $\bbX_{\B}$ by pulling back the action
along tensor autoequivalences, and we call module categories in the same orbit
\emph{pullback equivalent}. Pullback equivalent module categories have
equivalent dual categories, so each pullback equivalence class lies in a single
class of Theorem~\ref{thm:intro-C}.

Braided autoequivalences preserve $\Pic(\B)$, and $\partial$ turns pullback
into conjugation: $\partial_{{}_F\M}\cong F\circ\partial_{\M}\circ F^{-1}$ for
$F\in\Aut_{\otimes}^{\br}(\B)$, where ${}_F\M$ is the pullback of $\M$ along
$F$ (Lemma~\ref{lem:boundary-map}). Hence, if $\partial$ is surjective, for
instance if $\B$ is factorizable, it maps the $\Aut_{\otimes}^{\br}(\B)$-orbits
in $\Pic(\B)$ onto the conjugacy classes of $\Aut_{\otimes}^{\br}(\B)$, so
there are infinitely many such orbits as soon as there are infinitely many
conjugacy classes (Proposition~\ref{prop:inf-orbit-criteria}). For small
quantum groups, all of $\Pic(\B)$ lies in the class of
Theorem~\ref{thm:intro-C} that contains the regular module category.

\begin{introthm}\label{thm:intro-E}
  Let $\kk$ have characteristic $0$, let $\fg$ be a finite-dimensional
  simple Lie algebra, and let $q\in\kk$ be a primitive $l$-th root of unity
  for an odd integer $l>1$ with $\gcd(l,\det A_{\fg})=1$, where $A_{\fg}$ is
  the Cartan matrix, and with $\gcd(l,3)=1$ in type $G_2$. Then
  $\B=\Rep(u_q(\fg))$ has infinitely many $\Aut_{\otimes}(\B)$-orbits of
  indecomposable exact left module categories, all with the same dual
  category as the regular module category.
\end{introthm}

This is Theorem~\ref{thm:uqg-Aut-orbits-infinite} together with Proposition~\ref{prop:uqg-regular-orbit}, with the braiding on
$\B=\Rep(u_q(\fg))$ of \cite[\S2.3]{davydov2018autoequivalences}.
There every tensor autoequivalence of $\B$ is braided
\cite[Theorem~5.4]{davydov2018autoequivalences}, and the pullback classes of
module categories with the same dual category as the regular one correspond
to the conjugacy classes of $\Aut_{\otimes}(\B)$
(Proposition~\ref{prop:uqg-regular-orbit}). By
\cite[Lemma~5.5]{davydov2018autoequivalences}, this group contains the
adjoint group $G^{\ad}(\kk)$ as a normal subgroup of finite index, whose
semisimple conjugacy classes are already infinite in number.

\subsection{The method}
Everything rests on Theorem~\ref{thm:intro-A}(a).
For fusion categories, one recovers a module category from a
Lagrangian algebra $\bbL$ by decomposing the algebra $\Forg_{\C}(\bbL)$ in
$\C$ into indecomposable summands and taking modules over one of them. This
is how the inverse map is constructed in the proof of
\cite[Theorem~4.1]{davydov2013picard}. Without semisimplicity, the algebra
$\Forg_{\C}(\bbL)$ need not be exact. Already for $\bbL=\Induc_{\C}(\unit)$,
the $\C$-module category of $\Forg_{\C}(\bbL)$-modules in $\C$ is
equivalent to $\C\boxtimes\C^{\op}$, which is exact only if $\C$ is semisimple
(Remark~\ref{rem:DN-gap}).

We prove surjectivity without decomposing $\Forg_{\C}(\bbL)$. By a
theorem of Ostrik and Utiralova \cite[Theorem~5.13]{ostrik2026non}, which
rests on the characterization of exact algebras as finite products of simple
algebras \cite[Theorem~7.1]{coulembier2025simple}, every simple commutative
algebra in $\Z(\C)$ is a subalgebra of an adjoint algebra $\bbA_{\M}$. A
Lagrangian algebra has the same FP dimension $\FPdim(\C)$ as $\bbA_{\M}$, so
this inclusion is an isomorphism. Ostrik and Utiralova prove
\cite[Corollary~5.14]{ostrik2026non} by the same comparison of FP
dimensions. For injectivity we use the equivalence
$\C_{\Forg_{\C}(\bbA_{\M})}\simeq\M\boxtimes\M^{\op}$ of left $\C$-module
categories. Then, the $\C$-module subcategory generated by any simple object of
$\M\boxtimes\M^{\op}$ is equivalent to $\M$, so $\bbA_{\M}$ determines $\M$.

For Theorem~\ref{thm:intro-B}, a braided equivalence $T:\Z(\C)\to\Z(\D)$
gives the Lagrangian algebra $T^{-1}(\bbA_{\D})$, hence a module category
$\M$ with $\bbA_{\M}\cong T^{-1}(\bbA_{\D})$, whose right $\D$-action comes
from the tensor equivalences
$\D\simeq\Z(\D)_{\bbA_{\D}}\simeq\Z(\C)_{\bbA_{\M}}\simeq\C^*_{\M}$
of~\eqref{eq:def-psi-right-action}, the middle one induced by $T^{-1}$. This inverts $\Z$ on 1-cells up to
isomorphism. On 2-cells, the bimodule autoequivalences of $\M$ are given, up
to isomorphism, by invertible objects $\boldsymbol\beta$ of $\Z(\C)$
(Lemma~\ref{lem:Aut-M-invertible}), and $\Z$ sends the one given by
$\boldsymbol\beta$ to the monodromy $\Omega_{\boldsymbol\beta}$ whiskered
with the braided equivalence attached to $\M$
(Lemma~\ref{lem:Z-of-invertible-object}). Every monoidal automorphism of the
identity functor of $\Z(\C)$ is the monodromy of a unique isomorphism class
of invertible objects \cite[Theorem~5.3]{shimizu2017ribbon}.

Theorem~\ref{thm:intro-C} follows from Theorem~\ref{thm:intro-B}: composing
with invertible bimodule categories identifies the fibers of
$\M\mapsto\C^*_{\M}$ with the $\BrPic(\C)$-orbits, and $\M\mapsto\bbA_{\M}$
is equivariant along $\BrPic(\C)\cong\Aut_{\otimes}^{\br}(\Z(\C))$
(Lemma~\ref{lem:adj-alg-and-partial-M}). Theorem~\ref{thm:intro-D} follows by
matching one-sided structures on $B(\M)$ with trivializations on $G^-(\oB)$
(Lemma~\ref{lem:Pic-trivialization}).

\subsection{Organization}
\S\ref{sec:preliminaries} collects preliminaries on finite tensor
categories, induced autoequivalences of centers, and Lagrangian algebras.
\S\ref{sec:module-categories} treats module categories: relative
Deligne products, exactness, the invertibility criterion
(Proposition~\ref{prop:invertibility-criterion}), and the Schauenburg
equivalence. \S\ref{sec:lagrangian-algebras} proves
Theorem~\ref{thm:intro-A}, \S\ref{sec:BrPic-2-groupoid} proves
Theorems~\ref{thm:intro-B} and~\ref{thm:intro-C}, and
\S\ref{sec:picard-groupoid-orbits} proves Theorems~\ref{thm:intro-D}
and~\ref{thm:intro-E}.

\subsection{Acknowledgements}
HY thanks the Max Planck Institute for Mathematics in Bonn for its
hospitality and financial support. KS is supported by JSPS KAKENHI Grant
Number JP24K06676.

\subsection{Use of AI tools.}
AI assistance was used for ideas in \S\ref{sec:picard-groupoid-orbits}, in particular the idea for Theorem~\ref{thm:intro-E} was suggested by ChatGPT. The rest of the mathematics was developed by the authors without AI assistance. AI tools were then used to double-check proofs and citations, to suggest shorter versions of some arguments, and to edit the prose. The authors checked every such suggestion and take full responsibility for the content.

\section{Preliminaries}
\label{sec:preliminaries}
Throughout this article, the base field $\kk$ is algebraically closed. We usually use $\C,\D$ for finite tensor categories and $\B,\B'$ when a
braiding is part of the data.
Given a functor $F$, we write $F^{\ra}$ for a right adjoint of $F$, if it exists. Given a category $\C$, we write $\C^{\op}$ for the opposite category of $\C$.
For a natural transformation $\alpha$, its $X$-component is denoted by $\alpha(X)$. We write $\circ$ and $\hcomp$ for vertical and horizontal composition of natural transformations, respectively.
Thus, for natural transformations $\alpha : G \Rightarrow G'$ and $\beta : F \Rightarrow F'$ with $GF$ defined, we have
\begin{equation*}
  (\alpha \hcomp \beta)(X) = \alpha(F'(X)) \circ G(\beta(X))
  : G F(X) \to G' F'(X)
\end{equation*}
for an object $X$ of the source of $F$. When the distinction is clear from the context, we omit the subscript and write $\hcomp$ simply as $\circ$.

\subsection{Finite tensor categories}

We use the standard terminology of finite abelian categories, finite tensor categories,
and finite module categories. We refer the reader to \cite{etingof2016tensor} for the basic theory of these.
A \emph{left dual} of an object $X$ of a monoidal category is an object $X^*$ together with morphisms $\ev_X:X^*\otimes X\to\unit$ and $\coev_X:\unit\to X\otimes X^*$ satisfying the zigzag identities \cite[Definition~2.10.1]{etingof2016tensor}. Right duals ${}^*X$ are defined similarly.
Given a finite tensor category $\C$, we write $\C^{\rev}$ for the same category with the reversed tensor product $X\otimes^{\rev} Y=Y\otimes X$.

A braided tensor category is a finite tensor category $\B$ equipped with a braiding map $c_{X,Y} : X\otimes Y \xrightarrow{\sim} Y \otimes X$ for all $X, Y \in \B$. If $\B$ is braided, then $\oB$ denotes the braided category with the same underlying tensor category and the inverse braiding $\overline{c}_{X,Y}=c_{Y,X}^{-1}$.

A \emph{tensor functor} is a $\kk$-linear strong monoidal functor between finite tensor categories that is exact and faithful \cite[Definition~4.2.5]{etingof2016tensor}. We write $\Aut_{\otimes}(\C)$ for the group of isomorphism classes of tensor autoequivalences of $\C$, and $\Aut_{\otimes}^{\br}(\B)$ for the subgroup of braided ones when $\B$ is braided.

We call a tensor functor $F:\C\to \D$ {\em surjective} if every object
$Y\in \D$ is isomorphic to a subquotient of $F(X)$ for some $X\in \C$.
This is equivalent to the condition of \cite[Definition~1.8.3]{etingof2016tensor} that every simple object $S\in\D$ be a subquotient of some $F(X)$. Indeed, $F$ then maps projective objects to projective objects \cite[Theorem~6.1.16]{etingof2016tensor}. If $P$ is a projective cover of $\unit$, then $F(P)\otimes S$ is projective \cite[Proposition~4.2.12]{etingof2016tensor} and maps onto $S$, so the projective cover of $S$ is a direct summand of $F(P)\otimes S$. It is therefore a projective subquotient of $F(P\otimes X)$, hence a direct summand of it, since projective objects are injective \cite[Proposition~6.1.3]{etingof2016tensor}. As every object of $\D$ is a quotient of a direct sum of projective covers of simple objects, the claim follows.

\subsubsection{Drinfeld centers}
\label{subsubsec:drinfeld-centers}

The Drinfeld center $\Z(\C)$ of a finite tensor category $\C$ is the category whose objects are pairs $(V,\sigma)$, where $V\in \C$ and $\sigma:V\otimes \id_{\C}\to \id_{\C}\otimes V$
is a natural isomorphism satisfying the half-hexagon axiom. This fixes the direction of half-braidings. This direction is opposite to \cite{etingof2016tensor} and the same as \cite{schauenburg2001monoidal}. Since half-braidings are invertible, $(V,\sigma)\mapsto(V,\sigma^{-1})$ identifies $\Z(\C)$ with the center of \cite{etingof2016tensor} as a tensor category over $\C$, with the braiding replaced by its reverse. The Drinfeld center is a typical example of a braided tensor category: for $\B=\Z(\C)$, the braiding is $c_{(X,\sigma),(Y,\rho)} = \sigma_Y$.

We write $\Forg_{\C}:\Z(\C)\to \C$
for the forgetful functor. Let $\Induc_{\C}:\C\to \Z(\C)$ be its right adjoint.
The functor $\Induc_{\C}$ is lax monoidal as the right adjoint of the strong
monoidal functor $\Forg_{\C}$. The forgetful functor $\Forg_{\C}$ is
surjective \cite[Corollary~7.13.11]{etingof2016tensor}.

\subsubsection{Nondegeneracy}
\label{subsubsec:nondegeneracy}
Given a braided finite tensor category $\B$, we have two braided tensor functors
\[
G^+:\B \to \Z(\B), \quad X \mapsto (X,c_{X,-})
\quad \text{and} \quad
G^-: \oB \to \Z(\B), \quad X \mapsto (X,c_{-,X}^{-1}).
\]
Combining these we get a braided tensor functor
\[
G:=G_\B : \B\boxtimes \oB \to \Z(\B), \qquad X\boxtimes Y \mapsto \left( X\otimes Y, \sigma = (c_{X,-}\otimes \id_Y)\circ (\id_X\otimes c_{-,Y}^{-1}) \right)
\]
We say that $\B$ is \emph{factorizable} if $G$ is an equivalence. On the other hand, $\B$ is called \emph{nondegenerate} if $\Z_{(2)}(\B)\simeq \Vect$, where $\Z_{(2)}(\B)$ is the M\"uger center of $\B$ given by
\[
\Z_{(2)}(\B) = \{ X \in \B \ \mid \ c_{Y,X}\circ c_{X,Y} = \id_{X\otimes Y}, \ \forall \ Y\in\B \}.
\]
By \cite[Theorem~1.1]{shimizu2019non}, factorizability and nondegeneracy are equivalent.

\subsection{Induced autoequivalences of Drinfeld centers}
\label{subsec:induced-autoequivalences}
Given a tensor equivalence $F:\C\to \D$ with quasi-inverse $F^{-1}$, we write
$\widetilde{F}:\Z(\C)\to \Z(\D)$ for the induced braided tensor equivalence.
More precisely, $\widetilde{F}$ sends $(V,\sigma)\in\Z(\C)$ to $F(V)$ with the half-braiding
\[
  F(V)\otimes Y\cong F\bigl(V\otimes F^{-1}(Y)\bigr)
  \xrightarrow{\,F(\sigma_{F^{-1}(Y)})\,}
  F\bigl(F^{-1}(Y)\otimes V\bigr)\cong Y\otimes F(V)
  \qquad(Y\in\D).
\]
Up to canonical isomorphism, $\widetilde{\id}=\id$ and
$\widetilde{G\circ F}\cong\widetilde G\circ\widetilde F$. Since
$\Forg_{\D}\circ\widetilde F=F\circ\Forg_{\C}$, passing to right adjoints
(that is, taking the mate of this equality) gives the canonical isomorphism
\begin{equation}
  \label{eq:F-tilde-and-induction}
  \widetilde{F}\circ \Induc_{\C}
  \cong
  \Induc_{\D}\circ F
\end{equation}
of lax monoidal functors.

We now determine when an induced autoequivalence is isomorphic to the identity.
We denote by $\Inv(\C)$ the group of isomorphism classes of invertible objects
of $\C$. For $\beta\in \Inv(\C)$, let
\begin{equation}
  \label{eq:def-inner-auto}
  \inn_{\beta}:\C\to \C,
  \qquad
  X\mapsto \beta\otimes X\otimes \beta^*,
\end{equation}
be the inner tensor autoequivalence implemented by $\beta$, with tensor
structure induced by $\ev_\beta$ and unit constraint $\coev_\beta$.

We also need the following categorified Rosenberg--Zelinsky exact sequence \cite[Theorem~5.4]{shimizu2023pivotal} (see \cite[Theorem~1.3]{galindo2017tensor} for the fusion category case). For a finite tensor category $\C$, the sequence
\begin{equation}
  \label{eq:RZ-exact-sequence}
  \Inv(\Z(\C))
  \to \Inv(\C)
  \to \Aut_{\otimes}(\C)
  \to \Aut_{\otimes}^{\mathrm{br}}(\Z(\C))
\end{equation}
is exact. The first map forgets the half-braiding, the second sends $\beta$ to
$\inn_{\beta}$, and the last sends a tensor autoequivalence $F$ to
$\widetilde{F}$.

\begin{lemma}
  \label{lem:tilde-isomorphic}
  Let $F,G:\C\to\D$ be tensor equivalences between finite tensor categories. Then
  $\widetilde{F}\cong\widetilde{G}$ as tensor functors if and only if
  $G\cong \inn_{\beta}\circ F$ as tensor functors for some invertible object
  $\beta\in \D$.
\end{lemma}
\begin{proof}
  Set $T=G\circ F^{-1}:\D\to\D$. Then
  $\widetilde{F}\cong \widetilde{G}$ if and only if $\widetilde{T}$ is
  isomorphic to the identity on $\Z(\D)$. By the exactness of
  \eqref{eq:RZ-exact-sequence}, this holds if and only if
  $T\cong \inn_{\beta}$ for some invertible object $\beta\in\D$.
\end{proof}

Let $\D$ be a finite tensor category, and let $\beta \in \D$ be invertible. By the exactness of \eqref{eq:RZ-exact-sequence} for $\D$, $\widetilde{\inn_\beta} \cong \id$. We now make this isomorphism explicit. For $\mathbf{X} = (X, \sigma_X) \in \Z(\D)$, put
\begin{equation*}
    \omega_{\beta}(\mathbf{X})
    = (\sigma_X(\beta) \otimes \id_{\beta^*}) \circ (\id_X \otimes \coev_{\beta}) : X \to \beta \otimes X \otimes \beta^*
\end{equation*}
These morphisms form an isomorphism $\omega_{\beta} : \id_{\Z(\D)} \to \widetilde{\inn_{\beta}}$
of tensor functors \cite[Lemma~4.3]{shimizu2023pivotal}, where it is
denoted $\chi^{\beta}$ (equation~(4.8) there).

\subsection{Algebras in tensor categories}
For an algebra $A$ in a finite tensor category $\C$, we write $\C_A$ for the
category of right $A$-modules in $\C$. It is a finite left $\C$-module category.

Let $\B$ be a braided finite tensor category and $A$ a commutative algebra in
$\B$. Then $G^-(A)=(A,c_{-,A}^{-1})$ is a commutative algebra in $\Z(\B)$, so
$\B_A$ is a monoidal category under $\otimes_A$: a right $A$-module $M$ is an
$A$-bimodule with left action
\[ A\otimes M\xrightarrow{c_{M,A}^{-1}}M\otimes A\to M . \]
With our convention for the Drinfeld center, this is the structure for which the free module functor $\B\to\B_A$ is central. For $A=F^{\ra}(\unit)$ with $F$ central, see \cite[Proposition~6.5 and Corollary~6.7]{shimizu2024commutative}. This category contains the following subcategory of \emph{local modules}
\[
\B_A^\loc = \{ (M,a)\in\B_A \ \mid \ a \circ c_{A,M} \circ c_{M,A} = a \}.
\]

We now recall the exactness properties of algebras that replace semisimplicity
in the later arguments. An algebra $A$ in a finite tensor category $\C$ is called \emph{exact} if for any projective $P\in \C$ and any object $M\in\C_A$, the object $P\otimes M\in\C_A$ is projective \cite[Definitions~3.1 and~3.18]{etingof2004finite}. An algebra
is called \emph{indecomposable} if it is nonzero and not isomorphic to a direct product of
two nonzero algebras.

A nonzero algebra is called \emph{simple} if it does not admit any nontrivial
two-sided ideals. The next theorem, which is \cite[Theorem~7.1]{coulembier2025simple}, records the characterization of exact algebras that we use repeatedly.

\begin{theorem}\label{thm:CSZ}
  Let $A$ be an algebra in a finite tensor category $\C$. Then, $A$ is exact if and only if $A$ is a finite direct product of simple algebras. In particular, $A$ is indecomposable exact if and only if it is simple. \qed
\end{theorem}

\subsection{Frobenius--Perron dimensions}

We refer the reader to \cite[\S4]{etingof2016tensor} for the basic theory of
Frobenius--Perron (FP) dimension. We write $\FPdim_{\C}(X)$ for the FP dimension of an object $X$ of a finite tensor category $\C$.

\begin{remark}\label{rem:FPdim}
Let $X\xrightarrow{i} Y$ be a monomorphism in $\C$ with $\FPdim_{\C}(X) = \FPdim_{\C}(Y)$. Then $X\cong Y$. Indeed, FP dimension is additive on short exact sequences and positive on non-zero objects. So equality of FP dimensions implies that the quotient $Y/X$ is zero and $i$ is an isomorphism.
\end{remark}

\begin{definition}\label{def:lagrangian-algebra}
A simple commutative algebra $A$ in a nondegenerate braided finite
tensor category $\B$ is called {\em Lagrangian} if
\begin{equation}\label{eq:lagrangian-fpdim}
  \FPdim_{\B}(A)^2 = \FPdim(\B).
\end{equation}
\end{definition}

A Lagrangian algebra is exact by Theorem~\ref{thm:CSZ} and haploid, that is, $\Hom_{\B}(\unit,A)\cong\kk$ \cite[Lemma~5.9]{shimizu2024commutative}.
For a finite tensor category
$\C$, the center $\Z(\C)$ is factorizable
\cite[Proposition~8.6.3]{etingof2016tensor}, hence nondegenerate by
\cite[Theorem~1.1]{shimizu2019non}, and
$\FPdim(\Z(\C))=\FPdim(\C)^2$ \cite[Theorem~7.16.6]{etingof2016tensor}. Thus a Lagrangian algebra in $\Z(\C)$ is a
simple commutative algebra $A$ with $\FPdim_{\Z(\C)}(A)=\FPdim(\C)$.

\subsection{Central tensor functors}

A tensor functor $F:\B\to \C$ is called \emph{central} if $\B$ is braided and $F$ admits a lift to the Drinfeld center $\Z(\C)$, meaning that there is a braided tensor functor $F':\B\to \Z(\C)$ such that $F=\Forg_{\C}\circ F'$.

Central structures are rigid along surjective functors: a monoidal
isomorphism between the underlying functors automatically respects the
half-braidings.

\begin{lemma}\label{lem:braided-lift-unique}
  Let $F_1,F_2:\B\to\C$ be central tensor functors with lifts
  $F_1',F_2':\B\to\Z(\C)$, and assume that $F_1$ is surjective. Then
  every monoidal natural isomorphism $\theta:F_1\Rightarrow F_2$ is a monoidal
  natural isomorphism $F_1'\Rightarrow F_2'$; that is, each
  $\theta_X$ is a morphism in $\Z(\C)$.
\end{lemma}
\begin{proof}
  Write $\sigma^X$ and $\tau^X$ for the half-braidings of $F_1'(X)$ and
  $F_2'(X)$. We must show that
  \[
    (\id_V\otimes\theta_X)\circ\sigma^X_V=\tau^X_V\circ(\theta_X\otimes\id_V)
    \qquad(V\in\C).
  \]
  Both sides are natural in $V$, and $\otimes$ is exact in each variable, so
  the identity passes from $V$ to its subobjects and quotients. Since $F_1$ is
  surjective, it suffices to take $V=F_1(Y)$. In that case, since $F_1'$ and
  $F_2'$ are braided, $\sigma^X_{F_1(Y)}$ and $\tau^X_{F_2(Y)}$
  are $F_1(c_{X,Y})$ and $F_2(c_{X,Y})$ conjugated by the tensor constraints,
  and the identity follows from the naturality of $\tau^X$ at $\theta_Y$
  and the naturality and monoidality of $\theta$.
\end{proof}
Lemma~\ref{lem:braided-lift-unique} also applies when $F_1$ is only isomorphic to $\Forg_{\C}\circ F_1'$ as
a tensor functor: transporting the half-braidings along the isomorphism gives a
lift of $F_1$ itself.

Suppose that $F_i:\B \to \D_i$ for $i=1,2$ are surjective central tensor functors and set $A_i = F_i^\ra(\unit)$. The right adjoint $F_i^\ra$ is lax monoidal, so $A_i$ is an algebra. It is commutative, exact and haploid \cite[Corollary~6.7]{shimizu2024commutative}, and $F_i^\ra$ induces tensor equivalences
\[
K_1: \D_1 \simeq \B_{A_1} \quad \text{and} \quad K_2: \D_2 \simeq \B_{A_2}.
\]
Here $K_i(Y)=F_i^\ra(Y)$ carries its canonical right $A_i$-action, and
$K_i\circ F_i$ is isomorphic to the free module functor as tensor functors
\cite[Proposition~6.6(c) and Corollary~6.7]{shimizu2024commutative}.
We also fix a tensor equivalence $G: \D_1 \xrightarrow{\sim} \D_2$ and an isomorphism
\[ \tilde{\chi} : F_2^\ra \circ G \xRightarrow{\sim} F_1^\ra \]
of lax monoidal functors. Then $f = \tilde{\chi}_{\unit}\circ F_2^{\ra}(G_0):A_2\to A_1$,
where $G_0:\unit\to G(\unit)$ is the unit constraint of $G$, is an
algebra isomorphism, and restriction of scalars along $f$ is a tensor
equivalence $f^\star:\B_{A_1}\to\B_{A_2}$.

\begin{lemma}\label{lem:innF-chi-iso}
Under the setting above, the following diagram commutes up to a monoidal
natural isomorphism:
\[\begin{tikzcd}
	{\D_1} && {\B_{A_1}} \\
	{\D_2} && {\B_{A_2}}
	\arrow["{K_1}", from=1-1, to=1-3]
	\arrow["G"', from=1-1, to=2-1]
	\arrow["{f^{\star}}", from=1-3, to=2-3]
	\arrow["{K_2}"', from=2-1, to=2-3]
\end{tikzcd}\]
\end{lemma}
\begin{proof}
  By construction, $K_i(Y)$ is $F_i^\ra(Y)$ equipped with its canonical
  right $A_i$-action.  Since $\tilde{\chi}$ is monoidal, its component
  $\tilde{\chi}_Y:F_2^\ra(G(Y))\to F_1^\ra(Y)$ intertwines the canonical
  $A_2$-action on the source with the restriction, along
  $f$, of the canonical $A_1$-action on the target.
  Thus $\tilde{\chi}_Y$ defines a natural isomorphism
  $K_2\circ G\xRightarrow{\sim}f^\star\circ K_1$. It is monoidal, since the
  tensor constraints of $K_1$ and $K_2$ are induced by the lax monoidal
  structures of $F_1^\ra$ and $F_2^\ra$, which $\tilde\chi$ intertwines.
\end{proof}

\section{Module categories}
\label{sec:module-categories}
This section collects the facts on module and bimodule categories used later: relative Deligne products, exact module categories, invertible bimodule categories, and Schauenburg equivalence.

\subsection{Module and bimodule categories}
\label{ssec:module-bimodule-center-prelim}

A finite left module category over $\C$ is a finite abelian category $\M$ equipped
with a left action $\triangleright:\C\times \M\to \M$ that is $\kk$-linear and
right exact in each variable.
Finite right module categories are defined similarly. Module functors are
assumed to be linear and right exact. A left $\C$-module functor
$F:\M\to\N$ comes with a natural isomorphism
$F^l_{X,M}:F(X\triangleright M)\to X\triangleright F(M)$ satisfying the usual
coherence axioms, and $\Fun_{\C}(\M,\N)$ denotes the category of such
functors and module natural transformations. Bimodule functors carry, in
addition, a right constraint $F^r_{M,Y}:F(M\triangleleft Y)\to
F(M)\triangleleft Y$. For a finite left $\C$-module category $\M$, we write
\[
  \C_{\M}^*:=\Fun_{\C}(\M,\M)^{\rev}
\]
for the monoidal category of $\C$-module endofunctors of $\M$. Its tensor
product is the reversed composition.

A finite $\C$-$\D$-bimodule category is a $\C$-$\D$-bimodule category $\M$
such that the left action $\tr: \C\times \M\to \M$ and the right action
$\tl: \M\times \D\to \M$ are $\kk$-linear and right exact in each variable.
Let $\M$ be a finite $\C$-$\D$-bimodule category. The functor
$\id_{\M} \triangleleft X$ for $X \in \D$ is a linear right exact left
$\C$-module endofunctor of $\M$, with module structure given by the middle
associativity constraint $(W \triangleright M) \triangleleft X \cong W \triangleright (M \triangleleft X)$ ($W \in \C$, $X \in \D$, $M \in \M$). The right associativity constraint $(M \triangleleft X) \triangleleft Y \cong M \triangleleft (X \otimes Y)$ makes these
functors into a strong monoidal functor
\[
  a_{\M} : \D \to \C_{\M}^*,
  \quad X \mapsto \id_{\M} \triangleleft X,
\]
which we call the \emph{action functor} of the right $\D$-action.

Action functors change along bimodule equivalences as follows.
If $\N$ is another finite $\C$-$\D$-bimodule category and
$F:\M\to \N$ is an equivalence of $\C$-$\D$-bimodule categories, then
conjugation by $F$ gives the inner tensor equivalence
\[
  \inn_F:\C_{\M}^*\to \C_{\N}^*,
  \qquad
  G\mapsto F\circ G\circ F^{-1}.
\]
For $\N=\M$, $\inn_F$ is the functor $X\mapsto\beta\otimes X\otimes\beta^*$ of \eqref{eq:def-inner-auto} in $\C_{\M}^*$ for $\beta=F^{-1}$, since the product of $\C_{\M}^*$ is the reversed composition.
For composable equivalences $G:\cL\to\M$ and $F:\M\to\N$, there is a
canonical isomorphism
\begin{equation}
  \label{eq:inner-equiv-composition}
  \gamma_{F,G}:\inn_F\circ \inn_G \xrightarrow{\sim} \inn_{FG}
\end{equation}
of strong monoidal functors, induced by $(FG)^{-1}\cong G^{-1}F^{-1}$. The same
formulas apply to equivalences of left $\C$-module categories.
For $X \in \mathcal{D}$, we have isomorphisms
\begin{equation*}
  \id_{\mathcal{N}} \triangleleft_{\mathcal{N}} X
  \cong F F^{-1} ((-) \triangleleft_{\mathcal{N}} X)
  \cong F (F^{-1}(-) \triangleleft_{\mathcal{M}} X)
\end{equation*}
of left $\mathcal{C}$-module functors. Hence we obtain an isomorphism
\begin{equation}
  \label{eq:action-functor-naturality}
  a_{\N}\cong \inn_F\circ a_{\M}
\end{equation}
of strong monoidal functors.
Conversely, the action functor determines the right action up to
equivalence:

\begin{lemma}\label{lem:bimodule-criterion}
  Let $\M$ and $\N$ be finite $\C$-$\D$-bimodule categories, and let
  $F:\M\to\N$ be an equivalence of left $\C$-module categories. If there is
  an isomorphism $\theta:a_{\N}\cong\inn_F\circ a_{\M}$ of monoidal functors,
  then $F$, with a suitable right constraint, is an equivalence of
  bimodule categories.
\end{lemma}
\begin{proof}
  Choose an adjoint equivalence $(F,G,\eta,\varepsilon)$ of left $\C$-module categories, with
  module natural unit $\eta:\id_{\M}\xrightarrow{\sim}GF$, so that
  $(\inn_F\circ a_{\M})(Y)=F\circ a_{\M}(Y)\circ G$. Put
  \[
    F^r_{M,Y}=\theta_Y(F(M))^{-1}\circ F\bigl(a_{\M}(Y)(\eta_M)\bigr):
    F(M\triangleleft Y)\xrightarrow{\ \sim\ }F(M)\triangleleft Y.
  \]
  Each $F^r_{-,Y}:F\circ a_{\M}(Y)\Rightarrow a_{\N}(Y)\circ F$ is a
  morphism of left $\C$-module functors, since $\eta$ and $\theta_Y$ are
  module natural. This is the compatibility with the middle associators.
  The right module axioms follow from the monoidality of $\theta$, the
  tensor structure of $\inn_F$ and the triangle identities.
\end{proof}

Before introducing relative Deligne products, we record a simple observation about the Deligne product of module functors, which will be used in \S\ref{subsec:groupoid-lagrangian}.
Let $\C$ and $\D$ be finite tensor categories, and let $\M$ and $\N$ be a finite left $\C$-module category and a finite left $\D$-module category, respectively. Then $\M \boxtimes \N$ is a finite left $\C \boxtimes \D$-module category.

\begin{lemma}
  \label{lem:Lagrangian-piv-proof-2}
    Let $\M$ and $\N$ be as above, and let $F : \M \to \M$ and $G : \N \to \N$ be right exact module endofunctors over $\C$ and $\D$, respectively. We assume that $\M$ and $\N$ are non-zero.
    If $F \boxtimes G \cong \id_{\M \boxtimes \N}$ as $\C \boxtimes \D$-module functors on $\M \boxtimes \N$, then there are isomorphisms $F \cong \id_{\M}$ and $G \cong \id_{\N}$ of module functors over $\C$ and $\D$, respectively.
\end{lemma}
\begin{proof}
  Let $V \in \M$ and $W \in \N$. Since $\Bbbk$ is algebraically closed, the object $V \boxtimes W \in \M \boxtimes \N$ is simple if and only if both $V$ and $W$ are. Assume that $V$ and $W$ are simple. Then, since $F(V) \boxtimes G(W) = (F \boxtimes G)(V \boxtimes W) \cong V \boxtimes W$, both $F(V)$ and $G(W)$ are simple objects. Moreover, we have
  \begin{equation*}
    \Hom_{\M}(F(V), V) \otimes_{\Bbbk} \Hom_{\N}(G(W), W)
    \cong \Hom_{\M \boxtimes \N}(F(V) \boxtimes G(W), V \boxtimes W) \ne 0.
  \end{equation*}
  Thus, by Schur's lemma, we have $F(V) \cong V$ and $G(W) \cong W$. There is an equivalence $\M \simeq \M \boxtimes W$ of $\C$-module categories. Since $G(W) \cong W$, the diagram
  \begin{equation*}
      \begin{tikzcd}[column sep = 48pt]
          \M \arrow[r, "{\simeq}"] \arrow[d, "{F}"']
          & \M \boxtimes W \arrow[d, "{F \boxtimes G|_{\M \boxtimes W}}"] \\
          \M \arrow[r, "{\simeq}"] & \M \boxtimes W
      \end{tikzcd}
  \end{equation*}
  is commutative up to isomorphism. Since $F \boxtimes G|_{\M \boxtimes W} \cong \id_{\M \boxtimes W}$, we obtain $F \cong \id_{\M}$. The same argument gives $G \cong \id_{\N}$.
\end{proof}

\subsection{Relative Deligne products}
\label{subsubsec:relative-Deligne}
Let $\M$ be a finite right $\D$-module category and $\N$ a finite left
$\D$-module category. Their relative Deligne product is a finite linear
category $\M\boxtimes_{\D}\N$ with a $\D$-balanced right exact bilinear
functor $(M,N)\mapsto M\boxtimes_{\D}N$, such that for every linear category
$\cA$, composition with it is an equivalence from right exact linear functors
$\M\boxtimes_{\D}\N\to\cA$ to $\D$-balanced right exact bilinear functors
$\M\times\N\to\cA$ \cite[Definitions~3.1 and~3.2]{douglas2019balanced}. It
exists, and $(M,N)\mapsto M\boxtimes_{\D}N$ is exact in each variable
\cite[Theorem~3.3]{douglas2019balanced}. We write
\[
  b_{M,Y,N}:(M\triangleleft Y)\boxtimes_{\D}N\xrightarrow{\ \sim\ }
  M\boxtimes_{\D}(Y\triangleright N)
  \qquad(Y\in\D)
\]
for the balancing isomorphisms.
We will repeatedly use the following consequence of the universal property of the relative Deligne product:

\begin{lemma}\label{lem:descent}
  Let $F,G:\M\boxtimes_{\D}\N\to\cA$ be right exact linear functors.
  Restriction to the objects $M\boxtimes_{\D}N$ is a bijection from natural
  transformations $F\Rightarrow G$ to families of morphisms
  $F(M\boxtimes_{\D}N)\to G(M\boxtimes_{\D}N)$ that are natural in $M$ and
  $N$ and compatible with the balancing isomorphisms. The same holds for an
  iterated product $(\cL\boxtimes_{\C}\M)\boxtimes_{\D}\N$ and families
  indexed by the objects $L\boxtimes_{\C}M\boxtimes_{\D}N$ that are natural in
  $L$, $M$ and $N$ and compatible with the $\C$-balancing in $(L,M)$ and the
  $\D$-balancing in $(M,N)$.
\end{lemma}
The same holds for $\cL\boxtimes_{\C}(\M\boxtimes_{\D}\N)$, via the associativity equivalence below.
\begin{proof}
  The first assertion is the full faithfulness of the equivalence in the
  universal property. For the second, apply the first assertion twice, using
  that $\mathcal Q\mapsto\mathcal Q\boxtimes_{\D}N$ is exact: for fixed $N$, a
  natural transformation between the right exact functors
  $\mathcal Q\mapsto F(\mathcal Q\boxtimes_{\D}N)$ and
  $\mathcal Q\mapsto G(\mathcal Q\boxtimes_{\D}N)$ is determined by its
  components at the objects $L\boxtimes_{\C}M$.
\end{proof}

If $\M$ is a $\C$-$\D$-bimodule category and $\N$ is a $\D$-$\E$-bimodule
category, then $\M\boxtimes_{\D}\N$ is a $\C$-$\E$-bimodule category with
\[
  X\triangleright(M\boxtimes_{\D}N)=(X\triangleright M)\boxtimes_{\D}N,
  \qquad
  (M\boxtimes_{\D}N)\triangleleft Z=M\boxtimes_{\D}(N\triangleleft Z);
\]
the actions are induced by the universal property, and their constraints by
Lemma~\ref{lem:descent} \cite[Remark~3.7]{douglas2019balanced}. In the same
way, bimodule equivalences $F:\M\to\M'$ and $G:\N\to\N'$ induce a bimodule
equivalence $F\boxtimes_{\D}G : \M \boxtimes_{\D} \N \to \M' \boxtimes_{\D} \N'$. The universal property also gives the unit
and associativity equivalences
\[
  \C\boxtimes_{\C}\M\simeq\M\simeq\M\boxtimes_{\D}\D,
  \qquad
  (\cL\boxtimes_{\C}\M)\boxtimes_{\D}\N\simeq\cL\boxtimes_{\C}(\M\boxtimes_{\D}\N),
\]
given on objects by $X\boxtimes_{\C}M\mapsto X\triangleright M$,
$M\boxtimes_{\D}Y\mapsto M\triangleleft Y$, and the identity of
$L\boxtimes_{\C}M\boxtimes_{\D}N$. With these, finite tensor categories,
finite bimodule categories, bimodule functors and their natural
transformations form a 3-category
\cite[Theorem~2.2.18]{douglas2018dualizable}.

\subsection{Exact module categories}
\label{subsec:exact-module-categories}

A finite left $\C$-module category $\M$ is \emph{exact} if $P\triangleright M$
is projective for all objects $M\in\M$ and all projective objects $P\in\C$,
and \emph{indecomposable} if it is non-zero and it is not equivalent to a direct sum of two non-zero $\C$-module categories. If $\id_{\M}$ is a simple object of $\C_{\M}^*$, then
$\M$ is indecomposable, since a decomposition $\M\simeq\M_1\oplus\M_2$ gives
$\id_{\M}\cong\id_{\M_1}\oplus\id_{\M_2}$.

For an indecomposable exact left $\C$-module category $\M$, the dual category $\C_{\M}^*$ is a finite tensor category \cite[\S7.12]{etingof2016tensor}. Finite tensor categories $\C$ and $\D$ are called \emph{(categorically) Morita equivalent} if $\D\simeq\C_{\M}^*$ as tensor categories for an exact left $\C$-module
category $\M$ \cite[Definition~7.12.17]{etingof2016tensor}. Such an $\M$ is
indecomposable by the observation above, since $\id_{\M}$ is the unit object
of $\C_{\M}^*\simeq\D$, which is simple.

By a \emph{module full subcategory} of a finite left $\C$-module category $\M$, we mean a full subcategory closed under finite direct sums, subquotients, and
the action of $\C$.

\begin{lemma}\label{lem:simple-generated}
  Let $\M$ be an indecomposable exact left $\C$-module category. Then the only non-zero module full subcategory of $\M$ is $\M$ itself.
\end{lemma}
\begin{proof}
  Let $\mathcal S$ be a nonzero module full subcategory of $\M$. It contains a
  simple object $M$. Let $N\in\M$ be simple. Since $\M$ is indecomposable and exact,
  $N$ is a subquotient of $X\triangleright M\in\mathcal S$ for some $X\in\C$
  \cite[Lemma~7.6.6 and Proposition~7.6.7]{etingof2016tensor}, so
  $N\in\mathcal S$.
  Now let $P_{\unit}$ be the projective cover of $\unit$.
  Then $P_{\unit} \triangleright N$ belongs to $\mathcal{S}$, is projective in $\M$ and surjects onto $N$.
  Thus it has the projective cover of $N$ as its direct summand.
  As every object is a quotient of a finite direct sum of projective covers of simple objects, it follows that $\mathcal S=\M$.
\end{proof}

\subsection{Invertible bimodule categories}
\label{subsubsec:invertible-bimodules}
We call a finite $\C$-$\D$-bimodule category $\M$ \emph{invertible} if there is a finite $\D$-$\C$-bimodule category
$\overline{\M}$ such that
\[
  \M\boxtimes_{\D}\overline{\M}\simeq \C,
  \qquad
  \overline{\M}\boxtimes_{\C}\M\simeq \D
\]
as bimodule categories. We have the following criterion for invertibility (due to \cite{fuchs2022spherical}), which is a generalization of \cite[Proposition~4.2]{etingof2010fusion} to finite tensor categories.

\begin{proposition}\label{prop:invertibility-criterion}
  A finite $\C$-$\D$-bimodule category $\M$ is invertible if and only if its action functor $a_{\M}:\D\to\C_{\M}^*$ is an equivalence. In particular, $\M$ is
  indecomposable and exact as a left $\C$-module category.
\end{proposition}
\begin{proof}
  If $\M$ is invertible, then $a_{\M}$ is an equivalence by \cite[Theorem~3.11(iv)]{fuchs2022spherical}, since an adjoint equivalence of $\M$ with an inverse bimodule category is a strong Morita context. If
  $a_{\M}$ is an equivalence, then $\id_{\M}\cong a_{\M}(\unit)$ is simple, so $\M$
  is indecomposable, and so is $A$ for $\M\simeq\C_A$ \cite[Theorem~2.2.12]{douglas2018dualizable}. Since
  $\C_{\M}^*\simeq{}_A\C_A$ is rigid, $A$ is exact
  \cite[Theorem~5.1]{shimizu2024commutative}, so $\M$ is exact (as in the proof of
  \cite[Theorem~3.11(i)]{fuchs2022spherical}).
  Then $\M$, with the right action of $\C_{\M}^*$, is an invertible
  $\C$-$\C_{\M}^*$-bimodule category \cite[Theorem~3.9 and Proposition~4.9]{fuchs2022spherical}, and restricting along the tensor equivalence $a_{\M}$ shows that $\M$ is invertible.
\end{proof}

Proposition~\ref{prop:invertibility-criterion} and the definition of Morita equivalence imply that finite tensor categories $\C$ and $\D$ are Morita equivalent if and only if there is an invertible $\C$-$\D$-bimodule category.

\subsection{The Schauenburg equivalence}
\label{subsec:Schauenburg-equiv}

Let $\M$ be an indecomposable exact left $\C$-module category. The center $\Z(\C)$ acts on
$\M$ through the forgetful functor $\Forg_{\C}:\Z(\C)\to \C$. For
$\mathbf{X}=(X,\sigma)$ in $\Z(\C)$, the endofunctor $\Schau_{\M}(\mathbf X)=X\triangleright \id_{\M}$
is a $\C$-module functor with the structure
\begin{gather*}
  (W \triangleright \Schau_{\M}(\mathbf{X}))(M)
  = W \triangleright (X \triangleright M)
  \cong (W \otimes X) \triangleright M \\
  \xrightarrow{\quad \sigma_W^{-1} \quad}
  (X \otimes W) \triangleright M
  \cong X \triangleright (W \triangleright M)
  = \Schau_{\M}(\mathbf{X})(W \triangleright M)
\end{gather*}
for $W \in \C$ and $M \in \M$, where $\cong$ means the module associator, and an object of $\Z(\C_{\M}^*)$ with the half-braiding \cite[\S3]{shimizu2020further}
\begin{equation}
  \label{eq:Sch-equiv-half-br}
  \Schau_{\M}(\mathbf{X}) \otimes F
  = F (X \triangleright \mathord{?})
  \xrightarrow{\quad s_{X, ?} \quad}
  X \triangleright F (\mathord{?})
  = F \otimes \Schau_{\M}(\mathbf{X})
\end{equation}
for $F \in \C_{\M}^*$, where $s_{X, M} : F(X \triangleright M) \to X \triangleright F(M)$ is the module structure of $F$ and $\otimes$ is the monoidal product of $\C_{\M}^*$.
This gives a braided tensor equivalence \cite[Theorem~3.13]{shimizu2020further}
\begin{equation}
  \label{eq:def-Schauenburg-equiv}
  \Schau_{\M}:\Z(\C)\to \Z(\C_{\M}^*),
  \qquad
  \mathbf{X}\mapsto \mathbf{X}\triangleright \id_{\M}.
\end{equation}
When $\M\simeq \C_A$ for an algebra $A$ in $\C$, this is identified with
Schauenburg's equivalence $\Z(\C)\simeq \Z({}_A\C_A)$,
$\mathbf V\mapsto A\otimes\mathbf V$ \cite{schauenburg2001monoidal},
followed by the equivalence induced by
${}_A\C_A\simeq\C_{\M}^*$, $X\mapsto(-)\otimes_A X$. See the proof of
\cite[Theorem~3.13]{shimizu2020further}. We therefore call $\Schau_{\M}$ the {\em Schauenburg
equivalence} attached to $\M$.

Since the tensor product of $\C_{\M}^*$ is the reversed composition, the
tensor structure of $\Schau_{\M}$ is
\begin{equation}
  \label{eq:Schauenburg-tensor-structure}
  \Schau_{\M}(\mathbf X)\otimes\Schau_{\M}(\mathbf X')
  =X'\triangleright(X\triangleright(-))
  \cong(X'\otimes X)\triangleright(-)
  \xrightarrow{\ c_{\mathbf X',\mathbf X}\triangleright(-)\ }
  (X\otimes X')\triangleright(-)
  =\Schau_{\M}(\mathbf X\otimes\mathbf X'),
\end{equation}
where $c_{\mathbf X',\mathbf X}=\sigma'_X$ is the braiding of $\Z(\C)$, for
$\mathbf X'=(X',\sigma')$. We fix this tensor structure. It is the one for
which the half-braidings $\sigma_V$ form the unitor of the pseudofunctor in
Lemma~\ref{lem:pseudofunctor-Z-2}, and $\Schau_{\M}$ is braided for it: by
\eqref{eq:Sch-equiv-half-br}, the braiding
$\Schau_{\M}(\mathbf X)\otimes\Schau_{\M}(\mathbf X')\to
\Schau_{\M}(\mathbf X')\otimes\Schau_{\M}(\mathbf X)$ is
$c_{\mathbf X',\mathbf X}\triangleright(-)$.

The Schauenburg equivalence has the following naturality (see \cite[Theorem~3.13]{shimizu2020further}):

\begin{lemma}
An equivalence $F:\M\to \N$ of left $\C$-module categories gives an isomorphism of tensor functors
\begin{equation}
  \label{eq:Schauenburg-equiv-naturality}
  \alpha_F:\Schau_{\N}\xrightarrow{\ \sim\ }\widetilde{\inn_F}\circ \Schau_{\M},
\end{equation}
where $\inn_F:\C_{\M}^* \to \C_{\N}^*$ is the inner tensor equivalence of
\S\ref{ssec:module-bimodule-center-prelim} and $\widetilde{\inn_F}$ is the
induced equivalence between Drinfeld centers
(\S\ref{subsec:induced-autoequivalences}).
\end{lemma}
\begin{proof}
We fix an adjoint equivalence $(F^{-1},F,\eta,\varepsilon)$ of left $\C$-module categories. Then, for $\mathbf{X} \in \Z(\C)$, we define $\alpha_F$ to be the composition
\begin{gather*}
  \Schau_{\N}(\mathbf{X})
  = \mathbf{X} \triangleright_{\N} \id_{\N}
  \cong F F^{-1}(\mathbf{X} \triangleright_{\N} (\mathord{?}))
  \cong F(\mathbf{X} \triangleright_{\M} F^{-1}(?)) \\
  = F \circ (\mathbf{X} \triangleright_{\M} \id_{\M}) \circ F^{-1}
  = (\widetilde{\inn_F} \circ \Schau_{\M})(\mathbf{X}),
\end{gather*}
where the left $\C$-module structure of $F^{-1}$ is used at the second isomorphism.

It is straightforward to check that $\alpha_F(\mathbf X)$ is a morphism of $\Z(\C_{\N}^*)$, is natural in $\mathbf X$, and is monoidal. The half-braidings and module structures of both sides are composites of the module constraints of $F$ and $F^{-1}$ and of the maps $\sigma_W^{\pm1}\triangleright(-)$. These maps commute with the constraints by naturality in the $\C$-variable, and $\eta$ and $\varepsilon$ are module transformations. For monoidality, the tensor structure of $\Schau_{\N}$ is \eqref{eq:Schauenburg-tensor-structure} for $\N$, and that of $\widetilde{\inn_F}\circ\Schau_{\M}$ is the tensor structure of $\widetilde{\inn_F}$ composed with \eqref{eq:Schauenburg-tensor-structure} for $\M$.
\end{proof}

The isomorphism $\alpha_F$ is the one used to define $\chi_F$ in \eqref{eq:def-chi-F-tilde}. When $\M = \N$, the isomorphism $\alpha_F$ is given as follows:

\begin{lemma}
  \label{lem:Sch-equiv-alpha}
  Let $\M$ be an indecomposable exact left $\C$-module category. Then we have
  \begin{equation*}
    \alpha_F = \left(
      \Schau_{\M}
      = \id_{\Z(\C_{\M}^*)} \circ \Schau_{\M}
      \xrightarrow{\quad \omega_{F^{-1}} \hcomp \id \quad}
      \widetilde{\inn_F} \circ \Schau_{\M} \right).
  \end{equation*}
  for every invertible object $F \in \C_{\M}^*$.
\end{lemma}
\begin{proof}
  The claim is equivalent to the statement that the equation
  \begin{equation*}
    \alpha_F(\mathbf{X}) = \omega_{F^{-1}}(\Schau_{\M}(\mathbf{X}))
  \end{equation*}
  holds for all objects $\mathbf{X} \in \Z(\C)$.
  Let $(F^{-1}, F, \eta, \varepsilon)$ be the adjoint equivalence fixed in the definition of $\alpha_F$, so that $\eta : \id_{\M} \to F \circ F^{-1}$ and $\varepsilon : F^{-1} \circ F \to \id_{\M}$ are isomorphisms of left $\C$-module functors. We take $(F, \eta, \varepsilon)$ as a left dual object of $F^{-1} \in \C_{\M}^*$. Then, by definition, the right hand side of the equation of the claim is given by the composition
  \begin{equation*}
    \Schau_{\M}(\mathbf{X})
    \xrightarrow{\quad \id \otimes \eta \quad}
    \Schau_{\M}(\mathbf{X}) \otimes F^{-1} \otimes F
    \xrightarrow{\quad \sigma(F^{-1}) \quad} F^{-1} \otimes \Schau_{\M}(\mathbf{X}) \otimes F
  \end{equation*}
  for $\mathbf{X} \in \Z(\C)$, where $\otimes$ is the tensor product of $\C_{\M}^*$ and $\sigma$ is the half-braiding of $\Schau_{\M}(\mathbf{X})$ given by~\eqref{eq:Sch-equiv-half-br}. This is exactly the same as $\alpha_F(\mathbf{X})$. The proof is done.
\end{proof}

\section{Lagrangian algebras in the Drinfeld center}
\label{sec:lagrangian-algebras}
Throughout this section, $\C$ is a finite tensor category. To every
indecomposable exact left $\C$-module category $\M$ we attach its adjoint
algebra $\bbA_{\M}$, a Lagrangian algebra in $\Z(\C)$
(\S\ref{subsec:adjoint-algebra}). The assignment $\M\mapsto\bbA_{\M}$
induces a bijection between equivalence classes of such module categories and
isomorphism classes of Lagrangian algebras in $\Z(\C)$
(Theorem~\ref{thm:classification-Lagrangian}). In
\S\ref{subsec:frobenius-lagrangian} we show that $\bbA_{\M}$ is Frobenius
if and only if $\M$ is unimodular, and, for pivotal $\C$, that $\bbA_{\M}$ is
symmetric Frobenius if and only if $\M$ is pivotal and unimodular.
\S\ref{subsec:groupoid-lagrangian} shows that every algebra morphism
$\bbA_{\M}\to\bbA_{\N}$ is induced by a module equivalence $\M\to\N$, unique
up to isomorphism. This is the input for \S\ref{sec:BrPic-2-groupoid}.

\subsection{The adjoint algebra}
\label{subsec:adjoint-algebra}
Recall the Schauenburg equivalence $\Schau_{\M}$ of
\S\ref{subsec:Schauenburg-equiv}. The tensor functor
$\Forg_{\C_{\M}^*}\circ\Schau_{\M}:\Z(\C)\to\C_{\M}^*$ is central and
surjective, since $\Schau_{\M}$ is an equivalence and $\Forg_{\C_{\M}^*}$ is
surjective. Its right adjoint is the lax monoidal functor
$\Schau_{\M}^{-1}\circ\Induc_{\C_{\M}^*}$.

\begin{definition}
  For an indecomposable exact left $\C$-module category $\M$, define
  \begin{equation}
    \label{eq:def-adjoint-alg-module-category}
    \bbA_{\M}
    :=
    (\Schau_{\M}^{-1}\circ \Induc_{\C_{\M}^*})(\unit_{\C_{\M}^*})
    \in \Z(\C),
  \end{equation}
  and call it the adjoint algebra of $\M$.
\end{definition}

For $\M=\C$, this recovers $\Induc_{\C}(\unit)$, the adjoint algebra of $\C$
\cite[Remark~6.10]{shimizu2024commutative}. More generally, $\bbA_{\mathcal M}$ is a Lagrangian algebra in $\Z(\C)$ \cite[\S6.2.2]{shimizu2024commutative}.

Up to isomorphism, $\bbA_{\M}$ depends only on the equivalence class of $\M$.
Indeed, for an equivalence $F:\M\to \N$ of left $\C$-module categories,
\eqref{eq:F-tilde-and-induction} and \eqref{eq:Schauenburg-equiv-naturality}
give an isomorphism
\begin{equation}
  \label{eq:def-chi-F-tilde}
  \Schau_{\N}^{-1}\circ \Induc_{\C_{\N}^*}\circ \inn_F
  \cong
  \Schau_{\N}^{-1}\circ \widetilde{\inn_F}\circ \Induc_{\C_{\M}^*}
  \cong
  \Schau_{\M}^{-1}\circ \Induc_{\C_{\M}^*}
\end{equation}
of lax monoidal functors. Evaluating at the unit and using
$\inn_F(\unit_{\C_{\M}^*})\cong \unit_{\C_{\N}^*}$ gives an algebra
isomorphism $\bbA_{\N}\to\bbA_{\M}$. We denote by
\begin{equation}
  \label{eq:def-chi-F}
  \chi_F:\bbA_{\M}\xrightarrow{\sim}\bbA_{\N}
\end{equation}
the inverse of this isomorphism. It depends only on the isomorphism class of the module functor $F$, since an isomorphism $F\cong F'$ of module functors induces an isomorphism
$\inn_F\cong\inn_{F'}$ of tensor functors compatible with \eqref{eq:F-tilde-and-induction} and
\eqref{eq:Schauenburg-equiv-naturality}.

The isomorphisms \eqref{eq:def-chi-F} are functorial in $F$:

\begin{lemma}
  \label{lem:chi-F-G}
  For $\C$-module equivalences $F : \M \to \N$ and $G : \cL \to \M$ between indecomposable exact left $\C$-module categories, we have
  \begin{equation*}
    \chi_{F} \circ \chi_{G} = \chi_{F G}.
  \end{equation*}
\end{lemma}
\begin{proof}
  Put $R_{\mathcal M}=\Schau_{\mathcal M}^{-1}
  \circ \Induc_{\mathcal C_{\mathcal M}^*}$, and similarly for
  $\mathcal L$ and $\mathcal N$. Let
  $c_F:R_{\mathcal M}\to R_{\mathcal N} \circ \inn_F$ denote the inverse of
  \eqref{eq:def-chi-F-tilde}. With
  $\gamma=\gamma_{F,G}$ from \eqref{eq:inner-equiv-composition}, we have
  \[
    c_{FG}=(R_{\mathcal N}\gamma)\circ(c_F\inn_G)\circ c_G.
  \]
  To verify it, expand the two comparisons defining $c_F$. On the
  Schauenburg side, it reduces to the identity
  $\alpha_{FG}=(\widetilde{\gamma}\,\Schau_{\cL})\circ(\widetilde{\inn_F}\,\alpha_G)\circ\alpha_F$
  of isomorphisms $\Schau_{\N}\to\widetilde{\inn_{FG}}\circ\Schau_{\cL}$, which
  holds because the module constraint of $FG$ is the composite of those of
  $G$ and $F$. On the induction side, it is the compatibility of the
  isomorphisms \eqref{eq:F-tilde-and-induction} with composition, which holds
  because iterated and direct mates of the square of forgetful functors have
  the same adjunct.

  Write $j_F:\unit\to\inn_F(\unit)$ for the tensor unit map. Then
  $\chi_F=R_{\mathcal N}(j_F^{-1})c_F(\unit)$ and
  $j_{FG}=\gamma_{\unit}\inn_F(j_G)j_F$.
  Substitute these two identities into the displayed equality and use
  naturality of $c_F$ at $j_G^{-1}$. This gives
  $\chi_{FG}=\chi_F\chi_G$, as claimed.
  Figure~\ref{fig:proof-lem-chi-F-G} depicts this calculation. In the diagram, $(\heartsuit)$ is the displayed identity,
  $(\spadesuit)$ is the identity for $j_{FG}$, and (nat.) is the naturality
  of $c_F$.
\end{proof}

\begin{figure}
  \begin{equation*}
    \begin{tikzcd}[column sep = -40pt, row sep = 32pt]
      \bbA_{\cL}
      \arrow[rr, "{\chi_G}"]
      \arrow[d, equal]
      & & \bbA_{\M}
      \arrow[rr, "{\chi_F}"]
      \arrow[d, equal]
      & & \bbA_{\N}
      \arrow[d, equal]
      \\[-24pt]
      (\Schau_{\cL}^{-1} \circ \Induc_{\C_{\cL}^*})(\unit)
      \arrow[rd, "{\eqref{eq:def-chi-F-tilde}}"]
      \arrow[rr, phantom, "{\text{(def.)}}"]
      & & (\Schau_{\M}^{-1} \circ \Induc_{\C_{\M}^*})(\unit)
      \arrow[ld, <-]
      \arrow[rd, "{\eqref{eq:def-chi-F-tilde}}"]
      \arrow[rr, phantom, "{\text{(def.)}}"]
      & & (\Schau_{\N}^{-1} \circ \Induc_{\C_{\N}^*})(\unit)
      \arrow[ld, <-]
      \\
      & \quad (\Schau_{\M}^{-1} \circ \Induc_{\C_{\M}^*} \circ \inn_G)(\unit)
      \arrow[rd, "{\eqref{eq:def-chi-F-tilde}}"']
      \arrow[rr, phantom, "{\text{(nat.)}}"]
      & & (\Schau_{\N}^{-1} \circ \Induc_{\C_{\N}^*} \circ \inn_F)(\unit) \quad
      \arrow[ld, <-]
      \\
      \arrow[rr, phantom, "{(\heartsuit)}"]
      & & (\Schau_{\N}^{-1} \circ \Induc_{\C_{\N}^*} \circ \inn_F \circ \inn_G)(\unit)
      \arrow[d, "{\eqref{eq:inner-equiv-composition}}"]
      & &
      \arrow[ll, phantom, "{(\spadesuit)}"]
      \\[-16pt]
      & & (\Schau_{\N}^{-1} \circ \Induc_{\C_{\N}^*} \circ \inn_{F G})(\unit)
      \arrow[to=2-1, rounded corners, <-, "{\eqref{eq:def-chi-F-tilde}}",
      to path={ -- (\tikztostart -| \tikztotarget) -- (\tikztotarget) \tikztonodes }]
      \arrow[to=2-5, rounded corners,
      to path={ -- (\tikztostart -| \tikztotarget) -- (\tikztotarget) \tikztonodes }]
    \end{tikzcd}
  \end{equation*}
  \caption{Proof of Lemma~\ref{lem:chi-F-G}}
  \label{fig:proof-lem-chi-F-G}
\end{figure}

Applying \cite[Proposition~6.6(c) and Corollary~6.7]{shimizu2024commutative},
as recalled before Lemma~\ref{lem:innF-chi-iso}, to the central surjective
functor $\Forg_{\C_{\M}^*}\circ\Schau_{\M}$ gives a tensor equivalence
\begin{equation}
  \label{eq:adjoint-alg-module-in-ZC}
  \Kappa_{\M}:\C_{\M}^*\to \Z(\C)_{\bbA_{\M}},
  \qquad
  \Kappa_{\M}\circ\Forg_{\C_{\M}^*}\circ\Schau_{\M}\cong F_{\bbA_{\M}},
\end{equation}
where $F_{\bbA_{\M}}:\mathbf V\mapsto\mathbf V\otimes\bbA_{\M}$ is the free
module functor.

Let $\mathcal{M}$ and $\mathcal{N}$ be indecomposable exact left $\mathcal{C}$-module
categories, and let $F:\mathcal{M}\to\mathcal{N}$ be a module equivalence.
Restriction of scalars along $\chi_F^{-1}$ gives a tensor equivalence
\[
  (\chi_F^{-1})^\star:
  \mathcal{Z}(\mathcal{C})_{\bbA_{\mathcal{M}}}
  \longrightarrow
  \mathcal{Z}(\mathcal{C})_{\bbA_{\mathcal{N}}}.
\]
Under the equivalences \eqref{eq:adjoint-alg-module-in-ZC} for $\mathcal{M}$
and $\mathcal{N}$, this functor corresponds to $\inn_F$:

\begin{lemma}
  \label{lem:functors-induced-by-isomorphisms}
  Let $F : \mathcal{M} \to \mathcal{N}$ be as above. Then the diagram
  \begin{equation*}
    \begin{tikzcd}[column sep = 48pt]
      \mathcal{C}_{\mathcal{M}}^*
      \arrow[r, "{\Kappa_{\mathcal{M}}}"]
      \arrow[d, "{\inn_F}"']
      & \mathcal{Z}(\mathcal{C})_{\bbA_{\mathcal{M}}}
      \arrow[d, "{(\chi_F^{-1})^{\star}}"] \\
      \mathcal{C}_{\mathcal{N}}^*
      \arrow[r, "{\Kappa_{\mathcal{N}}}"]
      & \mathcal{Z}(\mathcal{C})_{\bbA_{\mathcal{N}}}
    \end{tikzcd}
  \end{equation*}
  commutes up to a monoidal natural isomorphism.
\end{lemma}
\begin{proof}
  Apply Lemma~\ref{lem:innF-chi-iso} to
  $F_1=\Forg_{\C_{\M}^*}\circ\Schau_{\M}$,
  $F_2=\Forg_{\C_{\N}^*}\circ\Schau_{\N}$ and $G=\inn_F$, with
  $\tilde\chi$ the isomorphism \eqref{eq:def-chi-F-tilde}. Then
  $f=\chi_F^{-1}$ by \eqref{eq:def-chi-F}, and $K_1=\Kappa_{\M}$,
  $K_2=\Kappa_{\N}$.
\end{proof}

\subsection{Classification of Lagrangian algebras}
We now establish a bijection between the isomorphism classes of Lagrangian algebras in $\Z(\C)$ and the equivalence classes of indecomposable exact $\C$-module categories.
We need a theorem of Ostrik and Utiralova \cite[Theorem~5.13]{ostrik2026non}, which we reformulate in terms of adjoint algebras:

\begin{theorem}
  \label{thm:simple-commutative-sub-adjoint}
  For every simple commutative algebra $\bbL$ in $\Z(\C)$, there exists an
  indecomposable exact left $\C$-module category $\M$ such that $\bbL$ is a subalgebra of $\bbA_{\M}$.
\end{theorem}
\begin{proof}
  By Theorem~\ref{thm:CSZ}, $\bbL$ is indecomposable exact, so
  \cite[Theorem~5.13]{ostrik2026non} gives an indecomposable exact algebra
  $A\in\C$ such that $\bbL$ is a subalgebra of $I_A(\unit_{{}_A\C_A})$, where
  $I_A$ is the lax monoidal right adjoint of the tensor functor
  $\Phi_A:\Z(\C)\to{}_A\C_A$, $\mathbf V\mapsto\mathbf V\otimes A$. (Their
  half-braidings are inverse to ours; since $(V,\gamma)\mapsto(V,\gamma^{-1})$
  is a monoidal isomorphism over $\C$ between the two centers, this does not
  matter.) Put $\M=\C_A$, an indecomposable exact module category. Via the half-braiding, $\Phi_A$ is $\mathbf V\mapsto A\otimes\mathbf V$, so under the equivalence ${}_A\C_A\simeq\C_{\M}^*$ it corresponds to $\Forg_{\C_{\M}^*}\circ\Schau_{\M}$ as tensor functors by the description of $\Schau_{\M}$ after \eqref{eq:def-Schauenburg-equiv}. Taking right adjoints and using \eqref{eq:def-adjoint-alg-module-category}, we get $I_A(\unit_{{}_A\C_A})\cong\bbA_{\M}$ as algebras.
\end{proof}

The following result is a non-semisimple analogue of
\cite[Proposition~4.8]{davydov2013witt}, with a different proof:

\begin{theorem}
  \label{thm:classification-Lagrangian}
  The assignment
  \[
    \M\mapsto \bbA_{\M}
  \]
  induces a bijection between equivalence classes of indecomposable exact left
  $\C$-module categories and isomorphism classes of Lagrangian algebras in
  $\Z(\C)$.
\end{theorem}

Let $\M$ be an indecomposable exact left $\C$-module category, and put $A_{\M}=\Forg_{\C}(\bbA_{\M})$. One of the key observations to establish this theorem is that there is an equivalence
\begin{equation}
\label{eq:WY25-Prop-4-2}
    \C_{A_{\M}} \simeq \M \boxtimes \M^{\op}
\end{equation}
of left $\C$-module categories, with $\C$ acting on the first factor. To see
this, recall that there is an isomorphism $A_{\M}\cong\int_{M\in\M}\iHom(M,M)$
\cite[Corollary~3.15]{shimizu2020further} of algebras, where $\iHom(M,N)\in\C$ is the
internal Hom, characterized by
$\Hom_{\C}(X,\iHom(M,N))\cong\Hom_{\M}(X\triangleright M,N)$
\cite[\S7.9]{etingof2016tensor}. Under this isomorphism the
multiplication, which comes from the lax monoidal structure of the right
adjoint, is induced by composition of internal Homs. The functor
\[
  \Rex(\M,\M)\to\C_{A_{\M}},\qquad F\mapsto\int_{M\in\M}\iHom(M,F(M)),
\]
is an equivalence
\cite[Proposition~4.12(2)]{walton2025nondegenerate}.
If we view $\mathcal{E} := \Rex(\M, \M)$ as a left $\C$-module category by the action $(X \triangleright F)(M) = X \triangleright F(M)$ for $X \in \C$, $F \in \mathcal{E}$ and $M \in \M$, then the above equivalence is in fact an equivalence of left $\C$-module categories, since $\iHom(M,X\triangleright N)\cong
X\otimes\iHom(M,N)$ naturally \cite[\S7.9]{etingof2016tensor}.
Finally, there is the Eilenberg--Watts equivalence
$\M\boxtimes\M^{\op}\simeq\Rex(\M,\M)$,
$U\boxtimes V\mapsto\Hom_{\M}(-,V)^*\otimes_{\kk}U$, which is $\C$-linear in $U$.
By composing the equivalences, we obtain \eqref{eq:WY25-Prop-4-2}.

\begin{proof}
  The algebra isomorphism $\chi_F$ of \eqref{eq:def-chi-F} shows that the
  assignment depends only on the equivalence class of $\M$.

  \emph{Surjectivity.} Let $\bbL$ be a Lagrangian algebra in $\Z(\C)$. By
  Theorem~\ref{thm:simple-commutative-sub-adjoint}, there are an indecomposable exact left
  $\C$-module category $\M$ and an algebra monomorphism
  $\bbL\hookrightarrow \bbA_{\M}$. Both algebras have FP dimension
  $\FPdim(\C)$, so $\bbL\cong \bbA_{\M}$ by Remark~\ref{rem:FPdim}.

  \emph{Injectivity.} Suppose that $\bbA_{\M}\cong \bbA_{\N}$ as algebras in $\Z(\C)$. By \eqref{eq:WY25-Prop-4-2}, we have equivalences
  \begin{equation}
  \label{eq:classification-Lagrangian-proof-1}
  \M \boxtimes \M^{\op} \simeq \C_{A_{\M}} \simeq \C_{A_{\N}} \simeq \N \boxtimes \N^{\op}
  \end{equation}
  of left $\C$-module categories. Choose a simple object $T \in \M \boxtimes \M^{\op}$, and let $T'$ be the object of $\N \boxtimes \N^{\op}$ corresponding to $T$ under \eqref{eq:classification-Lagrangian-proof-1}. Since $\Bbbk$ is algebraically closed, we may assume that $T = V \boxtimes W^{\op}$ and $T' = X \boxtimes Y^{\op}$ for simple objects $V, W \in \M$ and $X, Y \in \N$.

  Consider the functor $F_W : \M \to \M \boxtimes \M^{\op}$, $F_W(M) = M \boxtimes W^{\op}$, and write $\M \boxtimes W^{\op}$ for its image. Since $W$ is absolutely simple, $F_W$ is exact and fully faithful, its image is closed under subquotients (writing $\M\simeq\Lambda\text{-mod}$ and $\M^{\op}\simeq\Gamma\text{-mod}$, a submodule of $M\otimes_{\Bbbk}W$ has the form $M'\otimes_{\Bbbk}W$, since $\Gamma$ acts on $W$ through all of $\mathrm{End}_{\Bbbk}(W)$), and it is an equivalence $\M \simeq \M \boxtimes W^{\op}$ of $\C$-module categories. Likewise, $\N \simeq \N \boxtimes Y^{\op}$.

  The equivalence \eqref{eq:classification-Lagrangian-proof-1} identifies the smallest $\C$-module full subcategory of $\M \boxtimes \M^{\op}$ containing $V \boxtimes W^{\op}$ with the smallest one of $\N \boxtimes \N^{\op}$ containing $X \boxtimes Y^{\op}$. Since $\M \boxtimes W^{\op}$ is a module full subcategory equivalent to $\M$, Lemma~\ref{lem:simple-generated} shows that the first of these is $\M \boxtimes W^{\op}$. By the same argument, the second is equivalent to $\N\boxtimes Y^{\op}$. Since these two subcategories are equivalent, we conclude $\M \simeq \N$ as $\C$-module categories.
\end{proof}

\begin{remark}\label{rem:Lagrangian}
  The inverse map sends a Lagrangian algebra $\bbL\cong\bbA_{\M}$ to any indecomposable exact
  module full subcategory of $\C_{\Forg(\bbL)}\simeq\M\boxtimes\M^{\op}$
  \eqref{eq:WY25-Prop-4-2}. Indeed, such a subcategory contains some
  $V\boxtimes W^{\op}$ with $V,W$ simple, so by Lemma~\ref{lem:simple-generated},
  applied to $\M\boxtimes W^{\op}$ and to the subcategory itself, it equals
  $\M\boxtimes W^{\op}\simeq\M$.
  In particular, let $A$ be a simple quotient algebra of $\Forg(\bbL)$ in
  $\C$. One exists since $\Forg(\bbL)$ is a nonzero algebra of finite length.
  Restriction along $\Forg(\bbL)\to A$ identifies $\C_A$ with a module full
  subcategory of $\C_{\Forg(\bbL)}$, which is indecomposable exact by
  Theorem~\ref{thm:CSZ}. Hence $\C_A\simeq\M$ as left $\C$-module categories.
  So the simple quotients of $\Forg(\bbL)$ are well defined up to Morita
  equivalence: any two are Morita equivalent, and each has $\M$ as its
  category of modules. Compare the proof of \cite[Theorem~5.13]{ostrik2026non}.
\end{remark}

\begin{remark}\label{rem:forgetful-injective-Lagrangian}
The forgetful functor $\Forg_{\C} : \Z(\C) \to \C$ is injective on isomorphism classes of Lagrangian algebras. Indeed, if two Lagrangian algebras have isomorphic underlying algebras in $\C$, then \eqref{eq:classification-Lagrangian-proof-1} identifies the corresponding categories. By the injectivity argument above, the corresponding module categories are equivalent. The assertion then follows from Theorem~\ref{thm:classification-Lagrangian}.
\end{remark}

\subsection{Frobenius property of the Lagrangian algebra}
\label{subsec:frobenius-lagrangian}

Next we discuss two finer properties of $\M$, unimodularity and pivotality, and how they are visible on $\bbA_{\M}$.

If $A$ is an algebra in a finite tensor category $\C$, then $A^*$ has a canonical structure of $A^{**}$-$A$-bimodule.
A \emph{Frobenius algebra} is an algebra $A$ in $\C$ such that $A \cong A^*$ as right $A$-modules.
When $\C$ is pivotal, we can view $A^*$ as an $A$-bimodule. The algebra $A$ is said to be \emph{symmetric Frobenius} if $A^* \cong A$ as $A$-bimodules.

We first discuss when the Lagrangian algebra $\bbA_{\M}$ is Frobenius. Let $\C$ be a finite tensor category, not necessarily pivotal.
An indecomposable exact left $\C$-module category $\M$ is \emph{unimodular} \cite[Definition~1.1]{yadav2023unimodular}
if the finite tensor category $\C^*_{\M}$ is unimodular in the sense of \cite{etingof2004analogue}.

\begin{proposition}\label{prop:Lagrangian-uni}
$\bbA_\M$ is Frobenius if and only if $\M$ is unimodular.
\end{proposition}
\begin{proof}
By \cite[Theorem~A]{yadav2023unimodular}, which builds on \cite{shimizu2016unimodular}, $\bbA_\M$ is Frobenius if and only if the finite tensor category $\Fun_{\C}(\M,\M)=(\C^*_\M)^{\rev}$ is unimodular. Reversing the tensor product does not affect unimodularity, so this is the definition of unimodularity of $\M$.
\end{proof}

\begin{remark}
  The problem of characterizing finite tensor categories that are not Morita
  equivalent to a unimodular finite tensor category was posed in
  \cite[\S4.5]{yadav2023unimodular}. By
  Theorem~\ref{thm:classification-Lagrangian} and Proposition~\ref{prop:Lagrangian-uni}, these are precisely the
  categories $\C$ for which $\Z(\C)$ does not admit a Lagrangian algebra that is Frobenius. Over a field of characteristic zero, the representation
  category of the Taft algebra is one such example
  \cite[Theorem~4.45]{yadav2023unimodular}.
\end{remark}

Fuchs, Schaumann and Schweigert introduced a relative Serre functor for an exact module category \cite{fuchs2020eilenberg}. We recall its definition from \cite{shimizu2024quasi} in a slightly generalized form: for a finite left $\C$-module category $\M$, the relative Serre functor $\mathbb{S}_{\M} : \M \to \M$ is the right exact linear functor equipped with a natural isomorphism
\begin{equation*}
    \underline{\Hom}(P, M)^* \cong \underline{\Hom}(M, \mathbb{S}_{\M}(P))
    \quad (P \in \mathcal{P}, M \in \M),
\end{equation*}
where $\mathcal{P}$ is the full subcategory of $\M$ consisting of objects $P$ such that the functor $\underline{\Hom}(P, -)$ is exact. When $\M$ is exact, we have $\mathcal{P} = \mathcal{M}$, and this definition agrees with that of \cite{fuchs2020eilenberg}. Unlike the definition of \cite{fuchs2020eilenberg}, the definition of \cite{shimizu2024quasi} works for any finite left $\C$-module category.

There is a canonical isomorphism $\mathbb{S}_{\M}(X \triangleright M) \cong X^{**} \triangleright \mathbb{S}_{\M}(M)$ ($X \in \C$, $M \in \M$) making $\mathbb{S}_{\M}$ a `twisted' $\C$-module functor.
If $\C$ is equipped with a pivotal structure, then the twisted $\C$-module functor $\mathbb{S}_{\M}$ is made into an ordinary $\C$-module functor via the pivotal structure. We say that $\M$ is \emph{pivotal} \cite[Definition 6.2]{shimizu2024quasi} if there exists an isomorphism $\mathbb{S}_{\M} \cong \id_{\M}$ of $\C$-module functors.

\begin{proposition}\label{prop:Lagrangian-piv}
    Let $\C$ be a pivotal finite tensor category, and let $\M$ be an indecomposable exact left $\C$-module category. We equip $\Z(\C)$ with the pivotal structure induced by that of $\C$. Then the Lagrangian algebra $\bbA_{\M}$ is symmetric Frobenius if and only if $\M$ is pivotal and unimodular.
\end{proposition}

For the proof, we recall a formula for the relative Serre functor. Let $\C$ be a finite tensor category, not necessarily pivotal, and let $A$ be an algebra in $\C$. For $\M = \C_A$, the relative Serre functor is given by
\begin{equation}
  \label{eq:rel-Serre-formula}
    \mathbb{S}_{\M}(M) = M^{**} \otimes_{A^{**}} A^* \quad (M \in \M)
\end{equation}
and the twisted module structure of $\mathbb{S}_{\M}$ is induced by the canonical isomorphism $(X \otimes M)^{**} \cong X^{**} \otimes M^{**}$ ($X, M \in \C$) \cite[Corollary 5.12]{shimizu2024quasi}.

Together with the Eilenberg--Watts theorem, \eqref{eq:rel-Serre-formula} implies that $\C_A$ is pivotal if $A$ is symmetric Frobenius, and that every pivotal exact module category is equivalent to $\C_A$ for some symmetric Frobenius algebra $A$ \cite[Theorem 6.4]{shimizu2024quasi}.

Now let $\C$ and $\D$ be finite tensor categories, not necessarily pivotal, and let $\M$ and $\N$ be a finite left $\C$-module category and a finite left $\D$-module category, respectively. Then $\M \boxtimes \N$ is a finite left $\C \boxtimes \D$-module category, and the relative Serre functor is multiplicative with respect to the Deligne tensor product:

\begin{lemma}
  \label{lem:Lagrangian-piv-proof-1}
    $\mathbb{S}_{\M \boxtimes \N} \cong \mathbb{S}_{\M} \boxtimes \mathbb{S}_{\N}$ as twisted module functors.
\end{lemma}
\begin{proof}
    We may assume that $\M = \C_A$ and $\N = \D_B$ for some algebras $A$ and $B$ in $\C$ and $\D$, respectively \cite[Theorem~2.2.12]{douglas2018dualizable}. By \eqref{eq:rel-Serre-formula}, we have natural isomorphisms
    \begin{align*}
        \mathbb{S}_{\M \boxtimes \N}(M \boxtimes N)
        & \cong (M \boxtimes N)^{**} \otimes_{(A \boxtimes B)^{**}} (A \boxtimes B)^* \\
        & \cong (M^{**} \otimes_{A^{**}} A^*) \boxtimes (N^{**} \otimes_{B^{**}} B^*)
        \cong \mathbb{S}_{\M}(M) \boxtimes \mathbb{S}_{\N}(N)
    \end{align*}
    for $M \in \M$ and $N \in \N$. The resulting isomorphism $\mathbb{S}_{\M \boxtimes \N} \cong \mathbb{S}_{\M} \boxtimes \mathbb{S}_{\N}$ is compatible with the twisted module structures, since both are induced by the canonical isomorphisms $(X \otimes M)^{**} \cong X^{**} \otimes M^{**}$.
\end{proof}

Now, using Lemmas~\ref{lem:Lagrangian-piv-proof-1} and~\ref{lem:Lagrangian-piv-proof-2}, we can prove Proposition~\ref{prop:Lagrangian-piv}.

\begin{proof}[Proof of Proposition \ref{prop:Lagrangian-piv}]
  Recall that $\M$ is indecomposable exact. The `if' is exactly \cite[Corollary~6.15(ii)]{yadav2026Frobenius}. Conversely, assume that $\mathbb{A}_{\M}$ is a symmetric Frobenius algebra. Since the forgetful functor $\Forg_{\C} : \Z(\C) \to \C$ is a tensor functor preserving the pivotal structures, the algebra $A_{\M} = \Forg_{\C}(\mathbb{A}_{\M})$ is symmetric Frobenius, and hence $\C_{A_{\M}}$ is a pivotal $\C$-module category. By \eqref{eq:WY25-Prop-4-2}, $\M \boxtimes \M^{\op} \simeq \C_{A_{\M}}$ as $\C$-module categories. Lemmas \ref{lem:Lagrangian-piv-proof-1} and \ref{lem:Lagrangian-piv-proof-2}, applied with $\D=\Vect$ and $\N=\M^{\op}$, then give $\mathbb{S}_{\M} \cong \id_{\M}$ as module functors, so $\M$ is pivotal. Unimodularity follows from Proposition~\ref{prop:Lagrangian-uni}, because a symmetric Frobenius algebra is Frobenius.
\end{proof}

\subsection{The groupoid of Lagrangian algebras}
\label{subsec:groupoid-lagrangian}

For \S\ref{sec:BrPic-2-groupoid}, we need more than the isomorphism classes of Lagrangian algebras. We also need to consider the algebra isomorphisms between them.

Given two algebras $A$ and $B$ in the same monoidal category, write
$\Alg(A,B)$ for the set of algebra morphisms from $A$ to $B$. To understand this set for Lagrangian algebras, we need to recall the central Hopf comonad, see for instance the description in \cite[\S\S3.1 and~3.4]{shimizu2023pivotal}.
Let $\D$ be a finite tensor category.
For $V \in \D$, we define
\begin{equation*}
  Z(V) = \int_{X \in \D} X \otimes V \otimes X^*
\end{equation*}
and denote by $\pi_{X}(V) : Z(V) \to X \otimes V \otimes X^*$ ($V, X \in \D$) the universal dinatural transformation. The assignment $V \mapsto Z(V)$ extends to an endofunctor on the category $\D$. Moreover, the functor $Z$ is a comonad, with comultiplication and counit
\begin{equation*}
  \delta_Z(V) : Z(V) \to ZZ(V),
  \quad \varepsilon_Z(V) : Z(V) \to V
  \quad (V \in \D)
\end{equation*}
defined via the universal dinatural transformation $\pi_{-}(-)$.

The space $\Hom_{\D}(V, Z(W))$ is canonically identified with the space of natural transformations from $V \otimes \id_{\D}$ to $\id_{\D} \otimes W$. Moreover, a morphism $V \to Z(V)$ makes $V$ a comodule over the comonad $Z$ if and only if the corresponding natural transformation makes $V$ an object of $\Z(\D)$. This observation allows us to identify $\Z(\D)$ with the category of comodules over $Z$. Under this identification, the functor $\Induc_{\D}$ is identified with the cofree $Z$-comodule functor, that is,
\begin{equation*}
  \Induc_{\D}(V) = (Z(V), \delta_Z(V))
  \quad (V \in \D).
\end{equation*}
For an invertible object $\beta\in\D$ and $\mathbf{X}=(X,\sigma_X)\in\Z(\D)$, the isomorphism $\omega_\beta(\mathbf{X})$ of \S\ref{subsec:induced-autoequivalences} is also given by
\begin{equation}
  \label{eq:def-omega-by-coaction}
  \omega_\beta(\mathbf{X})
  = \pi_\beta(X) \circ \rho_X,
\end{equation}
where $\rho_X : X \to Z(X)$ is the $Z$-coaction corresponding to the half-braiding.

Let $\bbA_{\D} = \Induc_{\D}(\unit)$ be the adjoint algebra of $\D$.
The automorphism group of the algebra $\bbA_{\D} \in \Z(\D)$ has been
determined in \cite{shimizu2023pivotal}: there is a group isomorphism
\begin{equation}
  \label{eq:auto-of-adj-alg-1}
  \Inv(\D) \to \Alg_{\Z(\D)}(\bbA_{\D}, \bbA_{\D}),
  \quad [\beta] \mapsto \chi'_{\beta} := Z(\coev_{\beta}^{-1} \circ \pi_{\beta}(\unit)) \circ \delta_Z(\unit).
\end{equation}
See \cite[Lemma~3.10 and Remark~3.11]{shimizu2023pivotal}. For an indecomposable exact left $\C$-module category $\M$, there is the following map which was defined in \eqref{eq:def-chi-F}:
\begin{equation}
  \label{eq:auto-of-adj-alg-3}
  \Inv(\C_{\M}^*)
  \to \Alg_{\Z(\C)}(\bbA_{\M}, \bbA_{\M}),
  \quad [F] \mapsto \chi_F.
\end{equation}
It is a group anti-homomorphism by Lemma~\ref{lem:chi-F-G}.
The next lemma relates \eqref{eq:auto-of-adj-alg-1} and
\eqref{eq:auto-of-adj-alg-3}.

\begin{lemma}
  \label{lem:auto-of-adj-alg-1}
  For any indecomposable exact left $\C$-module category $\M$, the diagram
  \begin{equation*}
    \begin{tikzcd}[column sep = 96pt]
      \Inv(\C_{\M}^*)
      \arrow[d, "{(-)^{-1}}"']
      \arrow[r, "{\text{\eqref{eq:auto-of-adj-alg-1} with $\D = \C_{\M}^*$}}"]
      & \Alg_{\Z(\C_{\M}^*)}
      (\bbA_{\C_{\M}^*}, \bbA_{\C_{\M}^*})
      \arrow[d, "{\Schau_{\M}^{-1}}"] \\
      \Inv(\C_{\M}^*)
      \arrow[r, "{\eqref{eq:auto-of-adj-alg-3}}"]
      & \Alg_{\Z(\C)}
      (\bbA_{\M}, \bbA_{\M})
    \end{tikzcd}
  \end{equation*}
  commutes. In other words, we have
\begin{equation*}
    \chi_{F}^{} = \Schau_{\M}^{-1}(\chi'_{F^{-1}})
    \quad (F \in \Inv(\C_{\M}^*)).
  \end{equation*}
\end{lemma}
\begin{proof}
  Put $\D=\C_{\M}^*$, $S=\Schau_{\M}$, and let $\beta=F^{-1}$ in $\Inv(\D)$.
  Thus the left dual object of $\beta$ is given by $\beta^*=F$ with the evaluation and coevaluation coming from the chosen adjoint equivalence.
  Since the tensor product of $\D$ is reversed composition, the conjugation by the module equivalence $F$ is the inner tensor functor implemented by $\beta$:
  \[
    G:=\inn_F=\beta\otimes(-)\otimes \beta^*,
    \qquad G_0=\coev_\beta.
  \]
Lemma~\ref{lem:Sch-equiv-alpha} gives
  $\alpha_F=\omega_\beta S$. Let $\eta^S$ and $\varepsilon^S$ be the unit and
  counit of the adjoint equivalence $S\dashv S^{-1}$, and let
  $\alpha'_F=(S^{-1}\widetilde G\varepsilon^S)\circ(S^{-1}\alpha_FS^{-1})\circ(\eta^SS^{-1})$
  be the mate of $\alpha_F$. The second isomorphism in
  \eqref{eq:def-chi-F-tilde} is $(\alpha'_F)^{-1}$. Substituting
  $\alpha_F=\omega_\beta S$ and using the naturality of $\omega_\beta$ and a
  triangle identity gives
  \[
    \alpha'_F=S^{-1}\omega_\beta:
    S^{-1}\longrightarrow S^{-1}\widetilde G.
  \]
  Let $\widetilde{\xi}_G:\widetilde{G}\circ\Induc_{\D}\to\Induc_{\D}\circ G$ be the
  isomorphism \eqref{eq:F-tilde-and-induction} for $F=G$.
  By the definition of $\chi_F$, it remains to prove
  \begin{equation}
    \label{eq:lem-auto-of-adj-alg-1-2}
    \chi'_\beta
    =\Induc_{\D}(G_0^{-1})\circ\widetilde{\xi}_G(\unit)
      \circ\omega_\beta(\Induc_{\D}(\unit)).
  \end{equation}
  This is an equality of maps to the cofree comodule
  $\Induc_{\D}(\unit)$, so we may compare their adjuncts in $\D$.
  Write $\varepsilon=\varepsilon_Z$. Since $\widetilde{\xi}_G$ is the mate of the
  equality $\Forg_{\D}\circ\widetilde{G}=G\circ\Forg_{\D}$, naturality of
  $\varepsilon$ and a triangle identity give
  $\varepsilon_{G(V)}\circ\Forg_{\D}(\widetilde{\xi}_G(V))=G(\varepsilon_V)$ for $V\in\D$.
  Using this, the adjunct of the right-hand side is
  \begin{align*}
    &G_0^{-1}\circ G(\varepsilon_{\unit})
       \circ\omega_\beta(\Induc_{\D}(\unit))\\
    &\quad=\coev_\beta^{-1}\circ
       (\id_\beta\otimes\varepsilon_{\unit}\otimes\id_{\beta^*})
       \circ\pi_\beta(Z(\unit))\circ\delta_Z(\unit)\\
    &\quad=\coev_\beta^{-1}\circ\pi_\beta(\unit).
  \end{align*}
  The first equality uses \eqref{eq:def-omega-by-coaction}; the second
  uses naturality of $\pi_\beta$ and
  $Z(\varepsilon_{\unit})\delta_Z(\unit)=\id_{Z(\unit)}$.
  This is also the adjunct of $\chi'_\beta$ by
  \eqref{eq:auto-of-adj-alg-1}. The adjunction therefore proves
  \eqref{eq:lem-auto-of-adj-alg-1-2}, and applying $S^{-1}$ gives the claim.
\end{proof}

The set of algebra morphisms between Lagrangian algebras is described as follows:

\begin{lemma}
  \label{lem:isomorphisms-between-Lagrangian}
  Let $\M$ and $\N$ be indecomposable exact left $\C$-module categories.
  \begin{enumerate}
  \item If $\M\not\simeq \N$, then $\Alg(\bbA_{\M},\bbA_{\N})$ is empty.
  \item If $\M\simeq \N$ as left $\C$-module categories, then the map
    \[
      \mathrm{Eq}_{\C}(\M,\N)\to \Alg(\bbA_{\M},\bbA_{\N}),
      \qquad
      [F]\mapsto \chi_F,
    \]
    is a bijection, where $\mathrm{Eq}_{\C}(\M,\N)$ denotes isomorphism
    classes of $\C$-module equivalences from $\M$ to $\N$.
  \end{enumerate}
\end{lemma}
\begin{proof}
(1) We show the contrapositive. Suppose that there exists an algebra morphism $\bbA_{\M}\to\bbA_{\N}$. Since Lagrangian algebras in $\Z(\C)$ are simple and have the same FP dimension, every algebra morphism between them is an isomorphism: its kernel is a proper two-sided ideal of a simple algebra, hence zero, so the morphism is a monomorphism between algebras of equal FP dimension, and Remark~\ref{rem:FPdim} shows that it is an isomorphism. So $\bbA_{\M}\cong\bbA_{\N}$ and $\M\simeq\N$ by Theorem~\ref{thm:classification-Lagrangian}.

(2) Let $G : \M \to \N$ be an equivalence of left $\C$-module
  categories, and set $\D=\C_{\M}^*$.
We consider the following chain of bijections:
\begin{gather*}
  \mathrm{Eq}_{\C}(\M, \N)
  \mathrel{\mathop{\cong}^{\text{(i)}}}
  \mathrm{Eq}_{\C}(\M, \M)
  = \Inv(\D)^{\op}
  \mathrel{\mathop{\cong}^{\text{(ii)}}}
  \Inv(\D) \\
  \mathrel{\mathop{\cong}^{\text{(iii)}}}
  \Alg_{\Z(\D)}(\bbA_{\D}, \bbA_{\D})
  \mathrel{\mathop{\cong}^{\text{(iv)}}}
  \Alg_{\Z(\C)}(\bbA_{\M}, \bbA_{\M})
  \mathrel{\mathop{\cong}^{\text{(v)}}}
  \Alg(\bbA_{\M}, \bbA_{\N}),
\end{gather*}
where (i) is induced by the equivalence $G$,
(ii) is given by taking the inverse,
(iii) is the bijection \eqref{eq:auto-of-adj-alg-1},
(iv) is induced by $\Schau_{\M}^{-1}$,
and (v) is induced by the isomorphism $\chi_G$.
For $[H]\in\mathrm{Eq}_{\C}(\M,\N)$, Lemmas~\ref{lem:chi-F-G} and
\ref{lem:auto-of-adj-alg-1} give
\begin{gather*}
  [H]
  \mathrel{\mathop{\longmapsto}^{\text{(i)}}} [G^{-1} H]
  \mathrel{\mathop{\longmapsto}^{\text{(ii)}}} [(G^{-1} H)^{-1}]
  \mathrel{\mathop{\longmapsto}^{\text{(iii)}}} \chi'_{(G^{-1} H)^{-1}}
  \mathrel{\mathop{\longmapsto}^{\text{(iv)}}}
  \Schau_{\M}^{-1}(\chi'_{(G^{-1} H)^{-1}}) \\
  \mathrel{\mathop{\longmapsto}^{\text{(v)}}}
  \chi_G \circ \Schau_{\M}^{-1}(\chi'_{(G^{-1} H)^{-1}})
  = \chi_G \circ \chi_{G^{-1} H} = \chi_H.
\end{gather*}
This means that the composition of the bijections (i)--(v) is the map $[F] \mapsto \chi_F$ of part~(2), which is therefore bijective.
\end{proof}

\section{Brauer--Picard groupoid}
\label{sec:BrPic-2-groupoid}
We show that taking the Drinfeld center is a fully faithful pseudofunctor from the Brauer--Picard 2-groupoid of finite tensor categories to the 2-groupoid of
braided finite tensor categories.

The Brauer--Picard 2-groupoid $\uBrPic$ (see \cite[\S4]{etingof2010fusion} for fusion categories and
\cite[\S2G]{davydov2013picard} for the group $\BrPic(\C)$) has finite tensor categories as 0-cells. For
0-cells $\C$ and $\D$, the category $\uBrPic(\C,\D)$ of 1-cells from $\C$ to
$\D$ has invertible $\C$-$\D$-bimodule categories as objects and isomorphism
classes of bimodule equivalences as morphisms. We write
$\uBrPic(\C):=\uBrPic(\C,\C)$. Composition of 1-cells is the relative Deligne
product, written in functional order,
\[
  \uBrPic(\D,\E)\times \uBrPic(\C,\D)\to \uBrPic(\C,\E),
  \qquad
  (\M,\N)\mapsto \M\circ\N:=\N\boxtimes_{\D}\M,
\]
with the unit and associativity equivalences of
\S\ref{subsubsec:relative-Deligne}. Every 1-cell is invertible up to
equivalence and every 2-cell is invertible, so $\uBrPic$ is a 2-groupoid.

The 2-groupoid $\uEqBr$ has braided finite tensor categories as 0-cells. For
0-cells $\mathcal A$ and $\mathcal B$, the objects of
$\uEqBr(\mathcal A,\mathcal B)$ are the braided tensor equivalences
$\mathcal A\to\mathcal B$, and its morphisms are the monoidal natural
isomorphisms between them. Composition of 1-cells is composition of functors.

\begin{theorem}
  \label{thm:BrPic}
  The assignment $\C\mapsto \Z(\C)$ extends to a fully faithful pseudofunctor
  \[
    \Z:\uBrPic\to \uEqBr.
  \]
\end{theorem}
Here fully faithful means that for all 0-cells $\C$ and $\D$, the functor
\[
  \Z_{\C,\D}:\uBrPic(\C,\D)\to\uEqBr(\Z(\C),\Z(\D))
\]
is an equivalence of categories. We construct $\Z$ in
\S\ref{subsec:def-pseudofunctor-Z}. In
\S\ref{subsec:bijectivity-on-1-cells} we show that $\Z$ is bijective on
isomorphism classes of 1-cells. Full faithfulness is proved in
\S\ref{subsec:fully-faith-Z}. In \S\ref{subsec:BrPic-applications}
we obtain the isomorphism $\BrPic(\C)\cong\Aut_{\otimes}^{\mathrm{br}}(\Z(\C))$ and classify the finite tensor categories Morita equivalent to $\C$ by orbits of Lagrangian algebras (Theorem~\ref{thm:MME}).

\subsection{The construction of the pseudofunctor}
\label{subsec:def-pseudofunctor-Z}

We define $\Z$ on 1-cells by \eqref{eq:def-partial-M}, on 2-cells through
the central action isomorphism \eqref{eq:Phi-central-action}
(Lemma~\ref{lem:pseudofunctor-Z-1}), and construct its compositor and
unitor in Lemma~\ref{lem:pseudofunctor-Z-2}. All three steps rest on the
uniqueness statement of Lemma~\ref{lem:central-action-unique}.

Given an invertible $\C$-$\D$-bimodule category $\M$, define, in the
notation of \S\ref{subsec:Schauenburg-equiv},
\begin{equation}
  \label{eq:def-partial-M}
  \Phi_{\M}
  :=
  \widetilde{a_{\M}^{-1}}\circ \Schau_{\M}
  :
  \Z(\C)\to \Z(\D),
\end{equation}
where $a_{\M}:\D\to \C_{\M}^*$ is the action functor. We set
$\Z_{\C,\D}(\M)=\Phi_{\M}$.
Since $a_{\M}$ is a tensor equivalence by
Proposition~\ref{prop:invertibility-criterion}, $\Phi_{\M}$ is a braided
tensor equivalence.

The isomorphisms
$\widetilde{a_{\M}}\circ\widetilde{a_{\M}^{-1}}\cong\widetilde{a_{\M}a_{\M}^{-1}}\cong\id$
of \S\ref{subsec:induced-autoequivalences} (for each $\M$ we fix a quasi-inverse $a_{\M}^{-1}$ and an isomorphism $a_{\M}a_{\M}^{-1}\cong\id$) give an isomorphism
$\widetilde{a_{\M}}\circ\Phi_{\M}\cong\Schau_{\M}$ of braided tensor
functors into $\Z(\C_{\M}^*)$. Its underlying natural isomorphisms
\begin{equation}
  \label{eq:Phi-central-action}
  \kappa^{\M}_{\mathbf X}:
  \mathbf X\triangleright\id_{\M}
  \xrightarrow{\ \sim\ }
  \id_{\M}\triangleleft\Phi_{\M}(\mathbf X)
  \qquad(\mathbf X\in\Z(\C))
\end{equation}
say that acting by $\mathbf X$ on the left is the same as acting by
$\Phi_{\M}(\mathbf X)$ on the right. We call $\kappa^{\M}$ the
\emph{central action isomorphism} of $\M$. Since it underlies an
isomorphism of functors into $\Z(\C_{\M}^*)$, it has the following two
properties.
\begin{itemize}
\item Each $\kappa^{\M}_{\mathbf X}$ is an isomorphism of
  $\C$-$\D$-bimodule functors. Here $\mathbf X\triangleright\id_{\M}$ is a
  left $\C$-module functor through the half-braiding of $\mathbf X$, as in
  \S\ref{subsec:Schauenburg-equiv}, and a right $\D$-module functor through
  the middle associator. Dually, $\id_{\M}\triangleleft\Phi_{\M}(\mathbf X)$
  is a left $\C$-module functor through the middle associator and a right
  $\D$-module functor through the half-braiding of $\Phi_{\M}(\mathbf X)$.
  Right $\D$-linearity is the compatibility with the half-braidings against
  the objects $a_{\M}(Y)$, $Y\in\D$.
\item The family $\kappa^{\M}$ is a monoidal natural isomorphism
  $\Forg\circ\Schau_{\M}\Rightarrow a_{\M}\circ\Forg\circ\Phi_{\M}$ of
  tensor functors $\Z(\C)\to\C_{\M}^*$.
\end{itemize}
The central action isomorphism characterizes the tensor functor
$\Phi_{\M}$:

\begin{lemma}
  \label{lem:central-action-unique}
  Let $\M$ be an invertible $\C$-$\D$-bimodule category, let
  $T:\Z(\C)\to\Z(\D)$ be a braided tensor functor, and let
  $\kappa:\Forg\circ\Schau_{\M}\Rightarrow a_{\M}\circ\Forg\circ T$ be a
  monoidal natural isomorphism. Then there is a unique monoidal natural
  isomorphism $u:\Phi_{\M}\Rightarrow T$ such that
  \[
    \kappa_{\mathbf X}
    =(\id_{\M}\triangleleft u_{\mathbf X})\circ\kappa^{\M}_{\mathbf X}
    \qquad(\mathbf X\in\Z(\C)).
  \]
\end{lemma}
\begin{proof}
  The composite $\kappa\circ(\kappa^{\M})^{-1}$ is a monoidal isomorphism
  between the underlying tensor functors of the braided functors
  $\widetilde{a_{\M}}\circ\Phi_{\M}$ and $\widetilde{a_{\M}}\circ T$ from
  $\Z(\C)$ to $\Z(\C_{\M}^*)$. The first underlying functor is isomorphic to
  $\Forg\circ\Schau_{\M}$, which is surjective, because $\Schau_{\M}$ is an
  equivalence and the forgetful functor of a Drinfeld center is surjective
  \cite[Corollary~7.13.11]{etingof2016tensor}. By
  Lemma~\ref{lem:braided-lift-unique}, $\kappa\circ(\kappa^{\M})^{-1}$ is
  an isomorphism
  $\widetilde{a_{\M}}\circ\Phi_{\M}\Rightarrow\widetilde{a_{\M}}\circ T$ of
  functors into $\Z(\C_{\M}^*)$. Since $\widetilde{a_{\M}}$ is an
  equivalence, this isomorphism is $\widetilde{a_{\M}}(u)$ for a unique
  $u:\Phi_{\M}\Rightarrow T$, and $u$ is monoidal because
  $\widetilde{a_{\M}}$ is faithful and monoidal.
\end{proof}

Let $[F]:\M\Rightarrow \N$ be a 2-cell in $\uBrPic$, represented by an
equivalence $F:\M\to \N$ of $\C$-$\D$-bimodule categories with module
structures
\[
  F^l_{X,M}:F(X\triangleright M)\to X\triangleright F(M),
  \qquad
  F^r_{M,Y}:F(M\triangleleft Y)\to F(M)\triangleleft Y.
\]
The comparisons \eqref{eq:action-functor-naturality} and
\eqref{eq:Schauenburg-equiv-naturality} are built from $F^r$ and $F^l$,
respectively. Transporting $\kappa^{\M}$ along them gives a monoidal natural
isomorphism $\Forg\circ\Schau_{\N}\Rightarrow a_{\N}\circ\Forg\circ\Phi_{\M}$.
Lemma~\ref{lem:central-action-unique}, applied to $\N$ and $T=\Phi_{\M}$,
gives a monoidal isomorphism $\Phi_{\N}\Rightarrow\Phi_{\M}$. We take
$\Z_{\C,\D}([F]):\Phi_{\M}\Rightarrow\Phi_{\N}$ to be its inverse, which is
the unique monoidal natural isomorphism such that
\begin{equation}
  \label{eq:def-partial-2cell}
  \bigl(F(M)\triangleleft\Z_{\C,\D}([F])_{\mathbf X}\bigr)
  \circ F^r_{M,\Phi_{\M}(\mathbf X)}
  \circ F(\kappa^{\M}_{\mathbf X,M})
  =\kappa^{\N}_{\mathbf X,F(M)}\circ F^l_{X,M}
\end{equation}
for all $\mathbf X=(X,\sigma)\in\Z(\C)$ and $M\in\M$. In words, $F$ carries
the central action isomorphism of $\M$ to that of $\N$.

\begin{lemma}
  \label{lem:pseudofunctor-Z-1}
  The above construction defines a functor
  $\Z_{\C,\D}:\uBrPic(\C,\D)\to
  \uEqBr(\Z(\C),\Z(\D))$.
\end{lemma}
\begin{proof}
  By the uniqueness in Lemma~\ref{lem:central-action-unique}, it suffices to
  check that the relevant 2-cells satisfy \eqref{eq:def-partial-2cell}. We
  first check independence of the representative. If
  $\theta:F\Rightarrow F'$ is an isomorphism of bimodule functors, then
  \eqref{eq:def-partial-2cell} for $F$ implies \eqref{eq:def-partial-2cell}
  for $F'$, with the same 2-cell, by the naturality of $\kappa^{\M}$ and
  $\kappa^{\N}$ and the compatibility of $\theta$ with the module structures.
  So $\Z_{\C,\D}([F])$ depends only on $[F]$. Next we check compatibility with
  composition. If $G:\cL\to\M$ is another bimodule equivalence, then
  $(FG)^l_{X,L}=F^l_{X,G(L)}\circ F(G^l_{X,L})$ and
  $(FG)^r_{L,Y}=F^r_{G(L),Y}\circ F(G^r_{L,Y})$. Applying $F$ to
  \eqref{eq:def-partial-2cell} for $G$, using the naturality of $F^r$, and
  combining with \eqref{eq:def-partial-2cell} for $F$ at the objects $G(L)$,
  we obtain \eqref{eq:def-partial-2cell} for $FG$ with the 2-cell
  $\Z_{\C,\D}([F])\circ\Z_{\C,\D}([G])$. Finally, identity functors give
  identity 2-cells.
\end{proof}

Next we show that the functors $\Z_{\C,\D}$ are compatible with composition
of 1-cells. As in
\cite[\S5.1]{etingof2010fusion}, the idea is to move a central object through
the factors of a relative Deligne product.

\begin{lemma}
  \label{lem:pseudofunctor-Z-2}
  For a finite tensor category $\C$, let
  $\upsilon_{\C}:\id_{\Z(\C)}\Rightarrow\Phi_{\C}$ be the unique monoidal
  natural isomorphism such that
  \begin{equation}
    \label{eq:unitor}
    \kappa^{\C}_{\mathbf X,V}=\bigl(V\otimes(\upsilon_{\C})_{\mathbf X}\bigr)\circ\sigma_V
    \qquad(\mathbf X=(X,\sigma)\in\Z(\C),\ V\in\C),
  \end{equation}
  where $\C$ is the regular $\C$-$\C$-bimodule category. For 1-cells
  $\C\xrightarrow{\N}\D\xrightarrow{\M}\E$ in $\uBrPic$, there is a unique
  2-cell
  \[
    c_{\M,\N}:\Phi_{\M}\circ\Phi_{\N}\Rightarrow\Phi_{\M\circ\N}
  \]
  in $\uEqBr$ such that
  \begin{equation}
    \label{eq:compositor}
    \kappa^{\M\circ\N}_{\mathbf X,N\boxtimes_{\D}M}
    =\bigl((N\boxtimes_{\D}M)\triangleleft(c_{\M,\N})_{\mathbf X}\bigr)
    \circ\bigl(N\boxtimes_{\D}\kappa^{\M}_{\Phi_{\N}(\mathbf X),M}\bigr)
    \circ b_{N,\Phi_{\N}(\mathbf X),M}
    \circ\bigl(\kappa^{\N}_{\mathbf X,N}\boxtimes_{\D}M\bigr)
  \end{equation}
  for $\mathbf X\in\Z(\C)$, $N\in\N$ and $M\in\M$. With compositor $c$ and
  unitor $\upsilon$, $\Z:\uBrPic\to\uEqBr$ is a pseudofunctor.
\end{lemma}
\begin{proof}
  Put $\mathcal P=\M\circ\N=\N\boxtimes_{\D}\M$ and let
  $\mathbf X\in\Z(\C)$. The central action isomorphisms of $\N$ and $\M$
  and the balancing give isomorphisms
  \begin{align*}
    \mathbf X\triangleright(N\boxtimes_{\D}M)
    &=(\mathbf X\triangleright N)\boxtimes_{\D}M
    \cong(N\triangleleft\Phi_{\N}(\mathbf X))\boxtimes_{\D}M\\
    &\cong N\boxtimes_{\D}(\Phi_{\N}(\mathbf X)\triangleright M)
    \cong(N\boxtimes_{\D}M)\triangleleft\Phi_{\M}\Phi_{\N}(\mathbf X),
  \end{align*}
  natural in $N$ and $M$. They are compatible with the balancing at
  $Y\in\D$. Indeed, by the right $\D$-linearity of $\kappa^{\N}_{\mathbf X}$
  and the left $\D$-linearity of $\kappa^{\M}_{\Phi_{\N}(\mathbf X)}$, the
  two sides of the compatibility differ by the half-braiding of
  $\Phi_{\N}(\mathbf X)$ with $Y$, applied once before and once after the
  balancing. These two applications agree by the naturality of $b$ in its
  middle variable and the balancing axiom. By
  Lemma~\ref{lem:descent}, the displayed isomorphisms are therefore the
  components of a unique natural isomorphism
  $\lambda_{\mathbf X}:\mathbf X\triangleright\id_{\mathcal P}\Rightarrow
  \id_{\mathcal P}\triangleleft\Phi_{\M}\Phi_{\N}(\mathbf X)$.

  We check that the family $\lambda$ is a monoidal natural isomorphism
  $\Forg\circ\Schau_{\mathcal P}\Rightarrow
  a_{\mathcal P}\circ\Forg\circ\Phi_{\M}\Phi_{\N}$. Each of the required
  identities is an identity between natural transformations of right exact
  functors on $\mathcal P$, so by Lemma~\ref{lem:descent} it may be checked
  on objects $N\boxtimes_{\D}M$. There, naturality in $\mathbf X$ and left
  $\C$-linearity come from those of $\kappa^{\N}$, since $\C$ acts on the
  first factor. Monoidality comes from the monoidality of $\kappa^{\N}$ and
  $\kappa^{\M}$ and the balancing axiom. The tensor structures of
  $\Schau_{\N}$, $\Schau_{\M}$ and $\Schau_{\mathcal P}$ involve the
  braidings of the centers, so this step also uses that $\Phi_{\N}$ is braided.
  Now Lemma~\ref{lem:central-action-unique}, applied to $\mathcal P$ and the
  braided functor $\Phi_{\M}\Phi_{\N}$, gives a unique monoidal isomorphism
  $u:\Phi_{\mathcal P}\Rightarrow\Phi_{\M}\Phi_{\N}$ with
  $\lambda=(\id\triangleleft u)\circ\kappa^{\mathcal P}$. We put
  $c_{\M,\N}=u^{-1}$. Unwinding the definition of the family, this reads as
  \eqref{eq:compositor}. Conversely, \eqref{eq:compositor} determines
  $c_{\M,\N}$ by Lemma~\ref{lem:descent} and the uniqueness in
  Lemma~\ref{lem:central-action-unique}.

  The pseudofunctor axioms follow from the same uniqueness. Let
  $F:\N\to\N'$ and $G:\M\to\M'$ be bimodule equivalences. The naturality
  of $c$ with respect to $\Z([G])\hcomp\Z([F])$ and $\Z([F\boxtimes_{\D}G])$
  is an identity checked on objects $N\boxtimes_{\D}M$, where it follows
  from \eqref{eq:def-partial-2cell} for $F$ and $G$ and the naturality of
  $b$.

  For associativity, consider 1-cells
  $\B\xrightarrow{\cL}\C\xrightarrow{\N}\D\xrightarrow{\M}\E$. The two 2-cells
  $\Phi_{\M}\Phi_{\N}\Phi_{\cL}\Rightarrow\Phi_{\M\circ(\N\circ\cL)}$ in the
  associativity axiom (one of them composed with $\Z$ of the associativity
  equivalence) both satisfy the analogue of \eqref{eq:compositor} on objects
  $L\boxtimes_{\C}N\boxtimes_{\D}M$: the composite on the right moves
  $\mathbf X$ through the three factors in turn. Here we use
  \eqref{eq:def-partial-2cell} for the associativity equivalence, which is
  the identity on these objects. By the three-factor
  case of Lemma~\ref{lem:descent} and Lemma~\ref{lem:central-action-unique},
  the two 2-cells coincide.

  For the unit axioms, consider the regular bimodule category $\C$. The
  half-braidings $\sigma_V:X\otimes V\to V\otimes X$ of the objects
  $\mathbf X=(X,\sigma)$ form a monoidal natural isomorphism
  $\Forg\circ\Schau_{\C}\Rightarrow a_{\C}\circ\Forg$, and
  Lemma~\ref{lem:central-action-unique}, applied with $T=\id_{\Z(\C)}$, gives
  $\upsilon_{\C}^{-1}$, which proves the existence and uniqueness of
  $\upsilon_{\C}$ in \eqref{eq:unitor}. Under the unit equivalences
  $\C\boxtimes_{\C}\N\simeq\N\simeq\N\boxtimes_{\D}\D$, the composite that
  moves $\mathbf X$ first through $\C$ and then through $\N$, or first
  through $\N$ and then through $\D$, is $\kappa^{\N}_{\mathbf X}$, by the
  naturality of $\kappa^{\N}$ in $\mathbf X$ and left, or right, linearity of $\kappa^{\N}$. This proves the unit axioms.
\end{proof}

The pseudofunctor is compatible with adjoint algebras.

\begin{lemma}
  \label{lem:adj-alg-and-partial-M}
  For a 1-cell $\mathcal{M} : \mathcal{C} \to \mathcal{D}$ and an indecomposable exact left
  $\mathcal{D}$-module category $\mathcal{N}$, there is an isomorphism
  \begin{equation*}
    \Phi_{\mathcal{M}}^{-1}(\bbA_{\mathcal{N}})
    \cong \bbA_{\mathcal{M} \boxtimes_{\mathcal{D}} \mathcal{N}}
  \end{equation*}
  of algebras in $\mathcal{Z}(\mathcal{C})$.
\end{lemma}
\begin{proof}
  We first consider the case where $\mathcal{N} = \mathcal{D}$. Here
  $\bbA_{\D}\cong\Induc_{\D}(\unit_{\D})$ as algebras, by
  \eqref{eq:F-tilde-and-induction} for the equivalence $a_{\D}:\D\to\D^*_{\D}$
  and the isomorphism
  $\Schau_{\D}\cong\widetilde{a_{\D}}\circ\Phi_{\D}\cong\widetilde{a_{\D}}$
  given by the unitor $\upsilon_{\D}$ of Lemma~\ref{lem:pseudofunctor-Z-2}.
  By \eqref{eq:F-tilde-and-induction} with
  $F = a_{\mathcal{M}}$, we have
  isomorphisms of algebras in $\mathcal{Z}(\mathcal{C})$ as follows:
  \begin{gather*}
    \Phi_{\mathcal{M}}^{-1}(\bbA_{\mathcal{D}})
    \cong (\Schau_{\mathcal{M}}^{-1} \circ \widetilde{a_{\mathcal{M}}} \circ \Induc_{\mathcal{D}})(\unit_{\mathcal{D}})
    \cong (\Schau_{\mathcal{M}}^{-1} \circ \Induc_{\mathcal{C}_{\mathcal{M}}^*} \circ a_{\mathcal{M}})(\unit_{\mathcal{D}}) \\
    \cong (\Schau_{\mathcal{M}}^{-1} \circ \Induc_{\mathcal{C}_{\mathcal{M}}^*})(\unit_{\mathcal{C}_{\mathcal{M}}^*})
    \cong \bbA_{\mathcal{M}}
    \cong \bbA_{\mathcal{M} \boxtimes_{\mathcal{D}} \mathcal{D}},
  \end{gather*}
  where the last isomorphism is \eqref{eq:def-chi-F} for the unit
  equivalence $\M\boxtimes_{\D}\D\simeq\M$. Now we consider the general case.
  The left $\D$-module category $\N$ is an invertible
  $\D$-$\E$-bimodule category for $\E=\D_{\N}^*$ by
  Proposition~\ref{prop:invertibility-criterion}, since its action functor
  is the identity of $\D_{\N}^*$. Hence $\N\circ\M=\M\boxtimes_{\D}\N$ is an invertible $\C$-$\E$-bimodule category, so it is indecomposable exact as a left $\C$-module category and $\bbA_{\M\boxtimes_{\D}\N}$ is defined (Proposition~\ref{prop:invertibility-criterion}). By the first case, applied to $\N$ and to
  $\N\circ\M=\M\boxtimes_{\D}\N$, and by Lemma~\ref{lem:pseudofunctor-Z-2},
  we have isomorphisms
  \begin{equation*}
    \Phi_{\mathcal{M}}^{-1}(\bbA_{\mathcal{N}})
    \cong \Phi_{\mathcal{M}}^{-1}\Phi_{\mathcal{N}}^{-1}(\bbA_{\mathcal{E}})
    \cong \Phi_{\mathcal{M} \boxtimes_{\mathcal{D}} \mathcal{N}}^{-1}(\bbA_{\mathcal{E}})
    \cong \bbA_{\mathcal{M} \boxtimes_{\mathcal{D}} \mathcal{N}}
  \end{equation*}
  of algebras in $\mathcal{Z}(\mathcal{C})$.
\end{proof}

\subsection{Bijectivity on 1-cells up to isomorphism}
\label{subsec:bijectivity-on-1-cells}

We denote by $\BrPic$ the groupoid obtained from the 2-groupoid $\uBrPic$ by
identifying isomorphic 1-cells and truncating 2-cells. A similar notation is
used for $\uEqBr$. We establish a truncated version of
Theorem~\ref{thm:BrPic}, that is, the bijectivity of the map
\[
  \Phi:\BrPic(\C,\D)\to \EqBr(\Z(\C),\Z(\D)),
  \qquad
  [\M]\mapsto [\Phi_{\M}]
\]
for 0-cells $\C$ and $\D$ of $\uBrPic$. It is well defined by
Lemma~\ref{lem:pseudofunctor-Z-1}. We construct a map $\Psi$ in the
opposite direction from the classification of Lagrangian algebras, show
that it is well defined (Lemma~\ref{lem:psi-well-defined}), identify the
induced functor $\widetilde a$ of its right action
(Lemma~\ref{lem:tilde-right-action}), and conclude that $\Phi$ and $\Psi$ are
mutually inverse (Lemma~\ref{lem:BrPic-truncated}).

We construct the inverse of $\Phi$ as follows. Let $T:\Z(\C)\to \Z(\D)$ be a braided tensor equivalence. Then
$T^{-1}(\bbA_{\D})$ is a Lagrangian algebra in $\Z(\C)$, since $T^{-1}$ is a
braided tensor equivalence and $\bbA_{\D}$ is Lagrangian in $\Z(\D)$.
Thus, by
Theorem~\ref{thm:classification-Lagrangian}, there exists an indecomposable exact left
$\C$-module category $\M$ such that
$T^{-1}(\bbA_{\D})\cong \bbA_{\M}$ as algebras in $\Z(\C)$. We fix an
isomorphism $f:\bbA_{\M}\to T^{-1}(\bbA_{\D})$ of algebras in $\Z(\C)$ and
make $\M$ a $\C$-$\D$-bimodule category, invertible by
Proposition~\ref{prop:invertibility-criterion}, whose action functor is the
tensor equivalence
\begin{equation}
  \label{eq:def-psi-right-action}
  \D
  \xrightarrow{\ \Kappa_{\D}\ }
  \Z(\D)_{\bbA_{\D}}
  \xrightarrow{\ T^{-1}\ }
  \Z(\C)_{T^{-1}(\bbA_{\D})}
  \xrightarrow{\ f^{\star}\ }
  \Z(\C)_{\bbA_{\M}}
  \xrightarrow{\ \Kappa_{\M}^{-1}\ }
  \C_{\M}^*,
\end{equation}
where the tensor equivalence
$\Z(\D)_{\bbA_{\D}}\to \Z(\C)_{T^{-1}(\bbA_{\D})}$ induced by $T^{-1}$ is
denoted by $T^{-1}$ by abuse of notation, and $\Kappa_{\D}$ is the
equivalence \eqref{eq:adjoint-alg-module-in-ZC} for the surjective central
functor $\Forg_{\D}:\Z(\D)\to\D$ in place of $\Forg\circ\Schau_{\M}$, whose
right adjoint sends $\unit$ to $\Induc_{\D}(\unit)\cong\bbA_{\D}$ (see the
proof of Lemma~\ref{lem:adj-alg-and-partial-M}). We then
define
\[
  \Psi:\EqBr(\Z(\C),\Z(\D))\to \BrPic(\C,\D),
  \qquad
  [T]\mapsto [\M].
\]

\begin{lemma}
  \label{lem:psi-well-defined}
  The map $\Psi$ is well-defined.
\end{lemma}
\begin{proof}
  Let $T'$ be a braided tensor equivalence isomorphic to $T$, and let
  $\gamma:T'\to T$ be an isomorphism of braided tensor functors. Choose
  $\M'$ and $f':\bbA_{\M'}\to T'{}^{-1}(\bbA_{\D})$ as above. The isomorphism
  $\gamma$ induces an algebra isomorphism
  $g:T^{-1}(\bbA_{\D})\to T'{}^{-1}(\bbA_{\D})$. Then
  $h:=(f')^{-1}gf$ is an algebra isomorphism
  $\bbA_{\M}\to \bbA_{\M'}$. By
  Lemma~\ref{lem:isomorphisms-between-Lagrangian}(1), $\M\simeq\M'$, and then
  by part (2) of that lemma, $h=\chi_{\alpha}$ for some
  $\C$-module equivalence $\alpha:\M\to \M'$.

  The comparison of the right $\D$-actions is encoded by
  \[
    \begin{tikzcd}[column sep=20pt]
      \D \arrow[r, "{\Kappa_{\D}}"] \arrow[d, equal]
      & \Z(\D)_{\bbA_{\D}}
      \arrow[r, "{T^{-1}}"] \arrow[d, equal]
      & \Z(\C)_{T^{-1}(\bbA_{\D})}
      \arrow[r, "{f^{\star}}"] \arrow[d, "{(g^{-1})^{\star}}"]
      & \Z(\C)_{\bbA_{\M}}
      \arrow[r, "{\Kappa_{\M}^{-1}}"] \arrow[d, "{(h^{-1})^{\star}}"]
      & \C_{\M}^*
      \arrow[d, "{\inn_{\alpha}}"] \\
      \D \arrow[r, "{\Kappa_{\D}}"]
      & \Z(\D)_{\bbA_{\D}}
      \arrow[r, "{T'{}^{-1}}"]
      & \Z(\C)_{T'{}^{-1}(\bbA_{\D})}
      \arrow[r, "{f'{}^{\star}}"]
      & \Z(\C)_{\bbA_{\M'}}
      \arrow[r, "{\Kappa_{\M'}^{-1}}"]
      & \C_{\M'}^*.
    \end{tikzcd}
  \]
  The right square commutes up to isomorphism by
  Lemma~\ref{lem:functors-induced-by-isomorphisms}, the third square
  commutes because $f\circ h^{-1}=g^{-1}\circ f'$, and the second square commutes up
  to isomorphism by the functoriality of transport along $\gamma$. Hence
  $a_{\M'}\cong\inn_{\alpha}\circ a_{\M}$ for the two right $\D$-actions,
  so $\alpha$ is an equivalence of the resulting $\C$-$\D$-bimodule
  categories by Lemma~\ref{lem:bimodule-criterion}.
  Taking $T'=T$, this also shows that $[\M]$ does not depend on the choice of
  $\M$ and $f$.
\end{proof}

To prove that $\Phi$ and $\Psi$ are mutually inverse, we first relate the tensor equivalence defined by
\eqref{eq:def-psi-right-action} to the Schauenburg equivalence.

\begin{lemma}
  \label{lem:tilde-right-action}
  For the tensor equivalence $a:\D\to \C_{\M}^*$ defined by
  \eqref{eq:def-psi-right-action}, there is an isomorphism
  \[
    \widetilde{a}\cong \Schau_{\M}\circ T^{-1}
  \]
  of tensor functors.
\end{lemma}
\begin{proof}
  Both $\widetilde a$ and $\Schau_{\M}\circ T^{-1}$ are braided tensor
  functors $\Z(\D)\to\Z(\C_{\M}^*)$. By
  Lemma~\ref{lem:braided-lift-unique}, it suffices to show that their
  underlying functors $a\circ\Forg_{\D}$ and
  $\Forg\circ\Schau_{\M}\circ T^{-1}$ are isomorphic as tensor functors and
  surjective. We compare them after applying the tensor equivalence
  $\Kappa_{\M}$. Recall the free module functor $F_A$ of
  \eqref{eq:adjoint-alg-module-in-ZC}. It
  is surjective, since every module is a quotient of a free one. By
  \eqref{eq:adjoint-alg-module-in-ZC} and its analogue for $\Kappa_{\D}$, there are isomorphisms
  $\Kappa_{\M}\circ\Forg\circ\Schau_{\M}\cong F_{\bbA_{\M}}$ and
  $\Kappa_{\D}\circ\Forg_{\D}\cong F_{\bbA_{\D}}$ of tensor functors. Hence
  \[
    \Kappa_{\M}\circ a\circ\Forg_{\D}
    \cong f^{\star}\circ T^{-1}\circ F_{\bbA_{\D}}
    \cong f^{\star}\circ F_{T^{-1}(\bbA_{\D})}\circ T^{-1}
    \cong F_{\bbA_{\M}}\circ T^{-1}
    \cong \Kappa_{\M}\circ\Forg\circ\Schau_{\M}\circ T^{-1}
  \]
  as tensor functors. The first isomorphism combines
  \eqref{eq:def-psi-right-action} with
  $\Kappa_{\D}\circ\Forg_{\D}\cong F_{\bbA_{\D}}$, the second is the tensor
  structure of $T^{-1}$, and the third is given by
  $\id\otimes f:\mathbf W\otimes\bbA_{\M}\to
  f^{\star}(\mathbf W\otimes T^{-1}(\bbA_{\D}))$. The last is
  \eqref{eq:adjoint-alg-module-in-ZC}.
\end{proof}

\begin{lemma}
  \label{lem:BrPic-truncated}
  The maps $\Phi$ and $\Psi$ are mutually inverse.
\end{lemma}
\begin{proof}
  Let $T:\Z(\C)\to \Z(\D)$ be a braided tensor equivalence, and let $\M$ be
  the bimodule category representing $\Psi([T])$. By
  Lemma~\ref{lem:tilde-right-action}, the right action
  $a:\D\to \C_{\M}^*$ satisfies
  $\widetilde{a}\cong \Schau_{\M}\circ T^{-1}$. Thus
  \[
    \Phi_{\M}
    =
    \widetilde{a^{-1}}\circ \Schau_{\M}
    \cong
    T\circ \Schau_{\M}^{-1}\circ \Schau_{\M}
    \cong T,
  \]
  so $\Phi\Psi([T])=[T]$.

  Conversely, let $\M$ be an invertible $\C$-$\D$-bimodule category and set
  $T=\Phi_{\M}$. By Lemma~\ref{lem:adj-alg-and-partial-M},
  $T^{-1}(\bbA_{\D})\cong \bbA_{\M}$. Since $\Psi$ is well defined, we may
  compute $\Psi([T])$ using the left $\C$-module category $\M$ itself and an
  isomorphism $f:\bbA_{\M}\to T^{-1}(\bbA_{\D})$. Let $\M'$ be $\M$ equipped
  with the right $\D$-action constructed in
  \eqref{eq:def-psi-right-action}. Lemma~\ref{lem:tilde-right-action} gives
  $\widetilde{a_{\M'}}\cong\Schau_{\M}\circ\Phi_{\M}^{-1}\cong
  \widetilde{a_{\M}}$. By Lemma~\ref{lem:tilde-isomorphic}, this implies
  $a_{\M'}\cong \inn_{\alpha}\circ a_{\M}$ for some invertible object
  $\alpha\in \C_{\M}^*$. Since the tensor product of $\C_{\M}^*$ is reversed
  composition, $\inn_{\alpha}(G)=\alpha^*\circ G\circ\alpha$, that is,
  $\inn_{\alpha}$ is conjugation by the module equivalence $\alpha^*$, with
  quasi-inverse $\alpha$. By Lemma~\ref{lem:bimodule-criterion}, the module
  equivalence $\alpha^*:\M\to\M'$ is an equivalence of bimodule categories.
  Thus $\M'\simeq\M$ as bimodule categories and
  $\Psi\Phi([\M])=[\M]$.
\end{proof}

\begin{remark}\label{rem:DN-gap}
  For $\C=\D$, Lemmas~\ref{lem:pseudofunctor-Z-2} and~\ref{lem:BrPic-truncated} give the isomorphism
  $\BrPic(\C)\cong\Aut_{\otimes}^{\br}(\Z(\C))$ of
  \cite[Theorem~4.1]{davydov2013picard}, whose proof constructs the inverse
  map differently. For a braided autoequivalence $\alpha$ of $\Z(\C)$, it
  puts $L_{\alpha}=\alpha^{-1}(\Induc_{\C}(\unit))$, asserts that the algebra
  $\Forg_{\C}(L_{\alpha})$ in $\C$ is exact, and takes modules over an
  indecomposable summand. This assertion fails in general. Their argument shows that $\C_{\Forg_{\C}(L_{\alpha})}$ is exact over $\Fun_{\Z(\C)}(\C,\C)$ \cite[Lemma~3.2 and Remark~2.8]{davydov2013picard} and concludes that it is exact over $\C$. This last step fails, as follows. Let $\M$
  represent $\Psi([\alpha])$, so that $L_{\alpha}\cong\bbA_{\M}$ by
  construction. By \eqref{eq:WY25-Prop-4-2}, the left $\C$-module category
  $\C_{\Forg_{\C}(L_{\alpha})}$ is equivalent to $\M\boxtimes\M^{\op}$, with
  $\C$ acting on the first factor. Let $P$ be a projective cover of $\unit$.
  For nonzero $M\in\M$, the epimorphism $P\triangleright M\to\unit\triangleright
  M\cong M$ shows $P\triangleright M\neq0$, so for $N\in\M^{\op}$ the object
  $(P\triangleright M)\boxtimes N$ is projective only if $N$ is. Conversely,
  $\M\boxtimes\M^{\op}$ is semisimple if $\M$ is. Hence
  $\Forg_{\C}(L_{\alpha})$ is exact if and only if $\M$ is semisimple.
  For $\alpha=\id$ we have $\M\simeq\C$, so $\Forg_{\C}(\Induc_{\C}(\unit))$
  is exact only when $\C$ is semisimple. Our map $\Psi$ uses
  Theorem~\ref{thm:classification-Lagrangian} instead, which needs only
  that $L_{\alpha}$ is Lagrangian. By Remark~\ref{rem:Lagrangian}, the role
  of the summand is played by $\C_A$ for a simple quotient $A$ of
  $\Forg_{\C}(L_{\alpha})$.
\end{remark}

\subsection{Full faithfulness of $\Z$}
\label{subsec:fully-faith-Z}

We now complete the proof of Theorem~\ref{thm:BrPic}. By
Lemma~\ref{lem:BrPic-truncated}, the functor $\Z_{\C,\D}$ between groupoids
is essentially surjective, and $\Phi_{\M}\cong\Phi_{\N}$ implies
$\M\simeq\N$. Since Hom sets between isomorphic objects are torsors under the automorphism groups, compatibly with $\Z_{\C,\D}$,  it is an equivalence once the induced map
\[
  \Z_{\M,\M}:\Aut(\M)\to\Aut(\Phi_{\M})
\]
is bijective for every invertible $\C$-$\D$-bimodule category $\M$, where
$\Aut(\M)$ is the group of isomorphism classes of bimodule
autoequivalences of $\M$ and $\Aut(\Phi_{\M})$ is the group of monoidal
automorphisms of $\Phi_{\M}$.

\begin{lemma}
  \label{lem:Aut-M-invertible}
  Let $\M$ be an invertible $\C$-$\D$-bimodule category. Then
  \[
    \Inv(\Z(\C))\to\Aut(\M),
    \qquad
    [\boldsymbol\beta]\mapsto[\boldsymbol\beta\triangleright\id_{\M}],
  \]
  is a bijection, where $\boldsymbol\beta\triangleright\id_{\M}$ carries the
  bimodule structure described after \eqref{eq:Phi-central-action}.
\end{lemma}
\begin{proof}
  As in the first item after \eqref{eq:Phi-central-action}, a right
  $\D$-module structure on a left $\C$-module endofunctor $G$ of $\M$ is the
  same as a half-braiding of $G\in\C_{\M}^*$ against the objects
  $a_{\M}(Y)$, $Y\in\D$. Since $a_{\M}$ is an equivalence, such a
  half-braiding extends uniquely to all of $\C_{\M}^*$. This identifies the
  $\C$-$\D$-bimodule endofunctors of $\M$, and their natural transformations,
  with $\Z(\C_{\M}^*)$. Composing with the equivalence $\Schau_{\M}$, we see
  that every bimodule endofunctor of $\M$ is isomorphic to
  $\mathbf X\triangleright\id_{\M}$ for some $\mathbf X\in\Z(\C)$
  unique up to isomorphism. If the endofunctor is an autoequivalence, then the
  underlying object of $\Schau_{\M}(\mathbf X)$ is invertible in $\C_{\M}^*$, so
  $\Schau_{\M}(\mathbf X)$ is invertible in $\Z(\C_{\M}^*)$ (an inverse inherits a
  half-braiding), and hence $\mathbf X$ is invertible. Conversely,
  $\boldsymbol\beta\triangleright\id_{\M}$ is an autoequivalence for invertible
  $\boldsymbol\beta$.
\end{proof}

Next, let $c$ denote the braiding of $\Z(\C)$. For an invertible object
$\boldsymbol\beta$ and an object $\mathbf X$ of $\Z(\C)$, define the monodromy
$\Omega_{\boldsymbol\beta}(\mathbf X)\in\Aut_{\Z(\C)}(\mathbf X)$ by
\[
  \Omega_{\boldsymbol\beta}(\mathbf X)\otimes\id_{\boldsymbol\beta}
  =c_{\boldsymbol\beta,\mathbf X}\circ c_{\mathbf X,\boldsymbol\beta}.
\]
Then $\Omega_{\boldsymbol\beta}$ is a monoidal automorphism of
$\id_{\Z(\C)}$, and the monodromy map
$[\boldsymbol\beta]\mapsto\Omega_{\boldsymbol\beta}$ is a bijection from
$\Inv(\Z(\C))$ to the group of monoidal automorphisms of $\id_{\Z(\C)}$
\cite[Theorem~5.3 and equation~(5.7), with $F=G=\id_{\C}$]{shimizu2017ribbon}.

\begin{lemma}
  \label{lem:Z-of-invertible-object}
  Let $\M$ be an invertible $\C$-$\D$-bimodule category, let
  $\boldsymbol\beta=(\beta,\sigma^{\beta})$ be an invertible object of
  $\Z(\C)$, and put $F=\boldsymbol\beta\triangleright\id_{\M}$. Then
  \begin{equation}
    \label{eq:Z-of-invertible-object}
    \Z_{\C,\D}([F])
    =\id_{\Phi_{\M}}\hcomp\Omega_{\boldsymbol\beta},
    \qquad\text{that is,}\qquad
    \Z_{\C,\D}([F])_{\mathbf X}=\Phi_{\M}(\Omega_{\boldsymbol\beta}(\mathbf X)).
  \end{equation}
\end{lemma}
\begin{proof}
  Let $\mathbf X=(X,\sigma)\in\Z(\C)$ and $M\in\M$. The left constraint
  $F^l_{X,M}:\beta\triangleright(X\triangleright M)\to
  X\triangleright(\beta\triangleright M)$ is induced by
  $\sigma^{\beta}_X=c_{\boldsymbol\beta,\mathbf X}$, and $F^r$ is the middle
  associator. The left $\C$-linearity of $\kappa^{\M}_{\mathbf X}$, applied
  to the object $\beta\in\C$, says that
  \[
    F^r_{M,\Phi_{\M}(\mathbf X)}\circ F(\kappa^{\M}_{\mathbf X,M})
    =\kappa^{\M}_{\mathbf X,\beta\triangleright M}\circ s_{\beta,M},
  \]
  where $s_{\beta,M}:\beta\triangleright(X\triangleright M)\to
  X\triangleright(\beta\triangleright M)$, the module structure of
  $\mathbf X\triangleright\id_{\M}$ at $\beta$ in the notation of
  \S\ref{subsec:Schauenburg-equiv}, is induced by
  $\sigma_{\beta}^{-1}=c_{\mathbf X,\boldsymbol\beta}^{-1}$. Thus
  \eqref{eq:def-partial-2cell} for $u=\Z_{\C,\D}([F])$ becomes
  \[
    \bigl((\beta\triangleright M)\triangleleft u_{\mathbf X}\bigr)
    \circ\kappa^{\M}_{\mathbf X,\beta\triangleright M}
    =\kappa^{\M}_{\mathbf X,\beta\triangleright M}\circ F^l_{X,M}\circ s_{\beta,M}^{-1}.
  \]
  The automorphism $F^l_{X,M}\circ s_{\beta,M}^{-1}$ of
  $X\triangleright(\beta\triangleright M)$ is induced by
  $c_{\boldsymbol\beta,\mathbf X}\circ c_{\mathbf X,\boldsymbol\beta}$, which
  is $\Omega_{\boldsymbol\beta}(\mathbf X)\otimes\id_{\boldsymbol\beta}$ by
  definition. By the naturality of $\kappa^{\M}$ in
  $\mathbf X$, the right-hand side is therefore
  $\bigl((\beta\triangleright M)\triangleleft
  \Phi_{\M}(\Omega_{\boldsymbol\beta}(\mathbf X))\bigr)\circ
  \kappa^{\M}_{\mathbf X,\beta\triangleright M}$. Since
  $\beta\triangleright(-)$ is an equivalence and $a_{\M}$ is faithful, this
  proves \eqref{eq:Z-of-invertible-object}.
\end{proof}

\begin{proof}[Proof of Theorem~\ref{thm:BrPic}]
  It remains to show that $\Z_{\M,\M}$ is bijective. By
  Lemma~\ref{lem:Z-of-invertible-object}, the composite of $\Z_{\M,\M}$ with
  the bijection of Lemma~\ref{lem:Aut-M-invertible} is
  $[\boldsymbol\beta]\mapsto\id_{\Phi_{\M}}\hcomp\Omega_{\boldsymbol\beta}$.
  This is a bijection onto $\Aut(\Phi_{\M})$: the monodromy map is
  bijective, and whiskering with the equivalence $\Phi_{\M}$ is a bijection
  from monoidal automorphisms of $\id_{\Z(\C)}$ to those of $\Phi_{\M}$.
  Hence $\Z_{\M,\M}$ is bijective.
\end{proof}

\subsection{Equivalence of dual categories}
\label{subsec:BrPic-applications}
We apply Theorem~\ref{thm:BrPic} to the question of when two module
categories have tensor-equivalent dual categories. For a finite
tensor category $\C$, put
\begin{equation*}
  \bbX_{\C}
  :=
  \left\{
    \begin{array}{c}
      \text{nonzero indecomposable exact left}\\[-2pt]
      \text{$\C$-module categories}
    \end{array}
  \right\}\big/\simeq_{\C}.
\end{equation*}

\begin{definition}\label{def:monoidal-Morita-equivalence}
  Let $\M$ and $\N$ represent elements of $\bbX_{\C}$. We say that
  $\M$ and $\N$ are \emph{dual equivalent}, and write
  \[
    \M\sim_{\mathrm{de}}\N,
  \]
  if their dual finite tensor categories are tensor equivalent: $\C_{\M}^*\simeq_{\otimes}\C_{\N}^*$.
\end{definition}

In the 2-category of indecomposable exact $\C$-module categories, the
endomorphism category of $\M$ is
$\operatorname{Fun}_{\C}(\M,\M)=(\C_{\M}^*)^{\rev}$, with composition as
tensor product. Since reversing tensor products preserves tensor
equivalence, $\M\sim_{\mathrm{de}}\N$ means that $\M$ and $\N$ have
tensor-equivalent endomorphism categories. Theorem~\ref{thm:MME} shows that
this happens if and only if $\M$ and $\N$ lie in the same orbit of the action
of invertible $\C$-bimodule categories.

The Brauer--Picard group $\BrPic(\C)$ is the group of equivalence classes of
invertible $\C$-bimodule categories, with multiplication
$[\M][\N]=[\M\boxtimes_{\C}\N]$ \cite[\S8.27.2]{etingof2016tensor}. Since
$\M\boxtimes_{\C}\N=\N\circ\M$, Lemma~\ref{lem:pseudofunctor-Z-2} gives
$\Phi_{\M\boxtimes_{\C}\N}\cong\Phi_{\N}\circ\Phi_{\M}$. Thus the assignment $[\mathcal{P}] \mapsto [\Phi_\mathcal P]$ reverses the group multiplication. Consequently, the assignment $[\mathcal P] \mapsto [\Phi_\mathcal P^{-1}]$ is a group homomorphism. By Lemma~\ref{lem:BrPic-truncated}, it is in fact a group isomorphism
\begin{equation}
  \label{eq:BrPic-center-iso}
  \BrPic(\C)\xrightarrow{\ \sim\ }\Aut_{\otimes}^{\mathrm{br}}(\Z(\C)),
  \qquad
  [\mathcal P]\mapsto[\Phi_{\mathcal P}^{-1}].
\end{equation}

The isomorphism \eqref{eq:BrPic-center-iso} agrees with the one constructed in \cite[\S5.1]{etingof2010fusion} and \cite[Theorem~4.1]{davydov2013picard}. To see this, recall that the braided autoequivalence associated there with an invertible $\C$-bimodule category $\mathcal{P}$ sends $\mathbf{Z} \in \Z(\C)$ to the object whose left action on $\mathcal{P}$ agrees with the right action of $\mathbb{Z}$. By \eqref{eq:Phi-central-action}, this object is precisely $\Phi_{\mathcal{P}}^{-1}(\mathbf{Z})$.

The group $\BrPic(\C)$ acts on $\bbX_{\C}$ by
\[
  [\mathcal{P}]\cdot[\M]
  :=
  [\mathcal{P}\boxtimes_{\C}\M],
\]
which is well defined, since $\mathcal P\boxtimes_{\C}\M$ is an invertible
$\C$-$\C_{\M}^*$-bimodule category, hence indecomposable and exact as a left
$\C$-module category
(Proposition~\ref{prop:invertibility-criterion}).
Similarly,
$\Aut_{\otimes}^{\mathrm{br}}(\Z(\C))$ acts on the set
$\mathrm{Lag}(\Z(\C))$ of isomorphism classes of Lagrangian algebras by
$[F]\cdot[\bbL]:=[F(\bbL)]$. By Lemma~\ref{lem:adj-alg-and-partial-M},
$\bbA_{\mathcal P\boxtimes_{\C}\M}\cong\Phi_{\mathcal P}^{-1}(\bbA_{\M})$, so
the bijection of Theorem~\ref{thm:classification-Lagrangian} is
equivariant along the isomorphism \eqref{eq:BrPic-center-iso}.

\begin{theorem}\label{thm:MME}
  Let $\C$ be a finite tensor category. The dual equivalence classes in $\bbX_{\C}$ are
  exactly the $\BrPic(\C)$-orbits. The assignment $\M\mapsto\C_{\M}^*$
  induces a bijection from the set of these orbits to the set of
  tensor-equivalence classes of finite tensor categories Morita equivalent
  to $\C$, and $\M\mapsto\bbA_{\M}$ induces a bijection from the set of $\BrPic(\C)$-orbits of $\mathbb{X}_{\C}$ to the set
  of $\Aut_{\otimes}^{\mathrm{br}}(\Z(\C))$-orbits in $\mathrm{Lag}(\Z(\C))$.
\end{theorem}
\begin{proof}
  Every finite tensor category Morita equivalent to $\C$ is tensor equivalent
  to $\C_{\M}^*$ for some indecomposable exact left $\C$-module category $\M$ (see \S\ref{subsec:exact-module-categories}). We show that
  the fibers of the map
  \[
    [\M]\longmapsto[\C_{\M}^*]
  \]
  are precisely the $\BrPic(\C)$-orbits.

  If $\N\simeq\mathcal{P}\boxtimes_{\C}\M$ for an invertible
  $\C$-bimodule category $\mathcal{P}$, then
  $\mathcal{P}\boxtimes_{\C}\M=\M\circ\mathcal P$, with the right action
  inherited from $\M$, is an invertible $\C$-$\C_{\M}^*$-bimodule category.
  Hence $\C_{\N}^*\simeq_{\otimes}\C_{\M}^*$.

  Conversely, suppose that
  there is a tensor equivalence $T:\C_{\M}^*\xrightarrow{\sim}\C_{\N}^*$.
  Set
  $\D=\C_{\M}^*$ and regard $\N$ as a $\C$-$\D$-bimodule category
  $\N_T$ by transporting its right $\C_{\N}^*$-action along $T$.
  We recall that $\M$ is an invertible $\C$-$\D$-bimodule category. Let $\overline{\M}$ denote an inverse bimodule category. Then
  \[
    \mathcal{P}:=\N_T\boxtimes_{\D}\overline{\M}
  \]
  is an invertible $\C$-bimodule category, and the associativity and unit
  equivalences of \S\ref{subsubsec:relative-Deligne} give
  \[
    \mathcal{P}\boxtimes_{\C}\M
    \simeq
    \N_T\boxtimes_{\D}
    (\overline{\M}\boxtimes_{\C}\M)
    \simeq
    \N_T\boxtimes_{\D}\D
    \simeq
    \N
  \]
  as left $\C$-module categories. Thus the fibers of
  $[\M]\mapsto[\C_{\M}^*]$ are the dual equivalence classes and the
  $\BrPic(\C)$-orbits. This proves the first two statements.

  Theorem~\ref{thm:classification-Lagrangian} identifies $\bbX_{\C}$
  with the set of isomorphism classes of Lagrangian algebras, and this
  identification is equivariant along
  the isomorphism \eqref{eq:BrPic-center-iso}, as shown before the theorem.
  Hence it induces the last bijection.
\end{proof}

\begin{remark}\label{rem:mme-classes-not-finite}
  If $\C$ is a fusion category over a field of characteristic zero, then $\bbX_{\C}$ is finite \cite[Corollary~9.1.6(ii)]{etingof2016tensor}. For finite tensor categories, even the set of dual equivalence classes, can be infinite. An example is given as follows.
  Etingof and Gelaki construct three families of
  $32$-dimensional triangular Hopf algebras. Within each family the members
  are twist equivalent, while the duals of all members together fall into
  infinitely many twist-equivalence classes
  \cite[\S3, Theorem~4.1, and Corollary~4.3]{etingof2002families}.
  So some family $\{H_t\}$ has infinitely many twist classes of duals.
  Put
  $\C=\Rep(H_{t_0})\simeq\Rep(H_t)$. The fiber functor of each $\Rep(H_t)$
  makes $\Vect$ an element of $\bbX_{\C}$ (it is exact, as every object of $\Vect$ is projective) with dual $\Rep(H_t^*)^{\rev}$
  \cite[Example~7.12.26]{etingof2016tensor}, and these
  duals are pairwise tensor-inequivalent for twist-inequivalent $H_t^*$
  \cite[Theorem~6.1]{etingof2002families}.
\end{remark}

Having determined the $\BrPic(\C)$-orbits in $\mathbb{X}_{\C}$ in Theorem~\ref{thm:MME}, we now turn to their stabilizers, which can be described in terms of tensor autoequivalences of the dual categories.
Let
$\M$ be an indecomposable exact left $\C$-module category and put
$\D=\C_{\M}^*$, so that $\M$ is an invertible $\C$-$\D$-bimodule category
with action functor $a_{\M}=\id_{\D}$ and $\Phi_{\M}=\Schau_{\M}$. For a
tensor autoequivalence $\alpha$ of $\D$, let $\M_{\alpha}$ be $\M$ with the
right action $M\triangleleft_{\alpha}Y=\alpha(Y)(M)$. Its action functor is
$\alpha$. The canonical isomorphisms
$\widetilde{\alpha\circ\alpha'}\cong\widetilde{\alpha}\circ\widetilde{\alpha'}$
show that the map
\[
  \iota_{\M}:\Aut_{\otimes}(\C_{\M}^*)\to\Aut_{\otimes}^{\mathrm{br}}(\Z(\C)),
  \qquad
  [\alpha]\mapsto[\Schau_{\M}^{-1}\circ\widetilde{\alpha}\circ\Schau_{\M}],
\]
is a group homomorphism.

\begin{proposition}\label{prop:BrPic-stabilizer}
  Let $\M$ be an indecomposable exact left $\C$-module category, and identify $\BrPic(\C)$
  with $\Aut_{\otimes}^{\mathrm{br}}(\Z(\C))$ via \eqref{eq:BrPic-center-iso}. Then the stabilizer of
  $[\M]\in\bbX_{\C}$ is the image of $\iota_{\M}$, and so is the stabilizer
  of $[\bbA_{\M}]\in\mathrm{Lag}(\Z(\C))$. More precisely, an invertible $\C$-bimodule category $\mathcal{P}$ satisfies $[\mathcal P]\cdot[\M]=[\M]$ if and only if
  $\mathcal P\simeq\M_{\alpha}\boxtimes_{\D}\overline{\M}$ for some
  $\alpha\in\Aut_{\otimes}(\D)$, where $\overline{\M}$ is an inverse bimodule category of $\M$ viewed as an invertible $\C$-$\C_{\M}^*$-bimodule category.
  In this case, we have
  $[\Phi_{\mathcal P}^{-1}]=\iota_{\M}([\alpha])$.
\end{proposition}
\begin{proof}
  For the `if' part, let $\alpha$ be a tensor autoequivalence of $\D$.
  Since the action functor $\alpha$ of $\M_{\alpha}$ is an equivalence,
  $\M_{\alpha}$ is an invertible $\C$-$\D$-bimodule category
  (Proposition~\ref{prop:invertibility-criterion}). Hence
  $\mathcal P_{\alpha}:=\M_{\alpha}\boxtimes_{\D}\overline{\M}$ is an
  invertible $\C$-bimodule category, and the unit and associativity
  equivalences give
  \[
    \mathcal P_{\alpha}\boxtimes_{\C}\M
    \simeq\M_{\alpha}\boxtimes_{\D}(\overline{\M}\boxtimes_{\C}\M)
    \simeq\M_{\alpha}\boxtimes_{\D}\D
    \simeq\M_{\alpha}
  \]
  as $\C$-$\D$-bimodule categories. Since $\M_{\alpha}=\M$ as left
  $\C$-module categories, $[\mathcal P_{\alpha}]\cdot[\M]=[\M]$.

  Conversely, let $\mathcal P$ be an invertible $\C$-bimodule category, and let
  $F:\N\to\M$ be an equivalence of left $\C$-module categories, where
  $\N=\mathcal P\boxtimes_{\C}\M$. With the right action inherited from
  $\M$, $\N$ is an invertible $\C$-$\D$-bimodule category, so its action
  functor $a_{\N}$ is a tensor equivalence. Put
  $\alpha=\inn_F\circ a_{\N}$, a tensor autoequivalence of $\D$. This is the
  action functor of $\M_{\alpha}$, so $F$ is an equivalence
  $\N\to\M_{\alpha}$ of $\C$-$\D$-bimodule categories by
  Lemma~\ref{lem:bimodule-criterion}. Therefore
  \[
    \mathcal P
    \simeq\mathcal P\boxtimes_{\C}(\M\boxtimes_{\D}\overline{\M})
    \simeq\N\boxtimes_{\D}\overline{\M}
    \simeq\M_{\alpha}\boxtimes_{\D}\overline{\M}
    =\mathcal P_{\alpha}.
  \]

  It remains to compute $\Phi_{\mathcal P_{\alpha}}$. The left $\C$-module
  category underlying $\M_{\alpha}$ is $\M$, so
  \eqref{eq:def-partial-M} gives
  $\Phi_{\M_{\alpha}}=\widetilde{\alpha^{-1}}\circ\Schau_{\M}
  \cong\widetilde{\alpha}^{-1}\circ\Schau_{\M}$. Since
  $\mathcal P_{\alpha}\boxtimes_{\C}\M=\M\circ\mathcal P_{\alpha}$ is
  equivalent to $\M_{\alpha}$ as a bimodule category,
  Lemmas~\ref{lem:pseudofunctor-Z-2} and~\ref{lem:pseudofunctor-Z-1} give
  \[
    \Schau_{\M}\circ\Phi_{\mathcal P_{\alpha}}
    \cong\Phi_{\mathcal P_{\alpha}\boxtimes_{\C}\M}
    \cong\Phi_{\M_{\alpha}}
    \cong\widetilde{\alpha}^{-1}\circ\Schau_{\M}.
  \]
  Hence $[\Phi_{\mathcal P_{\alpha}}^{-1}]=\iota_{\M}([\alpha])$. The
  statement about $[\bbA_{\M}]$ follows, since the bijection of
  Theorem~\ref{thm:classification-Lagrangian} is equivariant along
  $\BrPic(\C)\cong\Aut_{\otimes}^{\mathrm{br}}(\Z(\C))$.
\end{proof}

The description of the stabilizer in Proposition \ref{prop:BrPic-stabilizer} yields the following identification with the group of outer tensor autoequivalences of the dual category:

\begin{corollary}\label{cor:BrPic-stabilizer-Out}
  Let $\M$ be an indecomposable exact left $\C$-module category and
  $\D=\C_{\M}^*$. The stabilizer of $[\M]\in\bbX_{\C}$ in $\BrPic(\C)$ is
  isomorphic to
  \[
    \mathrm{Out}_{\otimes}(\D)
    :=\Aut_{\otimes}(\D)/\{[\inn_{\beta}]:\beta\in\Inv(\D)\}.
  \]
\end{corollary}
\begin{proof}
  By Lemma~\ref{lem:tilde-isomorphic} (see also
  \cite[Theorem~5.4]{shimizu2023pivotal}), $\iota_{\M}([\alpha])$ is trivial
  if and only if $\alpha\cong\inn_{\beta}$ for an invertible object $\beta$
  of $\D$. So the kernel of $\iota_{\M}$ is the subgroup of classes of inner
  autoequivalences, and the claim follows from
  Proposition~\ref{prop:BrPic-stabilizer}.
\end{proof}

For a fusion category and $\M=\C$, Corollary~\ref{cor:BrPic-stabilizer-Out}
recovers the description of quasi-trivial
invertible bimodule categories by outer autoequivalences
\cite[\S4.3]{etingof2010fusion}. See also
\cite[Theorem~1.3 and Remark~6.2]{galindo2017tensor}.

\section{Picard groupoid and pullback orbits}
\label{sec:picard-groupoid-orbits}

Let $\B$ be a braided finite tensor category. The center contains two
canonical copies of $\B$, given by the braided embeddings $G^+$ and $G^{-}$ introduced in \S\ref{subsubsec:nondegeneracy}.
In this section,
we use the pseudofunctor $\Z$ constructed in \S\ref{sec:BrPic-2-groupoid} to describe
one-sided bimodule categories, those
whose right action comes from the left action through the braiding. They
form the Picard groupoid $\uPic(\B)$, and Theorem~\ref{thm:Pic} identifies it
with the groupoid of braided autoequivalences of $\Z(\B)$ trivialized on
$G^-(\oB)$ (\S\S\ref{subsec:picard-groupoid}
and~\ref{subsec:picard-center-groupoid}). Such an autoequivalence preserves
$G^+(\B)$, and its restriction there is the boundary map
$\partial:\Pic(\B)\to\Aut_{\otimes}^{\br}(\B)$ of Davydov and Nikshych
(\S\ref{subsec:picard-boundary}).

The group $\Aut_{\otimes}(\B)$ acts on $\bbX_{\B}$ by pullback, and
pullback-equivalent module categories have the same dual category, so
pullback orbits refine the MME classes of
\S\ref{subsec:BrPic-applications}. Since $\partial$ turns pullback
along braided autoequivalences into conjugation, a surjective $\partial$
gives infinitely many $\Aut_{\otimes}^{\br}(\B)$-orbits as soon as
$\Aut_{\otimes}^{\br}(\B)$ has infinitely many conjugacy classes
(Proposition~\ref{prop:inf-orbit-criteria}). For small quantum groups the
Picard locus is a single MME class whose pullback classes are the conjugacy
classes of an algebraic group, so infinitely many pullback classes lie in one
MME class (\S\ref{sec:uqsl2-picard},
Theorem~\ref{thm:uqg-Aut-orbits-infinite}).

\subsection{The Picard groupoid}
\label{subsec:picard-groupoid}

The braiding turns a left $\B$-module category $\M$ into a $\B$-bimodule
category $B(\M)$.
More precisely, $B(\M) = \M$ as a left $\B$-module category.
Its right action is $M\triangleleft X:=X\triangleright M$.
The middle associator is
\[
  X\triangleright(M\triangleleft Y)\cong(X\otimes Y)\triangleright M
  \xrightarrow{\ c_{X,Y}\ }(Y\otimes X)\triangleright M\cong
  (X\triangleright M)\triangleleft Y,
\]
and the right associator is
\[
  (M\triangleleft X)\triangleleft Y\cong(Y\otimes X)\triangleright M
  \xrightarrow{\ c_{X,Y}^{-1}\ }(X\otimes Y)\triangleright M
  =M\triangleleft(X\otimes Y).
\]
The bimodule axioms follow from the hexagon axioms and the naturality of the braiding. The middle associator is the same as the one used in \cite[\S2H, equation~(17)]{davydov2013picard}. We use a slightly different convention for the right associator: whereas \cite[equation (16)]{davydov2013picard} uses $c_{Y,X}$, we use $c_{X,Y}^{-1}$.
Our choice ensures that the action functor $a_{B(\M)}:\B \to \B^*_{\M}$ agrees with $\Forg \circ \Schau_{\M} \circ G^-$ as a tensor functor. Indeed, both functors send $Y\in \B$ to $Y\triangleright(-)$, with $\B$-module structure $c_{X,Y}$ and tensor structure $c_{X,Y}^{-1}$ (see \eqref{eq:Schauenburg-tensor-structure}).
If we adopted the right associator of \cite{davydov2013picard}, the induced tensor structure on $a_{B(\M)}$ would differ by the double braiding.
With our convention, a left $\B$-module functor $T : \M \to \N$ induces a bimodule functor $B(T) : B(\M) \to B(\N)$, whose right module constraint is given by its left module constraint under the identification $M\triangleleft X=X\triangleright M$.

A bimodule category is \emph{one-sided} if it is equivalent to $B(\M)$ for
some left $\B$-module category $\M$. Following
\cite[Remark~2.13]{davydov2013picard}, a \emph{one-sided structure}
on a $\B$-bimodule category $\mathcal P$, with underlying left module
category $\M$, is a natural isomorphism
\[
  d_{M,Y}:M\triangleleft Y\longrightarrow Y\triangleright M
  \qquad(M\in\mathcal P,\ Y\in\B)
\]
whose inverse $d^{-1}$ is a right constraint making the identity functor of
$\M$ a bimodule equivalence $B(\M)\to\mathcal P$. Explicitly, $d$ is compatible
with the middle associators of $B(\M)$ and $\mathcal P$ and with their right
associators. For our convention for $B(\M)$, these are the analogues of
equations (21) and (20) there. So $\mathcal P$ is one-sided if and only if it
admits a one-sided structure, and $B(\M)$ carries the identity one-sided
structure.

We denote by $\uPic(\B)$ the groupoid whose objects are finite left
$\B$-module categories $\M$ such that $B(\M)$ is invertible (they are
indecomposable and exact by Proposition~\ref{prop:invertibility-criterion}),
and whose morphisms are isomorphism classes of left $\B$-module equivalences.
For $\M,\N\in\uPic(\B)$, the balancing isomorphism
\[
  M\boxtimes_{\B}(X\triangleright N)\cong(X\triangleright M)\boxtimes_{\B}N
\]
is a right constraint for the identity functor
$B(\M)\boxtimes_{\B}B(\N)\to B(\M\boxtimes_{\B}\N)$, which is therefore a
bimodule equivalence (compare \cite[\S2H]{davydov2013picard}). The
compatibility with the right associators uses the balancing axiom for
$Y'\otimes Y$, by which the factors $c^{-1}_{Y,Y'}$ cancel. So
$\M\boxtimes_{\B}\N$ lies again in $\uPic(\B)$, and
the relative Deligne product makes $\uPic(\B)$ a monoidal groupoid. The unit
is $\B$: the identity functor with right constraint
$c_{M,Y}^{-1}:Y\otimes M\to M\otimes Y$ is a bimodule equivalence
$B(\B)\to\B$, compatible with the right associator $c_{X,Y}^{-1}$ of $B(\B)$
by the hexagon axioms. We call this monoidal groupoid the \emph{Picard
groupoid} of $\B$. Its monoid of isomorphism classes of objects is a group,
since every object is invertible: by Theorem~\ref{thm:Pic}, $\uPic(\B)$ is
monoidally equivalent to a monoidal groupoid whose objects are invertible
(\S\ref{subsec:picard-center-groupoid}). We call it the \emph{Picard group} $\Pic(\B)$. Forgetting the right action
identifies $\Pic(\B)$ with a subset of $\bbX_{\B}$, since objects of
$\uPic(\B)$ are indecomposable exact and its morphisms are the left module
equivalences.

Under composition of 1-cells, $\uBrPic(\B)$ is a monoidal groupoid. Since
$B(\M\boxtimes_{\B}\N)\simeq B(\M)\boxtimes_{\B}B(\N)=B(\N)\circ B(\M)$, the
assignment $\M\mapsto B(\M)$ is a monoidal functor
$\uPic(\B)\to\uBrPic(\B)^{\rev}$.
It is faithful, because an isomorphism $B(T)\cong B(T')$ of bimodule
functors is an isomorphism $T\cong T'$ of module functors. It is not full in
general: a bimodule equivalence $B(\M)\to B(\N)$ need not have equal left
and right constraints. On isomorphism classes, it identifies $\Pic(\B)$ with
the subgroup of one-sided classes in $\BrPic(\B)$.

\subsection{Braided autoequivalences fixing $\oB$}
\label{subsec:picard-center-groupoid}

The groupoid
$\uEqBr(\Z(\B)):=\uEqBr(\Z(\B),\Z(\B))$ is monoidal under composition of
1-cells. By Theorem~\ref{thm:BrPic} and Lemma~\ref{lem:pseudofunctor-Z-2},
the functor $\Z_{\B,\B}$ is a monoidal equivalence between them. Composing
it with $B$ gives a faithful functor $\M\mapsto\Phi_{B(\M)}$ from
$\uPic(\B)$ to $\uEqBr(\Z(\B))^{\rev}$, which is not full for the same
reason as $B$. We make it full by recording how $\Phi_{B(\M)}$ acts on
$G^-(\oB)$.

With $G^-$ as at the start of this section, we denote by
$\uAut^{\br}(\Z(\B);\oB)$ the
groupoid whose objects are pairs $(F,\eta)$, where $F$ is a braided tensor
autoequivalence of $\Z(\B)$ and $\eta:F\circ G^-\Rightarrow G^-$ is a
monoidal natural isomorphism, called a \emph{trivialization} of $F$. A
morphism $\alpha:(F,\eta)\to(F',\eta')$ is a monoidal natural isomorphism
$\alpha:F\Rightarrow F'$ such that $\eta'\circ(\alpha G^-)=\eta$.
Composition of functors gives the monoidal product
\[
  (F,\eta)\otimes(F',\eta')=\bigl(F\circ F',\ \eta\circ F(\eta')\bigr),
  \qquad \unit=(\id,\id_{G^-}),
\]
and forgetting the trivialization is a faithful monoidal functor
$\uAut^{\br}(\Z(\B);\oB)\to\uEqBr(\Z(\B))$. Its image on isomorphism classes
is the subgroup $\Aut_{\otimes}^{\br}(\Z(\B);\oB)$ of
$\Aut_{\otimes}^{\br}(\Z(\B))$ consisting of the classes $[F]$ with
$F\circ G^-\cong G^-$ as tensor functors.
Every object of $\uAut^{\br}(\Z(\B);\oB)$ is invertible: if $F^{-1}$ is a
quasi-inverse of $F$ with counit $\varepsilon:F^{-1}F\Rightarrow\id$, then
$\eta'=(\varepsilon G^-)\circ F^{-1}(\eta^{-1})$ is a trivialization of
$F^{-1}$. Since $\eta'\circ F^{-1}(\eta)=\varepsilon G^-$, the counit
$\varepsilon$ is an isomorphism $(F^{-1},\eta')\otimes(F,\eta)\to\unit$, and
since every object has a left inverse, every object is invertible.

The bridge between the two groupoids is the central action isomorphism
\eqref{eq:Phi-central-action}. For $\mathbf X=G^-(Y)$, a trivialization of
$\Phi_{\mathcal P}$ turns $\kappa^{\mathcal P}_{\mathbf X}$ into a
comparison between $Y\triangleright(-)$ and $(-)\triangleleft Y$, that is,
into the inverse of a one-sided structure.

\begin{lemma}\label{lem:Pic-trivialization}
  Let $\mathcal P$ and $\mathcal Q$ be invertible $\B$-bimodule categories.
  \begin{enumerate}
  \item The assignment
    $\eta\mapsto d^{-1}=(\id_{\mathcal P}\triangleleft\eta)\circ
    (\kappa^{\mathcal P}G^-)$ is a bijection from trivializations
    $\eta:\Phi_{\mathcal P}\circ G^-\Rightarrow G^-$ to one-sided structures
    $d$ on $\mathcal P$. For $\M\in\uPic(\B)$, we write $\eta_{\M}$ for the
    trivialization of $\Phi_{B(\M)}$ corresponding to the identity one-sided
    structure.
  \item Let trivializations $\eta$ of $\Phi_{\mathcal P}$ and $\eta'$ of
    $\Phi_{\mathcal Q}$ correspond to one-sided structures $d$ and $e$, and
    let $S:\mathcal P\to\mathcal Q$ be a bimodule equivalence. Then
    $\eta'\circ(\Z_{\B,\B}([S])\,G^-)=\eta$ if and only if
    \begin{equation}\label{eq:picard-structured-morphism}
      S^l_{Y,M}\circ S(d_{M,Y})=e_{S(M),Y}\circ S^r_{M,Y}
      \qquad(M\in\mathcal P,\ Y\in\B).
    \end{equation}
  \end{enumerate}
\end{lemma}
\begin{proof}
  (1) By Lemma~\ref{lem:braided-lift-unique}, applied to inverses, whose
  source $\Forg\circ G^-=\id_{\B}$ is surjective, a trivialization of
  $\Phi_{\mathcal P}$ is the same as a monoidal natural isomorphism
  $\Forg\circ\Phi_{\mathcal P}\circ G^-\Rightarrow\id_{\B}$ of tensor
  functors $\oB\to\B$. Let $\M$ be the underlying left module category of
  $\mathcal P$. Since $a_{\mathcal P}$ is fully faithful, composing
  with $\kappa^{\mathcal P}G^-$ identifies these with monoidal natural
  isomorphisms $\delta:\Forg\circ\Schau_{\mathcal P}\circ G^-\Rightarrow
  a_{\mathcal P}$ of tensor functors $\oB\to\B^*_{\mathcal P}$. The
  component $\delta_Y$ is a natural isomorphism
  $Y\triangleright(-)\Rightarrow(-)\triangleleft Y$ of left $\B$-module
  functors. Here $Y\triangleright(-)$ is a module functor through the
  half-braiding of $G^-(Y)$, and $(-)\triangleleft Y$ through the middle
  associator. For $d=\delta^{-1}$, the module condition on $\delta_Y$ says
  that $d^{-1}$ is compatible with the middle associators of $B(\M)$ and
  $\mathcal P$, since the half-braiding of $G^-(Y)$ gives the middle
  associator of $B(\M)$. Monoidality of $\delta$ says that $d^{-1}$ is
  compatible with the right associators, since the tensor structure of
  $\Forg\circ\Schau_{\mathcal P}\circ G^-$ is $c_{X,Y}^{-1}$ by
  \eqref{eq:Schauenburg-tensor-structure}, which is the right associator of
  $B(\M)$. Together these say that $d$ is a one-sided structure. For
  $\mathcal P=B(\M)$, the identity of $Y\triangleright M=M\triangleleft Y$
  is such a $\delta$, because $a_{B(\M)}=\Forg\circ\Schau_{\M}\circ G^-$ as
  tensor functors. It corresponds to the identity one-sided structure.

  (2) Put $u=\Z_{\B,\B}([S])$. Then
  $\zeta=\eta'\circ(uG^-)$ is a trivialization of
  $\Phi_{\mathcal P}$. Let $d^{\zeta}$ be the corresponding one-sided
  structure. Evaluate \eqref{eq:def-partial-2cell} at $\mathbf X=G^-(Y)$ and
  compose both sides with $S(M)\triangleleft\eta'_Y$. By (1)
  and the naturality of $S^r$, the result is
  \[
    S^r_{M,Y}\circ S\bigl((d^{\zeta}_{M,Y})^{-1}\bigr)
    =e^{-1}_{S(M),Y}\circ S^l_{Y,M}.
  \]
  Hence \eqref{eq:picard-structured-morphism} holds if and only if
  $d^{\zeta}=d$, that is, by (1), if and only if
  $\zeta=\eta$.
\end{proof}

We can now state the Picard version of Theorem~\ref{thm:BrPic}. For
$\M\in\uPic(\B)$, put $E(\M)=(\Phi_{B(\M)},\eta_{\M})$.

\begin{theorem}\label{thm:Pic}
  The assignment $\M\mapsto E(\M)$ extends to a monoidal equivalence
  \[
    E:\uPic(\B)\xrightarrow{\ \sim\ }\uAut^{\br}(\Z(\B);\oB)^{\rev},
    \qquad [T]\longmapsto\Z_{\B,\B}([B(T)]).
  \]
  In particular, $\uPic(\B)$ and $\uAut^{\br}(\Z(\B);\oB)$ are equivalent
  monoidal groupoids.
\end{theorem}
\begin{proof}
  A bimodule equivalence $S:B(\M)\to B(\N)$ satisfies
  \eqref{eq:picard-structured-morphism} for the identity one-sided
  structures exactly when $S^l=S^r$, that is, when $S=B(T)$ for its
  underlying module equivalence $T$. By Lemma~\ref{lem:Pic-trivialization}(2),
  $E$ is therefore a well-defined functor. It is fully faithful: by
  Theorem~\ref{thm:BrPic}, $\Z_{\B,\B}$ is a bijection from isomorphism
  classes of bimodule equivalences $B(\M)\to B(\N)$ to monoidal natural
  isomorphisms $\Phi_{B(\M)}\Rightarrow\Phi_{B(\N)}$. Those preserving
  $\eta_{\M}$ and $\eta_{\N}$ are the images of the $B(T)$, and $B$ is
  faithful.

  For essential surjectivity, let $(F,\eta)$ be an object of
  $\uAut^{\br}(\Z(\B);\oB)$. By Lemma~\ref{lem:BrPic-truncated}, there are
  an invertible $\B$-bimodule category $\mathcal P$ and a monoidal natural
  isomorphism $\alpha:F\Rightarrow\Phi_{\mathcal P}$. It is an isomorphism
  $(F,\eta)\to(\Phi_{\mathcal P},\eta\circ(\alpha^{-1}G^-))$. Let $d$ be the
  one-sided structure corresponding to $\eta\circ(\alpha^{-1}G^-)$, and let
  $\M$ be the underlying left module category of $\mathcal P$. The identity
  functor with right constraint $d^{-1}$ is a bimodule equivalence
  $S:B(\M)\to\mathcal P$, so $\M\in\uPic(\B)$. It satisfies
  \eqref{eq:picard-structured-morphism} for the identity structure on
  $B(\M)$ and $d$ on $\mathcal P$, so $\Z_{\B,\B}([S])$ gives an
  isomorphism $E(\M)\cong(F,\eta)$.

  For monoidality, let $\M,\N\in\uPic(\B)$. The compositor of
  Lemma~\ref{lem:pseudofunctor-Z-2}, followed by $\Z_{\B,\B}$ of the
  bimodule equivalence $B(\M)\boxtimes_{\B}B(\N)\simeq B(\M\boxtimes_{\B}\N)$
  of \S\ref{subsec:picard-groupoid}, is a monoidal natural isomorphism
  $\Phi_{B(\N)}\circ\Phi_{B(\M)}\Rightarrow\Phi_{B(\M\boxtimes_{\B}\N)}$. It
  is a morphism from $E(\N)\otimes E(\M)$ to $E(\M\boxtimes_{\B}\N)$: by
  \eqref{eq:compositor}, the compositor moves $\mathbf X=G^-(Y)$ through the
  two factors, and after the trivializations the two central action
  isomorphisms become identity one-sided structures by
  Lemma~\ref{lem:Pic-trivialization}(1).
  In detail, put $\mathcal P=B(\M)\boxtimes_{\B}B(\N)$, let
  $\gamma:\Phi_{B(\N)}\circ\Phi_{B(\M)}\Rightarrow\Phi_{\mathcal P}$ be the
  compositor, and let
  $\eta^{\mathcal P}=\eta_{\N}\circ\Phi_{B(\N)}(\eta_{\M})\circ(\gamma G^-)^{-1}$.
  By \eqref{eq:compositor}, the naturality of $\kappa^{B(\N)}$ and of the
  balancing $b$, and Lemma~\ref{lem:Pic-trivialization}(1) for $\eta_{\M}$
  and $\eta_{\N}$,
  \[
    (\id_{\mathcal P}\triangleleft\eta^{\mathcal P}_Y)\circ\kappa^{\mathcal P}_{G^-(Y),M\boxtimes N}
    =b_{M,Y,N}
    \qquad(M\in\M,\ N\in\N,\ Y\in\B),
  \]
  so $\eta^{\mathcal P}$ corresponds to the one-sided structure $b^{-1}$. The
  equivalence $S:\mathcal P\to B(\M\boxtimes_{\B}\N)$ is the identity on
  objects $M\boxtimes N$, with left constraint the identity and right
  constraint $b^{-1}$, so \eqref{eq:picard-structured-morphism} holds for $S$.
  By Lemma~\ref{lem:Pic-trivialization}(2), $\Z_{\B,\B}([S])\circ\gamma$ is
  compatible with the trivializations. For the unit, \eqref{eq:unitor} shows
  that $(\upsilon_{\B}G^-)^{-1}$ corresponds to the one-sided structure
  $c_{V,Y}$ on the regular bimodule category, and the equivalence
  $B(\B)\to\B$ has right constraint $c_{V,Y}^{-1}$, so
  Lemma~\ref{lem:Pic-trivialization}(2) applies again. The associator and
  unitors of $\boxtimes_{\B}$, induced by rebracketing and by the action
  functors, are compatible with these one-sided structures by the balancing
  axiom and the hexagon axioms. The
  coherence axioms hold in $\uEqBr(\Z(\B))$ by
  Lemma~\ref{lem:pseudofunctor-Z-2}, and hence in $\uAut^{\br}(\Z(\B);\oB)$,
  since forgetting trivializations is faithful.

  Finally, since every object of $\uAut^{\br}(\Z(\B);\oB)$ is invertible,
  taking inverses is a monoidal equivalence from
  $\uAut^{\br}(\Z(\B);\oB)^{\rev}$ to $\uAut^{\br}(\Z(\B);\oB)$.
\end{proof}

We now pass to isomorphism classes.

\begin{corollary}{(\textit{c.f.} \cite[Theorem~4.3]{davydov2013picard})}
  \label{cor:Pic-center-autoequivalences}
  The assignment $[\M]\mapsto[\Phi_{B(\M)}^{-1}]$ is a group isomorphism
  \[
    \Pic(\B)\xrightarrow{\ \sim\ }\Aut_{\otimes}^{\br}(\Z(\B);\oB).
  \]
\end{corollary}
\begin{proof}
  By Theorem~\ref{thm:Pic}, $[\M]\mapsto[E(\M)]^{-1}$ is a group isomorphism
  from $\Pic(\B)$ onto the group of isomorphism classes of
  $\uAut^{\br}(\Z(\B);\oB)$. It remains to show that forgetting
  trivializations is injective on isomorphism classes, that is, that two
  trivializations $\eta,\eta'$ of the same $F$ give isomorphic objects.

  The monoidal automorphism $\eta'\circ\eta^{-1}$ of $G^-$ is $G^-(z)$ for a
  monoidal automorphism $z$ of $\id_{\B}$, since $G^-$ is fully faithful. We
  extend $z$ to the center by $\widetilde z_{(V,\sigma)}=z_V$. To see that
  $z_V$ is a morphism in $\Z(\B)$, let $W\in\B$. Naturality and monoidality
  of $z$ give $(z_W\otimes z_V)\sigma_W=\sigma_W(z_V\otimes z_W)$, and
  naturality of $\sigma$ in $W$ lets us cancel $z_W$:
  \[
    (\id_W\otimes z_V)\sigma_W=\sigma_W(z_V\otimes\id_W).
  \]
  So $\widetilde z$ is a monoidal automorphism of $\id_{\Z(\B)}$. Put
  $\alpha=\widetilde zF$. Naturality of $\widetilde z$ gives
  $\eta\circ(\alpha G^-)=(\widetilde zG^-)\circ\eta=\eta'$, so $\alpha$ is
  an isomorphism $(F,\eta')\to(F,\eta)$.
\end{proof}

\subsection{The boundary map and pullback action}
\label{subsec:picard-boundary}

The centralizer of $G^-(\oB)$ in $\Z(\B)$ is the essential image of $G^+$.
Indeed, for $(V,\sigma)\in\Z(\B)$ and $X\in\B$, the double braiding of
$(V,\sigma)$ with $G^-(X)$ is $\sigma_X\circ c_{V,X}^{-1}$, which is the
identity for all $X$ exactly when $\sigma=c_{V,-}$. So a braided
autoequivalence of $\Z(\B)$ that fixes $G^-(\oB)$ also preserves $G^+(\B)$.
On Picard classes this restriction is the boundary map of Davydov and
Nikshych, which we now recall in our conventions.

For a left $\B$-module category $\M$, the \emph{alpha-induction functors}
\[
  \alpha^{+}_{\M}=\Forg\circ\Schau_{\M}\circ G^+,
  \qquad
  \alpha^{-}_{\M}=\Forg\circ\Schau_{\M}\circ G^-
\]
are tensor functors $\B\to\B^*_{\M}$. Recall that $\oB$ and $\B$ have the
same underlying tensor category. Explicitly, $\alpha^{\pm}_{\M}(X)$ is
$X\triangleright(-)$, with module structure
$Y\triangleright(X\triangleright M)\to X\triangleright(Y\triangleright M)$ given
by $c_{X,Y}^{-1}$ and by $c_{Y,X}$, respectively, and with tensor structure
$(Y\otimes X)\triangleright(-)\to(X\otimes Y)\triangleright(-)$ given by
$c_{Y,X}$ and by $c_{X,Y}^{-1}$, respectively, using
\eqref{eq:Schauenburg-tensor-structure}. The underlying module functors are
those of \cite[\S2I, equation~(22)]{davydov2013picard}.

The group $\Aut_{\otimes}(\B)$ of isomorphism classes of tensor
autoequivalences acts on the set $\bbX_\B$ of
\S\ref{subsec:BrPic-applications} by pullback. For
$F\in\Aut_{\otimes}(\B)$ and $\M\in\bbX_\B$, write $F\cdot\M={}_F\M$, where
${}_F\M$ has the same underlying category as $\M$ and the action
\[
  X \tr^F M = F^{-1}(X) \tr M;
\]
it is again indecomposable exact, since $F$ is an equivalence. Isomorphic
functors give equivalent module categories, so this is a well-defined
action. Two elements of $\bbX_\B$ are \emph{pullback equivalent} if they lie
in the same orbit.

\begin{lemma}\label{lem:boundary-map}
  Let $\M\in\uPic(\B)$.
  \begin{enumerate}
  \item There is a braided tensor autoequivalence $\partial_{\M}$ of $\B$,
    unique up to isomorphism, with
    $\Phi_{B(\M)}\circ G^+\cong G^+\circ\partial_{\M}$ as braided tensor
    functors. Moreover, $\alpha^-_{\M}$ is an equivalence, and
    $\partial_{\M}$ is, up to isomorphism, the unique tensor functor with
    $\alpha^+_{\M}\cong\alpha^-_{\M}\circ\partial_{\M}$.
  \item For $\M,\N\in\uPic(\B)$, there is an isomorphism
    $\partial_{\M\boxtimes_{\B}\N}\cong\partial_{\N}\circ\partial_{\M}$.
  \item For $F\in\Aut_{\otimes}^{\br}(\B)$, the pullback ${}_F\M$ lies in
    $\uPic(\B)$, and
    \begin{equation}\label{eq:picard-boundary-equivariance}
      \partial_{{}_F\M}\cong F\circ\partial_{\M}\circ F^{-1}.
    \end{equation}
  \end{enumerate}
\end{lemma}

\begin{proof}
  (1) Put $\mathcal P=B(\M)$. The trivialization $\eta_{\M}$ of
  $\Phi_{\mathcal P}$ fixes $G^-(\oB)$, so the braided autoequivalence
  $\Phi_{\mathcal P}$ preserves its centralizer $G^+(\B)$. As $G^+$ is fully
  faithful, $\Phi_{\mathcal P}\circ G^+\cong G^+\circ\partial_{\M}$ for a
  braided autoequivalence $\partial_{\M}$, unique up to isomorphism. By
  \S\ref{subsec:picard-groupoid}, $\alpha^-_{\M}=a_{\mathcal P}$ as
  tensor functors, and $a_{\mathcal P}$ is an equivalence by
  Proposition~\ref{prop:invertibility-criterion}. Restricting
  $\kappa^{\mathcal P}$ along $G^+$ gives
  \[
    \alpha^+_{\M}=\Forg\circ\Schau_{\mathcal P}\circ G^+
    \cong a_{\mathcal P}\circ\Forg\circ\Phi_{\mathcal P}\circ G^+
    \cong a_{\mathcal P}\circ\partial_{\M}
    \cong\alpha^-_{\M}\circ\partial_{\M}.
  \]
  Conversely, this isomorphism determines $\partial_{\M}$ up to isomorphism,
  since $\alpha^-_{\M}$ is an equivalence.

  (2) The bimodule equivalence
  $B(\M)\boxtimes_{\B}B(\N)\simeq B(\M\boxtimes_{\B}\N)$ of
  \S\ref{subsec:picard-groupoid} and
  Lemmas~\ref{lem:pseudofunctor-Z-1} and~\ref{lem:pseudofunctor-Z-2} give
  $\Phi_{B(\M\boxtimes_{\B}\N)}\cong\Phi_{B(\N)}\circ\Phi_{B(\M)}$, because
  $B(\M)\boxtimes_{\B}B(\N)=B(\N)\circ B(\M)$. Restrict to $G^+(\B)$ and
  use (1).

  (3) A left $\B$-module functor for the pulled-back action is the same as a
  module functor for the original action, with the structure transported
  along $F$. This identifies $\B^*_{{}_F\M}$ with $\B^*_{\M}$. Since $F$ is
  braided, it carries $c_{X,Y}$ and $c_{Y,X}^{-1}$ to the corresponding
  braidings, so this identification turns $\alpha^{\pm}_{{}_F\M}$ into
  $\alpha^{\pm}_{\M}\circ F^{-1}$. In particular $\alpha^-_{{}_F\M}$, which is
  the action functor of $B({}_F\M)$ by \S\ref{subsec:picard-groupoid},
  is an equivalence, so $B({}_F\M)$ is invertible by
  Proposition~\ref{prop:invertibility-criterion}, and
  \[
    \alpha^-_{{}_F\M}\circ F\partial_{\M}F^{-1}
    \cong\alpha^-_{\M}\circ\partial_{\M}\circ F^{-1}
    \cong\alpha^+_{\M}\circ F^{-1}
    \cong\alpha^+_{{}_F\M}.
  \]
  By the uniqueness in (1), $\partial_{{}_F\M}\cong F\partial_{\M}F^{-1}$.
\end{proof}

By Lemma~\ref{lem:boundary-map}(1), $\partial_{\M}$ is the analogue of the functor defined by \cite[equation~(23)]{davydov2013picard}, and by
Lemma~\ref{lem:pseudofunctor-Z-1} it depends only on the class of $\M$. The map
\[
  \partial:\Pic(\B)\to\Aut_{\otimes}^{\br}(\B),
  \qquad [\M]\mapsto[\partial_{\M}],
\]
is called the \emph{boundary map}. By Lemma~\ref{lem:boundary-map}(2), it reverses the
order of products, so $[\M]\mapsto[\partial_{\M}]^{-1}$ is a group
homomorphism. By Lemma~\ref{lem:boundary-map}(3), the subset
$\Pic(\B)\subset\bbX_\B$ of \S\ref{subsec:picard-groupoid} is stable
under the subgroup $\Aut_{\otimes}^{\br}(\B)$, and $\partial$ turns this
action into conjugation. Equation~\eqref{eq:picard-boundary-equivariance} is the
analogue, in our conventions, of the crossed-module compatibility of
\cite[Theorem~3.10, equation~(36)]{davydov2013picard}.

\begin{proposition}\label{prop:boundary-bijective}
  If $\B$ is factorizable, then $\partial$ is bijective.
\end{proposition}
This is \cite[Corollary~4.6]{davydov2013picard}, whose proof deduces it from
\cite[Theorem~4.3]{davydov2013picard}. We use Corollary~\ref{cor:Pic-center-autoequivalences} instead (compare
Remark~\ref{rem:DN-gap}).
\begin{proof}
  By Lemma~\ref{lem:boundary-map}(1), $[\partial_{\M}]^{-1}$ is the restriction to
  $G^+(\B)$ of $[\Phi_{B(\M)}^{-1}]$, and $[\M]\mapsto[\Phi_{B(\M)}^{-1}]$ is
  bijective onto $\Aut_{\otimes}^{\br}(\Z(\B);\oB)$ by
  Corollary~\ref{cor:Pic-center-autoequivalences}. So it suffices to show that restriction to $G^+(\B)$ is a bijection
  $\Aut_{\otimes}^{\br}(\Z(\B);\oB)\to\Aut_{\otimes}^{\br}(\B)$ when
  $G:\B\boxtimes\oB\to\Z(\B)$ is an equivalence. Its inverse sends $a$ to
  $[G(a\boxtimes\id)G^{-1}]$, which lies in $\Aut_{\otimes}^{\br}(\Z(\B);\oB)$
  because $G(a\boxtimes\id)G^{-1}G^-\cong G^-$, and restricts to $a$ because
  $G(a\boxtimes\id)G^{-1}G^+\cong G^+a$. The restriction map is injective because
  an autoequivalence $F$ with $FG^+\cong G^+$ and $FG^-\cong G^-$ satisfies
  $FG\cong G$: the two isomorphisms combine to a monoidal one, since $F$
  preserves the braiding between the two factors, by the universal property of the Deligne product $\B\boxtimes\oB$.
\end{proof}

We are interested in the number of orbits of $\bbX_{\B}$ by the action of $\Aut_{\otimes}^{\mathrm{br}}(\B)$, especially whether it is finite or not.
When the boundary map is surjective, its equivariance gives the following
criterion for infinitely many pullback orbits.

\begin{proposition}\label{prop:inf-orbit-criteria}
Suppose that $\partial$ is surjective. If the number of conjugacy classes in
$\Aut_{\otimes}^{\mathrm{br}}(\B)$ is infinite, then
$\Aut_{\otimes}^{\mathrm{br}}(\B)$ has infinitely many orbits in $\Pic(\B)$,
and hence in $\bbX_\B$.
\end{proposition}
\begin{proof}
Set $P = \Aut_{\otimes}^{\mathrm{br}}(\B)$. The subset
$\Pic(\B)\subset \bbX_\B$ is stable under $P$, so $\Pic(\B)/P$ is a subset
of $\bbX_\B/P$. By \eqref{eq:picard-boundary-equivariance}, the surjective
map $\partial$ sends $P$-orbits to conjugacy classes, so it induces a
surjection
\[ \Pic(\B)/ P \longrightarrow \{\text{conjugacy classes in } P\} .\]
If $P$ has infinitely many conjugacy classes, then $\Pic(\B)/P$ is
infinite, and hence so is $\bbX_\B/P$.
\end{proof}

\begin{corollary}\label{cor:factorizable-inf-orbits}
  Suppose that $\B$ is factorizable. If the number
  of conjugacy classes in $\Aut_{\otimes}^{\mathrm{br}}(\B)$ is infinite,
  then $\Aut_{\otimes}^{\mathrm{br}}(\B)$ has infinitely many orbits in
  $\Pic(\B)$, and hence in $\bbX_\B$.
\end{corollary}
\begin{proof}
  By Proposition~\ref{prop:boundary-bijective}, $\partial$ is bijective.
  Hence Proposition~\ref{prop:inf-orbit-criteria} applies.
\end{proof}

\subsection{Small quantum groups}
\label{sec:uqsl2-picard}
In this subsection we apply the results above to the category of modules
over a small quantum group, and show that it has infinitely many
pullback-equivalence classes of indecomposable exact left module categories,
all inside the MME class of the regular module category
(Proposition~\ref{prop:uqg-regular-orbit} and
Theorem~\ref{thm:uqg-Aut-orbits-infinite}).
Let $\kk$ be an algebraically closed field of characteristic $0$.
Let $\fg$ be a finite-dimensional simple Lie algebra, and let $G^{\ad}$ be
the adjoint simple algebraic group with Lie algebra $\fg$. Let $q\in\kk$ be
a primitive $l$-th root of unity for an odd integer $l>1$, and assume
\[ \gcd(l,\det A_\fg)=1,\]
where $A_\fg$ denotes the Cartan matrix. In type $G_2$, we assume that
$\gcd(l,3)=1$.

Let $u_q(\fg)$ denote the associated small quantum group of Lusztig
\cite{lusztig1990finite}, and let
\[ \B = \Rep(u_q(\fg))\]
be its category of finite-dimensional modules, with the braiding of
\cite[\S2.3]{davydov2018autoequivalences}, for the rest of this
subsection. Put
\[
  P=\Aut_{\otimes}^{\mathrm{br}}(\B).
\]
We use three facts from \cite{davydov2018autoequivalences}: $\B$ is
factorizable, so that
$G:\B\boxtimes\oB\to\Z(\B)$ is a braided equivalence; every tensor
autoequivalence of $\B$ is braided
\cite[Theorem~5.4]{davydov2018autoequivalences}, so that
$P=\Aut_{\otimes}(\B)$; and every braided autoequivalence of $\Z(\B)$
preserves the two factors $G^+(\B)$ and $G^-(\oB)$, up to isomorphism, by
the argument in the proof of
\cite[Proposition~5.9(ii)]{davydov2018autoequivalences}, which uses only
Lemma~2.5 and Proposition~5.3 there.
Since reversing the braiding does not change braided autoequivalences, the
restrictions to the two factors give a group isomorphism
\begin{equation}\label{eq:uqg-center-autos}
  \Aut_{\otimes}^{\br}(\Z(\B))\xrightarrow{\ \sim\ }P\times P,
  \qquad F\longmapsto(F|_{G^-(\oB)},\,F|_{G^+(\B)}),
\end{equation}
with inverse $(a,b)\mapsto G\circ(b\boxtimes a)\circ G^{-1}$.

Two kinds of invertible bimodule categories are easy to place in these
coordinates. For a tensor autoequivalence $F$ of $\B$, let $\Gamma_F$ be $\B$ with the regular left action and the right
action $Y\triangleleft X=Y\otimes F(X)$. Its action functor is $F$
followed by the equivalence $\B\simeq\B^*_{\B}$ given by right
multiplication, so $\Gamma_F$ is invertible by Proposition~\ref{prop:invertibility-criterion}.

\begin{lemma}\label{lem:graph-bimodules}
  Let $F$ be a braided tensor autoequivalence of $\B$, and let
  $\M\in\uPic(\B)$.
  \begin{enumerate}
  \item $\Phi_{\Gamma_F}\cong\widetilde{F^{-1}}$, so the class of
    $\Phi_{\Gamma_F}^{-1}$ has coordinates $(F,F)$ in
    \eqref{eq:uqg-center-autos}.
  \item The class of $\Phi_{B(\M)}^{-1}$ has coordinates
    $(1,\partial_{\M}^{-1})$.
  \item For every left $\B$-module category $\N$, there is an equivalence
    $\Gamma_F\boxtimes_{\B}\N\simeq{}_F\N$ of left $\B$-module categories.
  \end{enumerate}
\end{lemma}
\begin{proof}
  (1) For $\mathbf X=(X,\sigma)\in\Z(\B)$ and $Y\in\B$, the half-braiding
  $\sigma_Y:X\otimes Y\to Y\otimes X$, followed by the counit
  $X\cong FF^{-1}(X)$, is a natural isomorphism
  $\mathbf X\triangleright Y\cong Y\triangleleft\widetilde{F^{-1}}(\mathbf X)$.
  It is left $\B$-linear by the half-braiding axiom and monoidal in
  $\mathbf X$, so Lemma~\ref{lem:central-action-unique} gives
  $\Phi_{\Gamma_F}\cong\widetilde{F^{-1}}$, so $\Phi_{\Gamma_F}^{-1}\cong\widetilde F$. Since $F$ is braided,
  $\widetilde F$ maps $G^{\pm}(X)$ to $G^{\pm}(F(X))$, so it restricts to $F$
  on both factors.

  (2) The trivialization $\eta_{\M}$ makes the restriction of this class to
  $G^-(\oB)$ trivial, and Lemma~\ref{lem:boundary-map}(1) gives the
  restriction $\partial_{\M}^{-1}$ to $G^+(\B)$.

  (3) The functor $\Gamma_F\times\N\to\N$,
  $(Y,N)\mapsto F^{-1}(Y)\triangleright N$, is $\B$-balanced, since
  $F^{-1}(Y\otimes F(X))\triangleright N\cong
  F^{-1}(Y)\triangleright(X\triangleright N)$. The balanced functor induces an equivalence
  $\Gamma_F\boxtimes_{\B}\N\simeq\N$, whose inverse is
  $N\mapsto\unit\boxtimes_{\B}N$. This equivalence carries the left action
  $Z\triangleright(Y\boxtimes_{\B}N)=(Z\otimes Y)\boxtimes_{\B}N$ to
  $F^{-1}(Z)\triangleright(-)$, which is the action of ${}_F\N$.
\end{proof}

The regular orbit explains how infinitely many pullback classes can lie
in a single MME class.

\begin{proposition}\label{prop:uqg-regular-orbit}
There are group isomorphisms
\[
  \Pic(\B)\xrightarrow{\ \sim\ }P,\ [\M]\mapsto[\partial_{\M}]^{-1},\qquad
  \BrPic(\B)\cong P\times P.
\]
Under the last isomorphism, $\Pic(\B)$ is the subgroup $\{1\}\times P$. The
$\BrPic(\B)$-orbit of the regular module category is precisely $\Pic(\B)$,
with stabilizer $\Delta P$. As a homogeneous space it is
\[
  (P\times P)/\Delta P\cong P,
  \qquad (a,b)\Delta P\longmapsto ba^{-1},
\]
and in these coordinates the action is $(a,b)\cdot x=bxa^{-1}$. The
coordinate of $\M\in\Pic(\B)$ is $\partial_{\M}^{-1}$. Thus the Picard locus
is one MME class, and its pullback classes are the conjugacy classes of $P$.
\end{proposition}
The isomorphism $\Pic(\B)\cong P$ matches
\cite[Theorem~5.4]{davydov2018autoequivalences}, and
$\BrPic(\B)\cong P\times P$ is
\cite[Proposition~5.9(ii)]{davydov2018autoequivalences} without the
additional assumptions of \cite[Theorem~5.6]{davydov2018autoequivalences}.

\begin{proof}
The map $[\M]\mapsto[\partial_{\M}]^{-1}$ is bijective by
Proposition~\ref{prop:boundary-bijective}, since $\B$ is factorizable, and it
is a group homomorphism by Lemma~\ref{lem:boundary-map}(2). Composing the
group isomorphism
$[\mathcal P]\mapsto[\Phi_{\mathcal P}^{-1}]$ of
\S\ref{sec:BrPic-2-groupoid} with \eqref{eq:uqg-center-autos}
gives $\BrPic(\B)\cong P\times P$. By Lemma~\ref{lem:graph-bimodules}(2) and
the bijectivity of $\partial$, this isomorphism maps $\Pic(\B)$ onto $\{1\}\times P$.

The stabilizer of the regular module category $\B$ consists of the classes
of invertible bimodule categories $\mathcal P$ with $\mathcal P\simeq\B$ as
left module categories. Such a $\mathcal P$ is equivalent to $\B$ with a right
action by a tensor autoequivalence, by Lemma~\ref{lem:bimodule-criterion}:
transporting the action functor of $\mathcal P$ along a left module equivalence
$\B\simeq\mathcal P$ gives a tensor functor $\B\to\B^*_{\B}\simeq\B$;
that is, $\mathcal P\simeq\Gamma_F$. By Lemma~\ref{lem:graph-bimodules}(1),
the stabilizer is $\Delta P$. The factorization $(a,b)=(1,ba^{-1})(a,a)$
shows that every coset contains exactly one element of $\{1\}\times P$, so
$[\M]\mapsto[B(\M)\boxtimes_{\B}\B]=[\M]$ identifies $\Pic(\B)$ with the
orbit, with coordinate $\partial_{\M}^{-1}$ by
Lemma~\ref{lem:graph-bimodules}(2). Multiplying cosets gives the action
formula. By Lemma~\ref{lem:graph-bimodules}(3), the diagonal element
$(F,F)$ acts by pullback along $F$, which is conjugation by $F$, in
agreement with \eqref{eq:picard-boundary-equivariance}. Finally,
Theorem~\ref{thm:MME} identifies $\BrPic(\B)$-orbits with MME classes.
\end{proof}

The same factorization constructs representatives for all Brauer--Picard
classes.

\begin{corollary}\label{cor:uqg-brpic-representatives}
For $p_1,p_2\in P$, choose $\M\in\Pic(\B)$ with $\partial_{\M}=p_1p_2^{-1}$
and a representative $F$ of $p_1$, and let $\mathcal P_{p_1,p_2}$ be $B(\M)$
with the right action twisted by $F$: $M\triangleleft X=F(X)\triangleright M$. Then $\mathcal P_{p_1,p_2}$ is invertible with coordinates $(p_1,p_2)$. Its
class is independent of the choices, and
\[
  \mathcal P_{p_1,p_2}\boxtimes_{\B}\mathcal P_{q_1,q_2}
  \simeq\mathcal P_{p_1q_1,p_2q_2}.
\]
\end{corollary}
\begin{proof}
The unit equivalence $B(\M)\boxtimes_{\B}\B\simeq\M$ identifies
$B(\M)\boxtimes_{\B}\Gamma_F$ with $\mathcal P_{p_1,p_2}$. Its coordinates
are $(1,\partial_{\M}^{-1})(p_1,p_1)=(p_1,p_2)$ by
Lemma~\ref{lem:graph-bimodules}. Since
$\BrPic(\B)\cong P\times P$ is a group isomorphism, the class is determined
by its coordinates, which gives independence of the choices and the
product formula.
\end{proof}

The identity component of $P$ is $G^{\ad}$
\cite[Proposition~5.2]{davydov2018autoequivalences}, and
$P=\Gamma\ltimes(G^{\ad}\times Z)$ for the Dynkin diagram automorphism group
$\Gamma$ and a finite group $Z$
\cite[Lemma~5.5]{davydov2018autoequivalences}.
Under the additional order assumptions of
\cite[Theorem~5.6]{davydov2018autoequivalences}, one has
$P\cong\Gamma\ltimes G^{\ad}$. The results above use only the assumptions
stated at the beginning of this subsection.

\begin{theorem}\label{thm:uqg-Aut-orbits-infinite}
  The category $\Rep(u_q(\fg))$ has infinitely many pullback-equivalence
  classes of indecomposable exact left module categories.
\end{theorem}

\begin{proof}
Since $\B$ is factorizable and $P=\Aut_{\otimes}(\B)$, it suffices by
Corollary~\ref{cor:factorizable-inf-orbits} to show that $P$ has infinitely
many conjugacy classes. By Proposition~\ref{prop:uqg-regular-orbit}, the
resulting pullback classes all lie in the MME class of the regular module
category.

Put $H=G^{\ad}(\kk)$. Since $P=\Gamma\ltimes(G^{\ad}\times Z)$ with $\Gamma$
and $Z$ finite \cite[Lemma~5.5]{davydov2018autoequivalences}, $H$ is normal
of finite index in $P$.
A $P$-conjugacy class in $H$ is a union of at most $[P:H]$
$H$-conjugacy classes: conjugation by representatives of $P/H$ lists
all possible classes. Hence infinitely many $H$-conjugacy classes
imply infinitely many $P$-conjugacy classes.

Choose a maximal torus
\[
  T\subset G^{\ad},
  \qquad
  W=N_{G^{\ad}}(T)/T.
\]
Every semisimple element of $H$ is conjugate to an element of $T(\kk)$,
and two elements of $T(\kk)$ are conjugate in $H$ if and only if they
belong to the same $W$-orbit. Hence the semisimple conjugacy classes in
$H$ are parametrized by $T(\kk)/W$.

Since $\kk$ is algebraically closed, $T(\kk)\cong(\kk^{\times})^{\rank\fg}$
is infinite, whereas $W$ is finite. Thus $T(\kk)/W$ is infinite.
Therefore $H$, and hence $P$, has infinitely many conjugacy classes.
\end{proof}

\bibliographystyle{alpha}
\bibliography{references}

\end{document}